\documentclass{article}

\usepackage[scaled=0.92]{helvet}
\usepackage{amssymb,latexsym,amsmath,amsxtra}
\usepackage{amsfonts}
\usepackage{xspace}
\usepackage[final]{graphics}
\usepackage{graphicx,color,psfrag} 
\usepackage{subfig}
\usepackage{array}
\usepackage{enumitem}

\usepackage{datetime}

\usepackage{latexsym}

\usepackage{xspace}
\usepackage[colorlinks,bookmarks,urlcolor=blue,citecolor=blue,linkcolor=blue]{hyperref}

\renewenvironment{thebibliography}[1]{%
\begin{oldthebibliography}{#1}%
\setlength{\parskip}{0pt}%
\setlength{\itemsep}{0.5pt}%
}%
{%
\end{oldthebibliography}%
}

\renewcommand{\vec}[1]{\boldsymbol{#1}}
\newcommand{\vecs}[1]{\boldsymbol{#1}}
\newcommand{\paren}[1]{\left(#1\right)}
\newcommand{\brac}[1]{\left[#1\right]}
\newcommand{\E}{\mathbb{E}}

\newcommand{\avg}[1]{\E\brac{#1}}

\newcommand{\abs}[1]{\left|#1\right|}

\newcommand{\veta}{\vecs{\eta}}

\newcommand{\vx}{\vec{x}}

\newcommand{\vy}{\vec{y}}

\newcommand{\vQ}{\vec{Q}}

\def\R{\mathbb{R}}

\newcommand{\comment}[1]{{\color{black}{#1}}}

\definecolor{issuePJA_color}{rgb}{1.0,0.0,0.0}

\definecolor{commentPJA_color}{rgb}{1.0,0.0,0.8}

\definecolor{issueSI_color}{rgb}{0.0,0.7,0.0}

\definecolor{commentSI_color}{rgb}{0.1,0.1,0.6}

\newcommand{\mb}[1]{{\mathbf{#1}}}

\title{A First-Passage Kinetic Monte Carlo Method for Reaction-Drift-Diffusion Processes}

\author{ Ava J. Mauro \thanks{
    Department of Mathematics and Statistics, Boston University, 111
    Cummington St., Boston, MA 02215 (avamauro@bu.edu)}  
    \and 
    Jon Karl Sigurdsson \thanks{Department of Mathematics, University of California, Santa Barbara} 
    \and 
    Justin Shrake \thanks{Department of Mathematics, University of California, Santa Barbara}
    \and 
    Paul J. Atzberger \thanks{6712 South Hall, Department of
    Mathematics, University of California, Santa Barbara, CA 93106
    (atzberg@math.ucsb.edu)} 
    \and
    Samuel A. Isaacson \thanks{
    Department of Mathematics and Statistics, Boston University, 111
    Cummington St., Boston, MA 02215 (isaacson@math.bu.edu)}
}

\date{}

\begin{document}

\maketitle
\begin{abstract}
   Stochastic reaction-diffusion models are now a popular tool for
  studying physical systems in which both the explicit diffusion of
  molecules and noise in the chemical reaction process play important
  roles. The Smoluchowski diffusion-limited reaction model (SDLR) is
  one of several that have been used to study biological systems. 
  Exact realizations of the underlying stochastic process described by
  the SDLR model can be generated by the recently proposed
  First-Passage Kinetic Monte Carlo (FPKMC) method.  This exactness
  relies on sampling analytical solutions to one and two-body
  diffusion equations in simplified protective domains.

  In this work we extend the FPKMC to allow for drift arising from
  fixed, background potentials. As the corresponding Fokker-Planck
  equations that describe the motion of each molecule can no longer be
  solved analytically, we develop a hybrid method that discretizes the
  protective domains. The discretization is chosen so that the
  drift-diffusion of each molecule within its protective domain is
  approximated by a continuous-time random walk on a lattice.  New
  lattices are defined dynamically as the protective domains are
  updated, \comment{hence we will refer to our method as Dynamic
  Lattice FPKMC or DL-FPKMC. We focus primarily on the 
  one-dimensional case in this manuscript, and demonstrate} the 
  numerical convergence and accuracy of our method in
  \comment{this case} for both smooth and discontinuous potentials.
  \comment{We also present applications of our method, which 
  illustrate the impact of  drift on reaction kinetics.}  
\end{abstract} 

\thispagestyle{empty} 

\section{Introduction}

A fundamental challenge in cell biology is to understand how to
predict and control the dynamics of cellular
processes~\cite{AlbertsMOLECELLBIO}.  Stochasticity in the quantities
and movements of molecules can have significant effects on the
outcomes of cellular processes, particularly given the low copy numbers
of many signaling and regulatory proteins and mRNAs present in a cell.
For such species, the actual number and locations of molecules can
provide a more accurate and useful description than the local
concentration. The method we present in this paper allows for the
explicit simulation of the stochastically varying numbers and
locations of molecular species undergoing chemical reactions and
drift-diffusion.

Experimental studies and mathematical models have shown that
stochasticity in the chemical reaction process plays a role in many
cellular processes, for example gene
expression~\cite{CollinsNatureEukNoise,ElowitzScienceNoiseExper,RaserScience2004},
cell-fate decision
making~\cite{ElowitzNatureCompentence06,RajComptence07}, and signaling
pathways in development~\cite{AriasNoiseDevelop06}.  Mathematical
models of such processes frequently treat an individual cell as a
single well-mixed volume or a small number of well-mixed compartments,
such as the nucleus and the cytoplasm. However, the heterogeneous
spatial distribution of chemical species, as well as interactions with
internal membranes and organelles, often have significant effects on
cellular processes.  For example, after gene regulatory proteins enter
the nucleus through nuclear pores, the time required to find specific
DNA binding sites can be significantly influenced by the spatial
structure of DNA within the nucleus~\cite{IsaacsonPNAS2011}.
Similarly, the spatial distribution of components of cellular
signaling processes can play a decisive role in the successful
propagation of signals from the cell membrane to the
nucleus~\cite{MunozGarcia:2009hd, Kholodenko:2010jv}.

A number of stochastic reaction-diffusion models have been introduced
to understand the combined influence of noise due to the chemical
reaction process and spatial
diffusion~\cite{AndrewsMinCDE04,ElfIEESys04,WoldeEgfrdPNAS2010,IsaacsonPNAS2011}.
These mathematical models resolve the explicit spatial movement of
proteins and mRNAs within cells.  There are three such mathematical
models that have been commonly used: the spatially-continuous
Smoluchowski diffusion-limited reaction model
(SDLR)~\cite{SmoluchowskiDiffLimRx, KeizerJPhysChem82}; what we call
the Doi
model~\cite{TeramotoDoiModel1967,DoiSecondQuantA,DoiSecondQuantB}; and
the lattice-based reaction-diffusion master equation model
(RDME)~\cite{GardinerRXDIFFME,McQuarrieJAppProb}. 

These models often treat the movement of molecules as purely
diffusive; however, drift can also play a significant role in the
dynamics of cellular processes. Examples of sources of such drift
include active transport, variations in chemical potential, material
heterogeneities in the cytoplasm and nucleoplasm, and local
interactions with cellular structures. The incorporation of drift has
played a key role in developing models for molecular-motor based
active transport~\cite{PeskinOsterMotor1995, Atzberger2006}, movement
of proteins on DNA~\cite{MirnyDNASlide04}, and protein movement
subject to the influence of volume exclusion by
chromatin~\cite{IsaacsonPNAS2011}.  It has been shown how to extend
the RDME to incorporate drift due to
potentials~\cite{IsaacsonPNAS2011}.  In the present work we consider a
generalization of the SDLR model that allows for drift due to fixed
potentials.  For creating realizations of the stochastic process
described by this model, we present a new numerical method, combining
elements of the First-Passage Kinetic Monte Carlo (FPKMC)
\comment{method~\cite{KalosDSMC06,OppelstrupPRE2009,DonevJCP2010,WoldeEgfrdPNAS2010}}
and the lattice methods
of~\cite{ElstonPeskinJTB2003,IsaacsonPNAS2011}. We will refer to our method as Dynamic
Lattice FPKMC or DL-FPKMC.

\subsection{The SDLR and RDME Models}

In this subsection we compare the SDLR and RDME models to provide
motivation for our choice of the SDLR model. In the standard SDLR
approach~\cite{SmoluchowskiDiffLimRx}, the positions of molecules are
modeled as point particles undergoing Brownian motion.  The state of a
chemical system is given by the collection of stochastic processes for
the positions of each molecule of each chemical species at a given
time.  Bimolecular, or second-order, reactions between two molecules
are modeled as occurring either instantaneously, or \comment{a with an
intrinsic reaction-rate}, when the molecules' separation reaches a
specified reaction radius~\cite{KeizerJPhysChem82}.  Bimolecular
reactants are not allowed to approach closer than their reaction
radius.  Unimolecular, or first-order, reactions involving a single
molecule represent internal processes, such as decay or splitting, and
are assumed to occur with specified probabilities per unit time. (The
Doi model is similar, but allows bimolecular reactants to approach
arbitrarily close to each other. Bimolecular reactions then occur with
a fixed probability per unit time when the separation between two
molecules is \emph{less than} the
reaction-radius~\cite{TeramotoDoiModel1967,DoiSecondQuantA,DoiSecondQuantB}.)
The SDLR model can be extended to incorporate drift due to potentials.
The underlying reaction process remains unchanged, but molecules move
by a drift-diffusion process instead of pure Brownian motion.
\comment{We give a more detailed description of the SDLR model with
  drift in Section~\ref{S:overview}, where we also discuss the
  different bimolecular reaction mechanisms and where reaction
  products should be placed.}

In contrast to the spatially-continuous SDLR model, space in the RDME
is partitioned by a mesh into a collection of voxels. The diffusion of
molecules is approximated by a continuous-time random walk on the
mesh, with bimolecular reactions occurring with a fixed probability
per unit time for molecules within the same voxel.  Unimolecular
reactions are modeled in the same manner as in the SDLR model. The
state of a chemical system is then described by the collection of
stochastic processes for the number of each chemical species within
each voxel at a given time. Molecules are assumed to be well-mixed
within each voxel (i.e. uniformly distributed). The RDME can be
considered an extension of the chemical master equation (CME), a
standard non-spatial model for stochastic chemical
kinetics~\cite{GardinerRXDIFFME,McQuarrieJAppProb,GardinerHANDBOOKSTOCH,VanKampenSTOCHPROCESSINPHYS,RamaswamyNatCommun2012}.
One advantage of the RDME is the ability to construct a hierarchy of
more macroscopic approximations for which efficient numerical solution
methods have been developed~\cite{Rossinelli2008136}.

When the RDME is interpreted as an independent physical
model~\cite{ErbanChapman2009,RamaswamyJChemPhys2011} there is a
nonzero, lower bound on the lattice spacing that arises from the
assumption that the diffusive mixing timescale within a voxel is
faster than the timescale for a (well-mixed) bimolecular reaction to
occur~\cite{IsaacsonSJSC2006,ElfPNASRates2010}. There is also a
simultaneous upper bound on the lattice spacing to ensure the random
walk approximation of molecular diffusion is accurate.  Only when
these bounds are both satisfied is the RDME considered ``physically
valid''.  However, in the absence of nonlinear reactions, the solution
to the RDME should \emph{converge} to that of the SDLR model. This can
be seen by considering the particle-tracking representation of the
RDME derived in~\cite{IsaacsonRDMENote}. It is also standard to choose
the diffusive hopping rates in the RDME so that, in the absence of any
reactions, the solution to the RDME recovers the Brownian motion of
point particles as the lattice spacing is taken to
zero~\cite{IsaacsonSJSC2006,LotstedtFERDME2009}.  For these reasons,
in applications the RDME is often considered an approximation to the
more microscopic SDLR
model~\cite{RamaswamyJChemPhys2011,ElfPNASRates2010,Hellander:2012jk}.

For systems that include bimolecular reactions, interpreting the RDME
as an approximation to the SDLR model can be problematic. It has been
proven that in the continuum limit that the RDME lattice spacing is
taken to zero bimolecular reactions are lost in two or more
dimensions~\cite{IsaacsonRDMELims,IsaacsonRDMELimsII,Hellander:2012jk}.
As such, the time required for two molecules to react becomes
infinite.  The loss of bimolecular reactions is consistent with the
physical lower bound on the lattice spacing, and demonstrates that the
RDME can only be interpreted as an approximation to the SDLR model for
lattice spacings that are neither too large nor too small.  The error
of this approximation cannot be made arbitrarily small.
In~\cite{IsaacsonRDMELimsII} the two molecule $\textrm{A} + \textrm{B}
\to \varnothing$ reaction was studied in $\R^3$ in both the RDME and
SDLR model.  It was shown that for certain biologically relevant
parameter values the reaction time distribution in the RDME could
\emph{at best} approximate the reaction time distribution in the SDLR
model to within 5-10\% percent for an optimal choice of lattice
spacing.

Several recent efforts have derived renormalized bimolecular reaction
rates for use in the RDME that are designed to accurately capture
\emph{one specific statistic} of the SDLR model over a range of
sufficiently large lattice
spacings~\cite{ElfPNASRates2010,Hellander:2012jk}. For
example,~\cite{ElfPNASRates2010} matches the mean equilibration time
for the $\textrm{A} + \textrm{B} \leftrightarrows \textrm{C}$ reaction
in a system with one \textrm{A} molecule and one \textrm{B} molecule.
More recently, a modified convergent RDME (CRDME) that approximates
the Doi model was proposed in~\cite{IsaacsonCRDME}. Still, as of yet
there is no RDME-like approximation of the (microscopic) SDLR model in
which the approximation error can be made arbitrarily small.

There are several additional challenges to using RDME-like lattice
models to study cellular processes in realistic geometries.  Foremost,
\comment{existing methods for formulating RDMEs} approximate the
domain geometry by either unstructured
meshes~\cite{LotstedtFERDME2009} or Cartesian grid embedded boundary
methods~\cite{IsaacsonSJSC2006}. In the former it can be difficult to
construct meshes for which a continuous-time random walk approximation
of diffusion is well-defined at all mesh voxels when in three
dimensions~\cite{LotstedtFERDME2009}.  The latter tends to lose
accuracy in voxels cut by the boundary~\cite{IsaacsonSJSC2006}.
\comment{This may reduce the overall order of accuracy of the
  resulting spatial discretization that determines the spatial hopping
  rates in the RDME, thus reducing the accuracy of the RDME in
  approximating the Brownian motion of molecules, see ~\cite{IsaacsonSJSC2006}. }
  
For the preceding reasons, in this work we have chosen to focus on
developing a convergent numerical method for directly approximating
the SDLR model.  In particular, our approach does not constrain
molecules to remain on a fixed lattice as in the RDME.  While our
method is presented in one dimension to illustrate that it is
converging to the SDLR model, we expect that the method should be
well-suited for handling complex geometries in two and three
dimensions by \comment{treating boundaries as in the Walk on
  Rectangles algorithm~\cite{LejayExactWOS2006}. The latter was
  designed to generate exact realizations of the Brownian motion of
  molecules in domains with piecewise linear (planar) boundaries in 2D
  (3D).  The FPKMC can be interpreted as an extension of the Walk on
  Spheres~\cite{MullerWOSAlgo56} and Walk on
  Rectangles~\cite{LejayExactWOS2006} methods to systems that involve
  chemical reactions in simple geometries. We therefore expect it will
  be straightforward to adapt the FPKMC to handle complex geometries
  by incorporating the methods of~\cite{LejayExactWOS2006}.}

Exact numerical realizations of the stochastic processes described by
the RDME can be generated by the well-known Stochastic Simulation
Algorithm (SSA) method~\cite{GillespieJPCHEM1977}, which is a Kinetic
Monte Carlo (KMC) method~\cite{KalosKMC75}. More recently, the
First-Passage Kinetic Monte Carlo (FPKMC) method was developed to
generate exact realizations of the stochastic processes described by
the SDLR model~\cite{KalosDSMC06, OppelstrupPRE2009, DonevJCP2010,
  WoldeEgfrdPNAS2010}. \comment{The method developed 
  in~\cite{KalosDSMC06, OppelstrupPRE2009, DonevJCP2010} can 
  generate exact realizations of the SDLR model when molecules react 
  instantly upon reaching a fixed separation (a pure-absorption reaction). 
  In~\cite{WoldeEgfrdPNAS2010} this method was modified to generate
  exact realizations of the SDLR model with a more microscopic
  partial-absorption reaction mechanism.  There molecules (possibly)
  react based on a specified intrinsic reaction-rate upon collision.
  This modified version of the FPKMC algorithm was called Green's
  Function Reaction Dynamics (GFRD) in~\cite{WoldeEgfrdPNAS2010}.  In
  this manuscript we use FPKMC to refer to any of the methods
  of~\cite{KalosDSMC06,OppelstrupPRE2009,DonevJCP2010,WoldeEgfrdPNAS2010}.
  (Note, the earlier GFRD method of~\cite{WoldeGFRD05} generates approximations
  to the stochastic process described by the SDLR model, but differs from that
  in~\cite{WoldeEgfrdPNAS2010}.  The former is not exact, as it
  assumes the existence of a timestep below which a system of reacting
  molecules can be decoupled into one and two-body problems,
  see~\cite{WoldeGFRD05}.)}

An exact FPKMC method has also been introduced incorporating spatially
and temporally varying transition or annihilation rates for single
particles~\cite{SchwarzRieger2013}, which could be used to simulate
transitions from diffusive to ballistic modes in models of
intracellular transport~\cite{LoverdoNatPhys2008}. \comment{The newer
  GFRD method~\cite{WoldeEgfrdPNAS2010} has been modified to allow for
  advection due to a \emph{spatially-uniform, constant} velocity field
  along a one-dimensional track~\cite{GFRDAdvect2010}.}  It should be
noted that there are a number of alternative numerical methods,
\comment{such as Brownian Dynamics,} that have also been proposed for
approximating the stochastic processes described by the SDLR or Doi
models~\cite{AndrewsBrayPhysBio2004,ErbanChapman2011,
  HellanderGFRRD2011}.  Each of these methods is approximate, and
ultimately based on a discretization in time.

While it is conceptually simple to modify the SDLR model to include
drift due to potentials, \comment{it is still an open question how
  best to generate} numerical realizations of the underlying
stochastic processes described by the model. In this work we propose a
dynamic lattice version of the FPKMC method to allow for the inclusion
of drift due to a potential. While our Dynamic Lattice FPKMC method
(DL-FPKMC) no longer generates exact realizations of the SDLR model,
the error introduced in the method is controlled. The method can be
extended in a straightforward manner to also include spatially varying
diffusion coefficients.

\subsection{A DL-FPKMC Method for the SDLR Model with Drift}

\comment{In this subsection we present an overview introducing the DL-FPKMC
method. In the body of the paper, we will develop the details of the method
and analyze convergence.} 

The central idea of the original FPKMC method is to enclose one or two
molecules within a `protective domain' that isolates them from all
other molecules of the system.  A significant change in the state of
the system occurs only when a molecule leaves its protective domain or
a reaction occurs. In the cases discussed below, first-passage time
distributions for such events can be determined from information about
the underlying physical system.  The first-passage time distributions
can then be sampled to determine the time and type of state change
that occurs next. This event-driven approach provides an especially
efficient simulation algorithm by allowing each update of the
algorithm to span the time interval between significant state changes
(as opposed to proceeding over many small fixed time steps during
which the change in relevant state variables is minor).

In the rather special case of pure Brownian motion in simple domains
(spheres or rectangular regions), the first-passage time distributions
for a molecule to leave a protective domain~\cite{KalosDSMC06,
  OppelstrupPRE2009,DonevJCP2010,WoldeEgfrdPNAS2010} or for two
molecules to reach a threshold radius for
reaction~\cite{WoldeEgfrdPNAS2010} can be computed analytically by
solving the diffusion equation. The use of these expressions allows
for the generation of exact realizations of the stochastic process
described by the SDLR model with the FPKMC.  However, for many
situations in cell biology, pure Brownian motion does not provide the
most realistic description of the movement of molecules as a
consequence of active transport, chemical gradients, interactions with
cellular structures, etc. In such cases, significant drift terms are
inherent to the particle dynamics and can be modeled as arising from a
fixed potential field.  The DL-FPKMC method we develop extends the
FPKMC to allow for such drift.  Analytical expressions for the
first-passage time distributions from protective domains are no longer
possible with the addition of spatially varying drift. In DL-FPKMC, we
therefore approximate the drift-diffusion process each molecule
undergoes within a protective domain by a continuous-time random walk
on a discretized mesh.  The ``hopping rates'' for these walks are
\comment{determined from the Wang--Peskin--Elston}
finite-difference discretization of the Fokker-Planck
equation~\cite{ElstonPeskinJTB2003}. 
\comment{When new protective domains are created during the
  course of a simulation, they are dynamically meshed.}  For this
reason, our method can be interpreted as a dynamic-lattice master
equation model.  Unlike the standard RDME, it has the benefit of
converging to the SDLR model as the lattice spacing is reduced.

We present results demonstrating both the convergence and accuracy of
the DL-FPKMC method in one dimension as the mesh spacing in the
discretization is decreased.  In particular, we apply our algorithm to
the bimolecular reaction $A + B \to \varnothing$ where the molecules
of species $A$ and $B$ undergo drift-diffusion subject to various
types of potential functions (smooth, discontinuous, and constant).
Our results indicate that the method is approximately second-order
accurate for smooth potentials and approximately first-order accurate
for discontinuous potentials.  In this paper, we focus on convergence
of the DL-FPKMC in one dimension because of the availability of exact
analytical solutions and high-accuracy numerical solutions with which
to assess the error.  For a one-dimensional domain containing only one
molecule of species $A$ and one molecule of species $B$, the SDLR
model for the reaction system $A + B \to \varnothing$ can be described
by a single two-dimensional PDE for the probability density the
particles have not reacted and are located at specified positions.
This PDE can be solved numerically to high accuracy with general
potential fields using finite difference discretizations, and can
solved analytically when the potential is constant.  Having these
numerical and analytic solutions allows us to check both the accuracy
and the convergence of the DL-FPKMC simulation results. If the domain
containing the two molecules were instead two- or three-dimensional,
the corresponding PDE would be four- or six-dimensional, respectively,
and would therefore be challenging to solve to \emph{high accuracy} by
standard PDE discretization techniques. For the two-molecule system in
one dimension, we provide detailed numerical results demonstrating the
convergence of DL-FPKMC to the SDLR model.  For systems with more than
two molecules, where high-accuracy solutions to the equations for the
probability density of being in a given state are not available, we
show more qualitative convergence results.

The FPKMC method was originally presented as an efficient way to
simulate reaction-diffusion systems at low particle densities
``without all the hops'' by using larger
``superhops''~\cite{KalosDSMC06}.  While DL-FPKMC uses more hops than
FPKMC due to the random walk approximation of molecular motion, we
demonstrate that DL-FPKMC maintains efficiency at low particle
densities by requiring far fewer hops than fixed lattice methods with
comparable resolution.  By discretizing each individual protective
region, DL-FPKMC allows fine meshes to be used in localized regions
when needed for accuracy considerations. Examples where fine meshes
may be necessary include resolving bimolecular reactions, boundary
conditions, or rapidly varying potential fields.  For protective
domains in which such features are not present, coarser meshes can be
used. In this way the DL-FPKMC offers an alternative to the types of
global adaptive mesh methods that have been proposed for RDME-based
models~\cite{Bayati201113}.

The paper is organized as follows.  Section~\ref{S:overview} presents
our approach for incorporating drift into the SDLR model.  In
Section~\ref{S:FPKMCSect} we give an overview of the implementation
and steps of the FPKMC or DL-FPKMC algorithm to generate realizations
of the stochastic process described by the SDLR model. The general
discussion in Sections~\ref{S:overview} and~\ref{S:FPKMCSect} assumes
the molecules move in $\R^{n}$. The specific implementation we
develop in Section~\ref{S:oneDimMethod} restricts the molecules to
intervals in $\R$, but our technique can be extended to higher
dimensions through use of the Walk on Rectangles
method~\cite{LejayExactWOS2006}.  Section~\ref{S:oneDimMethod}
presents our numerical method for using a dynamic lattice to
incorporate drift into the FPKMC algorithm, and in
Section~\ref{S:methConv} we demonstrate the convergence and accuracy
of this DL-FPKMC method.  In Section~\ref{S:RunningTime} we provide a
running time analysis of DL-FPKMC, in which we demonstrate $O(N)$
scaling with the number of molecules in the system and compare
DL-FPKMC to a fixed lattice method.  Section~\ref{S:Applications}
\comment{presents several illustrative} applications \comment{of
  DL-FPKMC}.  In~\ref{S:ComparisonPotentials} we compare the effects
of drift due to several potentials on reaction time and location
statistics.  We conclude in~\ref{S:DNAproteinApplication} by
investigating a simplified model of a coupled protein-polymer fiber
system, in which two \comment{reacting molecules undergo}
drift-diffusion along a polymer, and may also unbind from the polymer
and diffuse in three dimensions.  We study the interaction between
\comment{polymer geometry, binding potentials along the polymer, and
  unbinding rate.}


\section{Incorporating Drift into the SDLR Model} \label{S:overview} 

\comment{In this section a modification of the SDLR model
  incorporating drift due to a potential is presented.}  The SDLR
reaction-diffusion model\comment{, with or without drift,} can be
described by a system of partial integro-differential equations
(PIDEs) for the probability densities of having a given number of
molecules of each chemical species at a specified set of positions,
similar to the stochastic reaction-diffusion PIDE models
in~\cite{IsaacsonRDMENote} and~\cite{DoiSecondQuantA,DoiSecondQuantB}.
Alternatively, one can consider the collection of stochastic processes
for the numbers of molecules and positions of each molecule of each
chemical species in the system. Due to the high-dimensionality of the
system of PIDEs for the probability densities, numerical methods for
solving the SDLR model, including FPKMC methods, are typically based
on Monte Carlo approaches that approximate the underlying stochastic
processes.  \comment{This section gives the mathematical formulation
  of the stochastic processes underlying the SDLR
  reaction-\emph{drift}-diffusion model.}  In the remainder of the
paper, we develop a Dynamic Lattice First-Passage Kinetic Monte Carlo
method (DL-FPKMC) for generating realizations of \comment{these}
stochastic processes.  The \comment{DL-FPKMC} realizations are
approximate, but the error is small and goes to zero as the lattice
spacing is decreased.

\subsection{\comment{Drift-Diffusion in the SDLR Model}}

In the modified SDLR model that includes drift, molecules are modeled
as points or hard spheres undergoing drift-diffusion processes.  In a
system with $K$ chemical species, we label the $k^{\textrm{th}}$
chemical species by $S^k$, $k = 1,\dots,K$.  We denote by 
\comment{$ \{ M^k(t): t \ge 0 \}$} the
stochastic process for the number of molecules of species $S^k$ at
time $t \comment{\ge 0}$.  The position vector of the $l^{\textrm{th}}$ molecule of
species $S^k$ at time $t \comment{\ge 0}$ is given by the vector stochastic process
\comment{$ \{ \vQ_{l}^{k}(t): t \ge 0 \}$, where}
$\vQ_{l}^{k}(t) \in \R^{n}$, $l = 1, \dots, M^k(t)$.

In the absence of any possible chemical reactions we assume the $l^{\textrm{th}}$ 
molecule of species $S^k$ with position $\vQ^k_l(t)$ undergoes 
diffusion with diffusion coefficient $D^k$ and experiences drift due 
to a potential $V^k \big(\vQ^k_l(t)\big)$. In this case, 
$\vQ^k_l(t)$ satisfies the stochastic differential equation (SDE)
\begin{equation} \label{eq:driftDiffSDEs}
  d \vQ^k_l(t) = \frac{-D^k}{k_{\textrm{B}}{T}}  \, \nabla V^k \big(\vQ^k_l(t)\big) \, dt + \sqrt{2 D^k} \, d \vec{W}^k_l(t),
\end{equation}
where $k_{\textrm{B}}$ is Boltzmann's constant, $T$ is absolute
temperature, $\nabla$ denotes the gradient in the coordinates of
$\vQ^k_l(t)$, and $\mb{W}^k_l(t)$ denotes the standard $n$-dimensional
Wiener process which describes Brownian motion. \comment{Note, this
  equation arises in the over-damped regime in which inertial forces
  are assumed to be negligible~\cite{GardinerHANDBOOKSTOCH}.}  Our
method for generating realizations of this stochastic process will be
presented in Section~\ref{S:oneDimMethod}.  The form of drift due to a
potential in Eq.~\eqref{eq:driftDiffSDEs} is useful for modeling
environmental interactions such as volume exclusion, attractive DNA
forces felt by regulatory proteins, and effective potentials felt by
molecular motors. It does not allow the possibility of interactions
between the diffusing molecules of the chemical system. To incorporate
such potentials into our method would be feasible for short-range pair
interactions, but for simplicity we assume no potential interactions
between molecules.

\subsection{\comment{Reactions in the SDLR Model}}

When bimolecular and unimolecular reactions are added into the system,
the molecules continue to move by the drift-diffusion
process in Eq.~\eqref{eq:driftDiffSDEs} under the additional constraint that
any pair of bimolecular reactants are not allowed to approach closer
than their corresponding reaction radius.  \comment{Upon reaching this
  reaction radius there are two alternative bimolecular reaction
  mechanisms that are commonly used in the SDLR model. In a pure
  absorption reaction the molecules react instantly upon reaching this
  separation~\cite{AndrewsBrayPhysBio2004}, while in a partial
  absorption reaction the molecules are allowed the additional
  possibility of having a non-reactive
  collision~\cite{ElfPNASRates2010,WoldeEgfrdPNAS2010}.  The latter
  model typically uses an intrinsic reaction rate constant to
  determine the probability of reaction upon collision
  (see~\cite{ElfPNASRates2010,KeizerJPhysChem82} for possible methods
  to relate macroscopic and intrinsic bimolecular reaction rates).
  Mathematically, a pure absorption reaction is modeled through a
  Dirichlet boundary condition in the corresponding system of
  PIDEs that describe the SDLR model, while a partial absorption
  reaction is modeled by a Robin boundary
  condition~\cite{KeizerJPhysChem82}. In this work, for simplicity we
  use the pure absorption mechanism.  While this is a more specialized
  model than the partial absorption reaction, a large number of
  modeling efforts have made use of the  Smoldyn Brownian Dynamics
  simulator~\cite{AndrewsBrayPhysBio2004}, which incorporates the pure absorption reaction
  mechanism (see~\cite{SmoldynWebsite} for
  references to these modeling studies).}
Unimolecular reactions representing internal processes are modeled as
occurring with exponentially distributed times based on a specified
reaction-rate constant. 

We will assume that reaction products are placed at the locations
specified in Table~\ref{tab:RxnProductLocations}.  In the cases of two
reaction products, for either unimolecular or bimolecular reactions,
the angular orientation of the product separation vector about the
center of mass is chosen randomly, \comment{as
  in~\cite{AndrewsBrayPhysBio2004}. More precise models for where to
  place the products would require extending the SDLR model to include
  the individual shapes, sizes, and/or rotational diffusion of
  molecules.}
 
When using the partial absorption Robin boundary condition-based
bimolecular reaction mechanism, the separation distance \comment{upon
  unbinding} can be chosen to be exactly the reaction radius.  In the
case that two products react by the pure absorption Dirichlet boundary
condition-based mechanism, such a choice would lead to an immediate
re-association reaction. This issue is generally handled through the
introduction of an unbinding radius~\cite{AndrewsBrayPhysBio2004}.
\comment{The unbinding radius is usually calibrated to enforce
  specified geminate recombination
  probabilities~\cite{AndrewsBrayPhysBio2004} (i.e. the
  probability that after unbinding a pair of molecules rebinds instead
  of diffusing away from each other).}  In the remainder, we restrict
  our focus to irreversible bimolecular reactions and assume reactions
  occur immediately when the separation between two reactants equals a
  specified reaction radius (a pure-absorption reaction). \comment{While 
  we focus on this special case for simplicity, a minor modification of
  our discretization procedure within protective domains would allow
  the use of a partial-absorption bimolecular reaction mechanism.}

\begin{table}
\vspace{.5cm}
\centering
\caption{\label{tab:RxnProductLocations} Placement of reaction products.}
{\renewcommand{\arraystretch}{1.5} 
\renewcommand{\tabcolsep}{0.2cm}
\begin{tabular}{| p{1.9cm} !{\vrule width 2pt} p{4.3cm} | p{8cm} |} 
  \hline  
& \emph{One Reaction Product} & \emph{Two Reaction Products} \\ 
  \noalign{\hrule height 2pt}
  \emph{Unimolecular Reaction} & 
  The product is placed at the same location as the reactant. & 
  The products are placed a specified distance apart, 
  with their center of mass at the location of the reactant. \\ 
  \hline
  \emph{Bimolecular Reaction} & The product 
  is placed at the center of mass of the two reactants.
  & The products are placed a specified distance apart, 
  with their center of mass at the same location as the 
  center of mass of the reactants. \\
\hline 
\end{tabular}  }
\vspace{.5cm}
\end{table}


\section{First-Passage Kinetic Monte Carlo Methods} 
\label{S:FPKMCSect} 

\comment{
This section gives an overview of First-Passage Kinetic Monte Carlo 
methods, including DL-FPKMC
and earlier methods, and then presents the steps for implementing
FPKMC or DL-FPKMC.}

\subsection{Overview of FPKMC Approaches}

First-Passage Kinetic Monte Carlo methods have been developed to
generate exact realizations of the SDLR model in the absence of
drift~\cite{KalosDSMC06,OppelstrupPRE2009,
  DonevJCP2010,WoldeEgfrdPNAS2010,SchwarzRieger2013}.  These novel algorithms are 
based on the Walk on Spheres method for solving exit time problems in
complicated geometries~\cite{MullerWOSAlgo56}.  They rely on being
able to derive exact analytical solutions of the diffusion equation in
spheres and rectangular solids.  In these FPKMC algorithms, a
spherical or rectangular region called a `protective domain' is drawn
around every molecule in the system, with the collection of protective
domains chosen to be disjoint.  The first-passage time for each
molecule, meaning the time when the molecule will first hit the
boundary of its protective domain, can be sampled exactly using the
corresponding analytical solution to the diffusion equation.  The
molecule that exits its protective domain first is updated to its exit
position, and a new protective domain is defined.  When two reactants
are sufficiently close that a sphere can be drawn that contains only
them, the corresponding two-body solution to the SDLR PIDEs can be
used to exactly sample a candidate time and location for their
reaction~\cite{WoldeEgfrdPNAS2010}.

A central feature of each of these prior methods is the assumption
that the molecular species undergo purely diffusive stochastic
dynamics within their protective domains. When drift due to a
potential is present in addition to diffusion, the probability
densities for the locations of one or two molecules within
  their protective domains are no longer described by the diffusion
equation, but rather by a Fokker-Planck equation.  Let
  $\Omega \subset \R^{n}$ denote the overall domain with boundary
  $\partial \Omega$, and let $U \subset \Omega$ label a protective
  domain with boundary $\partial U$. We are interested in the time a
  molecule first leaves $U$, leading to a zero Dirichlet boundary
  condition on $\partial U$. Let $V(\vx)$ denote the strength of the
  potential at $\vx \in \Omega$.  The probability density,
  $\rho(\vx,t)$, for a single molecule 
  \comment{with  diffusion coefficient $D$}
  to be at location $\vx$ within
  its protective domain $U$ at time $t > 0$ evolves according to the
  equations
\begin{align} \label{eq:FPE_Rn_singlePD}
\frac{\partial \rho( \vx ,t)}{\partial t} = D   \nabla \cdot  
\left(  \rho( \vx ,t) \frac{ \nabla  V(\vx )} {k_{\textrm{B}} T}  +  \nabla \rho(\vx ,t)  \right),   
&  \qquad \text{on }  U, \\
\rho (\vx ,t)= 0,  &  \qquad  \text{on }  \partial U \backslash  (\partial U \cap \partial \Omega),   \notag   \\
\rho(\vx,0) = \delta(\vx - \vx_0),  \notag  &
\end{align}   
where $\vx_0$ is the initial position of the molecule within the
protective domain.  If $\partial U$ intersects $\partial \Omega$, the
boundary conditions on $\partial U \cap \partial \Omega$ will agree
with the boundary conditions on $\partial \Omega$.  

For two bimolecular reactants within one protective domain $U
\comment{\subset \Omega \subset \R^{n}}$, the joint probability density 
$\rho( \vx, \vy ,t)$ for one molecule to be at location 
$\vx \comment{\in U}$ and the other molecule to be location 
$\vy \comment{\in U}$ at time $t > 0$ \comment{can be described by 
  a Fokker-Planck equation in $\R^{2n}$. Define the domain in 
  $\R^{2n}$ to be $W = \{  (\vx, \vy):  \vx \in U, \ \vy \in U,  \ 
  \text{and} \ || \vx - \vy|| > r_{\textrm{R}} \}$. $D_1$ and $D_2$ 
  will denote the respective diffusion coefficients of the two molecules. 
  We define $\mathcal{D}$ to be a $2n \times 2n$ diagonal matrix, 
  with the first $n$ elements of the diagonal equal to $D_1$ and the 
  next $n$ elements equal to $D_2$. Let} $V_{1}(\vx)$ and $V_2(\vy)$ 
\comment{be} the potential fields that impart drift to the $\vx$ molecule 
and the $\vy$ molecule respectively, \comment{and define $V(\vx, \vy) 
  = V_1(\vx) + V_2(\vy)$. The behavior on the boundary $ \partial W$ 
  may be different on each of the following three components: the 
  reactive boundary $ \partial W_\textrm{rxn} = \partial W \cap 
  \{|| \vx - \vy || = r_{\textrm{R}}\}$, the (possibly empty) intersection 
  with the overall domain boundary $ \partial W_\textrm{outer} = 
  (\partial W  \cap \{ (\vx, \vy): \vx \in \partial \Omega \,  \text{ or } \, 
  \vy \in \partial \Omega \} ) \backslash \partial W_\textrm{rxn} $,
  and the remaining component $ \partial W \backslash 
  (\partial W_\textrm{rxn} \cup \partial W_\textrm{outer} )$. Then, 
  $\rho( \vx, \vy ,t)$  satisfies}
\begin{align} \label{eq:FPE_Rn_pairPD}
\frac{\partial \rho( \vx, \vy ,t)}{\partial t} = \comment{ \nabla \cdot  \mathcal{D} }
\left(  \rho( \vx, \vy ,t) \frac{ \nabla  \comment{V(\vx , \vy)} } {k_{\textrm{B}} T}  +  \nabla \rho(\vx, \vy ,t)  \right),   
&  \qquad \text{on }  \comment{W}, \\
\rho(\vx,\vy,t) = 0,  &  \qquad  \text{on }   \comment{ \partial W_\textrm{rxn} },   \notag  \\
\rho (\vx, \vy, t)= 0,  &  \qquad  \text{on }  \comment{\partial W 
\backslash (\partial W_\textrm{rxn} \cup \partial W_\textrm{outer} )},  \notag    \\
\rho(\vx, \vy, 0) = \delta(\vx - \vx_0, \vy - \vy_0), \notag    
\end{align}
where $\vx_0$ and $\vy_0$ are the initial positions of the molecules
within the protective domain, and the gradient and divergence
operators are in the $(\vx,\vy)$ coordinates. \comment{Similar to}
the single particle case, if $\partial U$ intersects $\partial \Omega$, 
the boundary conditions on \comment{$\partial W_\textrm{outer} $}
will agree with the boundary conditions on $\partial \Omega$.  The 
boundary condition $\rho(\vx,\vy,t) = 0$ when $||\vx-\vy||
= r_{\textrm{R}}$ models the pure-absorption reaction mechanism. 
This reactive boundary condition could be modified to use a Robin 
partial-absorption mechanism if desired. In the remainder we
assume all molecules \comment{have the same diffusion coefficient and}
experience the same potential field, so that \comment{$D_1 = D_2$ and}
$V_1=V_2$. In the case that the potential $V$ is a constant function, 
Eqs.~\eqref{eq:FPE_Rn_singlePD} and~\eqref{eq:FPE_Rn_pairPD} 
reduce to diffusion equations.

Our method for using discretizations of these Fokker-Planck equations to
sample times and locations of first-passage and reaction events is
described in Section~\ref{S:oneDimMethod}.


\begin{figure}
\centering
\includegraphics[width=15cm]{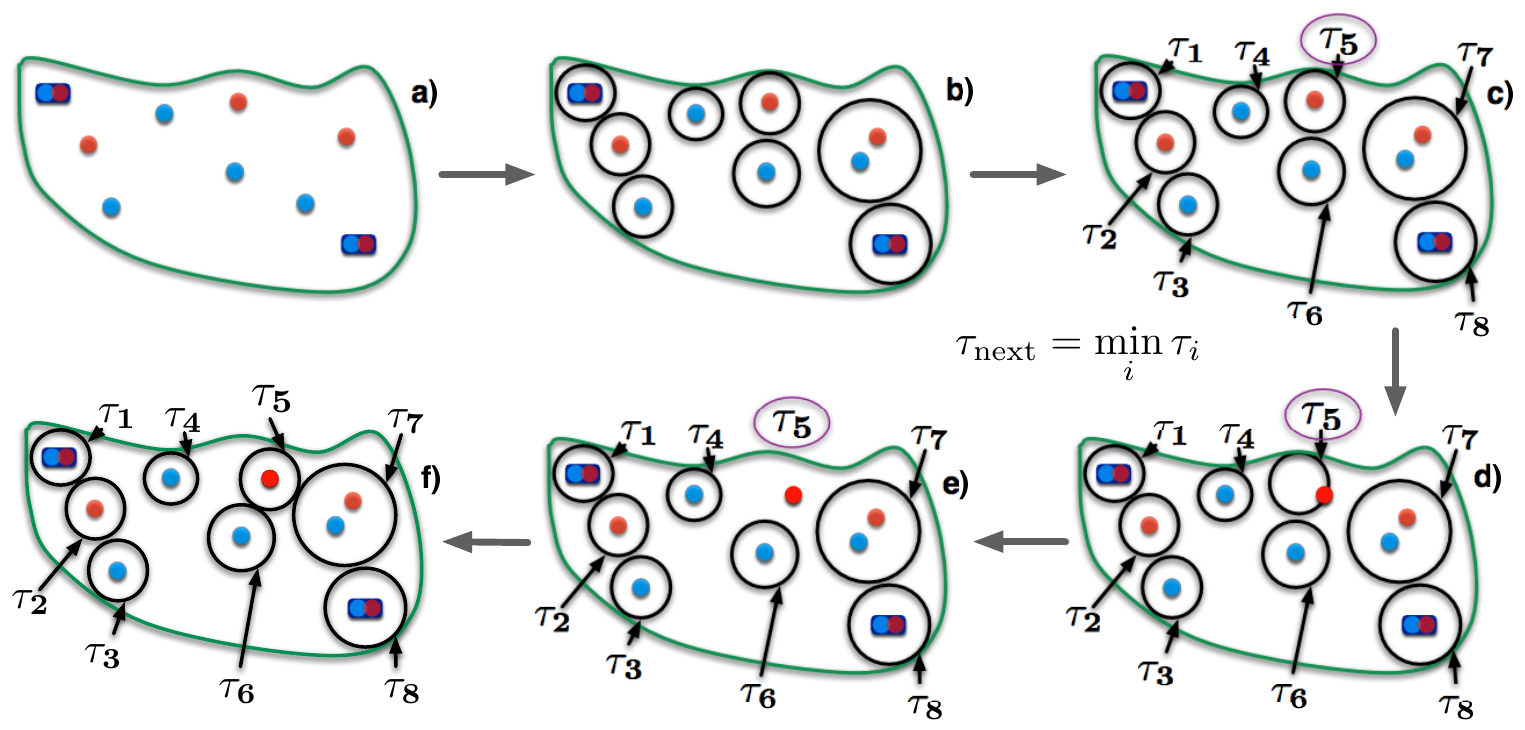}
\caption{\label{fig:fpkmcSchematic} \small Schematic of the FPKMC or
  DL-FPKMC algorithm as described in
  Subsection~\ref{StepsOfAlgorithm}. Here, $\tau_i$ denotes the next
  event time for the molecule or pair of molecules in the $i^{th}$
  protective domain, and $\tau_{\textrm{next}} = \tau_5$ is the time
  of a global event.}
\end{figure} 

\subsection{Main Steps of the FPKMC or DL-FPKMC Algorithm} \label{StepsOfAlgorithm}

In this subsection we describe the role of protective domains and the
processing of events in our implementation of the FPKMC or DL-FPKMC
algorithm.  We then list the main steps of the algorithm. Our
implementation is based on the FPKMC algorithm developed
in~\cite{KalosDSMC06,OppelstrupPRE2009,DonevJCP2010}, with some
modifications. The content of this section applies to both the FPKMC
and DL-FPKMC algorithms. The \comment{only theoretical} difference
between the two methods is the approach for sampling event times
within a protective domain. In FPKMC, which can only be used in the
case of constant $V(\vx)$, event times are sampled from \emph{exact}
solutions of the diffusion equation.  In DL-FPKMC, which allows for
arbitrary $V(\vx)$, event times are \emph{approximated} by generating
sample paths of continuous-time random walks on meshes within
protective domains (see Section~\ref{S:oneDimMethod}).

\subsubsection{\comment{Protective Domains}}
To apply the algorithm, every molecule in the system is placed in a
protective domain. In one dimension the protective domains are
intervals and in higher dimensions the protective domains are usually
rectangular or spherical regions. In general, the boundaries of
protective domains are absorbing. The boundary of a protective domain
can contain a portion of the boundary of the overall spatial domain,
in which case the protective domain boundary conditions will depend on
the overall domain boundary conditions.  We allow protective domains
to contain either one or two molecules. Protective domains containing
only one molecule are referred to as `single protective domains', and
those containing two molecules are referred to as `pair protective
domains'. Molecules in separate protective domains behave
independently. Each molecule undergoes drift-diffusion within its
protective domain, and may undergo unimolecular reactions. Two
molecules in the same protective domain may additionally participate
in bimolecular reactions. To maintain independence when bimolecular
reactants are in different protective domains, we require 
\comment{a separation of} at least one reaction radius.  For
non-reacting molecular species we allow for overlap to prevent the
size of protective domains from going to zero. \comment{Additional
  details on constructing and updating protective domains are provided
  in Appendix~\ref{S:ProtectiveDomainsAppendix}. }

\subsubsection{\comment{Events, Times, and Updating}}
Each event that may occur will have a type, time, and location.  The
two major event types are first passage from a protective domain and
reaction.  First passage from a protective domain occurs when a
molecule first reaches an absorbing boundary of its protective domain.
In DL-FPKMC with general $V(\vx)$ (resp. FPKMC with constant
  $V(\vx)$), times and locations for first-passage events from single
  protective domains are sampled from probability densities 
determined from  approximate (resp. exact) solutions of
  Eq.~\eqref{eq:FPE_Rn_singlePD}.  Similarly, for pair protective
  domains, solutions of Eq.~\eqref{eq:FPE_Rn_pairPD} are used to
  sample times and locations for first-passage events or for
  bimolecular reactants to first reach a separation of one reaction
  radius.  The time for a unimolecular reaction to occur is sampled
from an exponential distribution with a specified reaction rate, and a
corresponding reaction location is sampled from a
    ``no-passage'' probability density for the molecule involved (see
    below).

To facilitate the discussion of these \comment{first-passage and reaction} 
events, we use specific names 
for three times.  The `global time' will refer to when the most recent 
event has occurred, irrespective of its particular type or which 
molecules were involved.  An `individual time' and a `next event time' 
will be associated with each particular molecule. `Individual time' 
will refer to when the molecule was last updated,
and `next event time' will refer to the sampled time at which the 
molecule might next undergo an event. Individual times are less 
than or equal to the global time, and next event times are greater 
than the global time.

Usually, the individual time and location of a molecule will only be
updated when the molecule undergoes a major event (first-passage or
reaction). In this case, the time and location of the molecule will be
updated to the time and location of the event. However, a molecule can
also be updated to any specified time prior to its next event time, by 
sampling a new position for the molecule within the
  protective domain from the conditional probability density for the
  molecule to be at a position within the domain, at the specified
  time, and not yet have undergone a first-passage or reaction event.
  This procedure is called a `no-passage' update.
  
\subsubsection{\comment{Overall Algorithm}} \label{S:algoSteps}

The \comment{FPKMC or DL-FPKMC} algorithm is carried out according to the following steps:
\begin{enumerate}  [itemsep=0pt, topsep=1pt]
\item Protective domains are defined around each molecule or pair of 
   molecules, as shown in Figure~\ref{fig:fpkmcSchematic}b.   

\item The next individual event for each molecule or pair of molecules 
   is determined  by sampling an event type, time, and location.  In  
   Figure~\ref{fig:fpkmcSchematic}c, each next event time is labeled
   by a $\tau_i$.

\item To determine global events, the individual events are sorted in 
 a priority queue ordered from the shortest event time to the longest 
 event time.  For example, $\tau_5$ denotes the shortest event time 
 in Figure~\ref{fig:fpkmcSchematic}c.

\item The next global event is determined from the priority queue 
   using the next individual event with the shortest time. The global 
   time and the individual time(s) of the participating molecule(s) 
   are updated to the event time. In the case of a first-passage 
   event for a molecule to leave its protective domain, the molecule's 
   location is updated to the sampled first-passage location, as shown 
   in Figure~\ref{fig:fpkmcSchematic}d. If this molecule is in a pair 
   protective domain with another molecule, the other molecule is 
   no-passage updated to the new global time.
    In the case of a reaction event, the reaction products are placed 
    at or about the reaction location, as specified in Table~\ref{tab:RxnProductLocations}.

\item  Molecules in protective domains that are close to 
or overlap the newly updated 
molecules are no-passage updated to the new global time.  
 
\item New protective domains are constructed only for those molecules 
  that have undergone an update to reach the current global time, as 
  shown in  Figure~\ref{fig:fpkmcSchematic}f.  New events are sampled 
  for these updated molecules, and the event times are sorted into the 
  priority queue.  All other molecules and events remain unchanged.  

\item Steps 4 through 6 are then repeated.

\end{enumerate}
Note that Step 5 is used to keep the sizes of the protective domains
from becoming too small~\cite{DonevJCP2010}, in which case the
\emph{effective} time steps used in the FPKMC or DL-FPKMC methods
could become very short.

We remark that information about the state of any molecule in the system 
is available for any particular time in the simulation.
For instance, if one would like to sample the locations of all molecules 
at a specified time, this can obtained by taking the state of the system at
the largest global time before or equal to the specified time and then 
no-passage updating each molecule to the specified time.


\subsection{Protective Domain Changes during One Simulation}

\begin{figure} 
\centering
\includegraphics{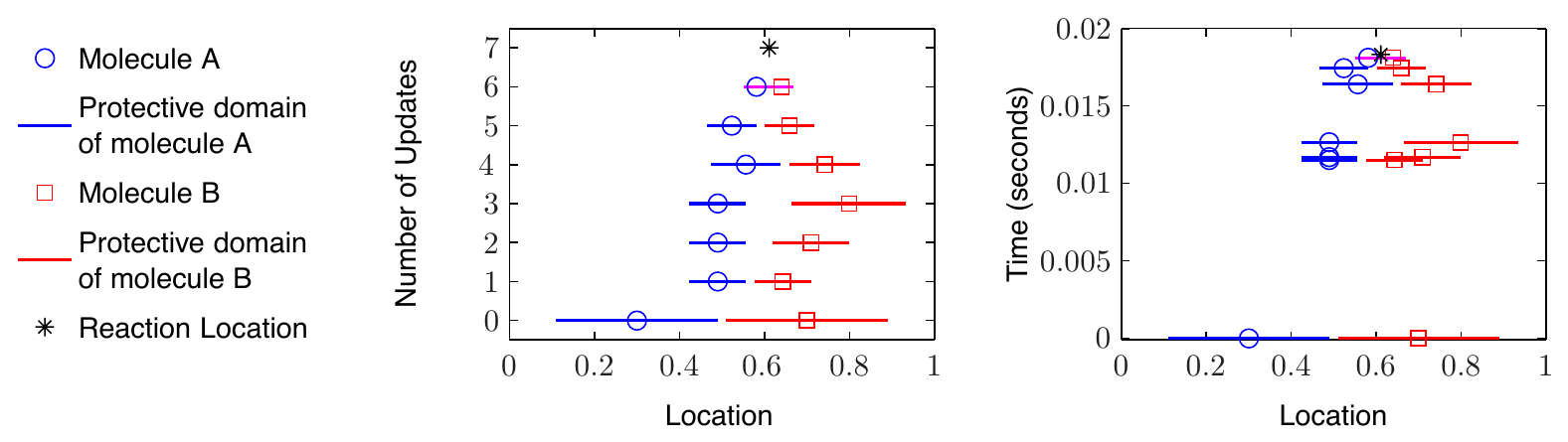}
\caption{\label{fig:SampleRun1} \small One simulation of the
  reaction $A + B \to \varnothing $, with one molecule each of $A$ and
  $B$ present initially and $V(x)=0$. Both panels show the same run of
  the simulation.  In the left panel, the vertical axis is the number
  of times the simulation cycled through Steps $4$ to $6$ of the
  algorithm.  In the right panel, the vertical axis is the time of the
  most recent event.}%
\end{figure}

During simulations, updates are made to the protective domains
sequentially as events occur changing the state of the system.  To
illustrate this process, we consider the simulation in one dimension
of the reaction $A + B \to \varnothing$ starting with one molecule
each of $A$ and $B$ and using our DL-FPKMC algorithm.  One 
simulation is shown in Figure~\ref{fig:SampleRun1}.

In the left panel of Figure~\ref{fig:SampleRun1}, the vertical axis is 
the number of times that the simulation cycled through Steps $4$ to
$6$ of the algorithm; we call this number $N_{\textrm{update}}$.  In 
the right panel, the vertical axis is the time of the most recent event.  
At $N_{\textrm{update}}=1$ in this particular run of the simulation, 
molecule $A$ is first-passage updated to the right endpoint of its 
initial protective domain.  This location is close to the left endpoint 
of molecule $B$'s protective domain, so molecule $B$ is no-passage 
updated and new protective domains are defined around each molecule.  
From $N_{\textrm{update}}=1$ to $N_{\textrm{update}}=3$, molecule 
$B$ is first-passage updated but does not come close to the protective 
domain of molecule $A$, so molecule $A$ is not updated. At 
$N_{\textrm{update}}=6$, the distance between molecules $A$ and 
$B$ is less than a specified pair threshold, so they are placed in a 
pair protective domain.  At $N_{\textrm{update}}=7$, the distance 
between the molecules reaches the reaction radius and the reaction 
occurs.


\section{\comment{Propagation of Molecules within Protective 
Domains  in One-Dimensional DL-FPKMC}} 
\label{S:oneDimMethod}

\comment{In this section we introduce a lattice discretization of the
  Fokker-Planck equation~\cite{ElstonPeskinJTB2003} which is used
  within each protective domain. The discretization is chosen to have
  the form of a master equation, so that the discretization weights
  can be interpreted as transition rates for continuous-time random
  walks by the molecules within each protective domain. In
  Subsection~\ref{S:SamplePaths} we describe how the Stochastic
  Simulation Algorithm (SSA)~\cite{GillespieJPCHEM1977,KalosKMC75} is
  then used to generate realizations of these random walks within each
  protective domain, giving the next event times and locations needed
  by the DL-FPKMC algorithm (see Subsection~\ref{S:algoSteps}).
  It should be stressed that our DL-FPKMC method modifies FPKMC by
  using this lattice method to propagate molecules \emph{within} their
  protective domains.}  To introduce our methods,
we consider the case where the simulation domain 
is a one-dimensional interval
\comment{(in higher dimensions, the same discretization can be used
in each coordinate)}.
We assume that bimolecular 
reactions occur instantaneously when the reactants' separation 
reaches the reaction radius. We allow the simulation domain 
to have reflecting, absorbing, or periodic boundaries. Reflecting 
boundaries are modeled using zero-flux boundary 
conditions, and absorbing boundaries are modeled using 
zero Dirichlet boundary conditions. Protective domains are 
proper subintervals of the overall domain.

\subsection{\comment{Lattice Discretization of the Fokker-Planck Equation} } 
\label{S:WPE_LatticeDiscretization}

\comment{In the case of} pure diffusion, the probability distributions for
 first-passage times, first-passage locations, and no-passage
locations can all be determined from analytic solutions of the
diffusion equation~\cite{KalosDSMC06, OppelstrupPRE2009, DonevJCP2010,
  WoldeEgfrdPNAS2010}.  In contrast, once drift is considered, such
analytic approaches are no longer possible in general. Instead one
must consider probability densities that satisfy Fokker-Planck
equations such as Eqs.~\eqref{eq:FPE_Rn_singlePD} 
and~\eqref{eq:FPE_Rn_pairPD}, 
in which we restrict to the case
where the drift arises from a spatially varying potential energy
function $V(x)$. 

\comment{To sample event times and locations}
in DL-FPKMC, we 
introduce approximations by treating the movement of each molecule
within its protective domain as a discrete-space continuous-time
Markov chain, more specifically a continuous-time random walk on
discrete mesh points.
Jump rates between neighboring mesh points are obtained using the
Wang--Peskin--Elston~\cite{ElstonPeskinJTB2003} (WPE) spatial
discretization of the one-dimensional Fokker-Planck equation
\begin{equation} \label{eq:FPE}
\frac{\partial \rho(x,t)}{\partial t} = D  \frac{\partial }{\partial x} 
\left(\rho(x,t) \frac{d V(x)}{d x}   + 
\frac{\partial \rho(x,t) }{\partial x} \right).
\end{equation}
We remark that the factor $(k_B T)^{-1}$ is absorbed into the
potential function $V$.  The WPE discretization weights then determine
the jump rates (i.e. probabilities per unit time) for
molecules to hop from one mesh point to another. In particular, the jump
rate for a molecule to hop from the mesh point $x_i$ to a neighboring
mesh point $x_j$, in the case of a uniform mesh of width $h$, is given
by
{\renewcommand{\arraystretch}{1.5}
\begin{equation} \label{eq:JumpRate}
a_{ij}  = \left\{
\begin{array}{cl} 
 \frac{D }{h^2} \frac{V(x_j) - V(x_i)} {\exp[V(x_j) - V(x_i) ] -1}  & \text{for }  V(x_i)\ne V(x_j) \\
 \frac{D }{h^2}    & \text{otherwise. }  \end{array} \right. 
\end{equation} } 
Let  $p_i(t)$ be the probability that a molecule is located at mesh 
point $x_i$ at time $t$. Then the time evolution of $p_i(t)$ is described
by the master equation
\begin{equation}  \label{eq:MasterEqn}
 \frac{d \, p_i(t)}{d t}  = 
a_{i-1,i} \, p_{i-1}(t)  -  (a_{i,i-1} + a_{i,i+1}) \, p_i(t) + a_{i+1,i} \, p_{i+1}(t).
\end{equation}
If $x_{i \pm 1}$ is an absorbing boundary, then $p_{i \pm 1}(t)=0$ 
in Eq.~\eqref{eq:MasterEqn}.
We shall extend Eqs.~\eqref{eq:JumpRate} and~\eqref{eq:MasterEqn} for 
non-uniform discretizations \comment{in Subsection~\ref{S:NonUnif}}
 and in Appendix~\ref{S:DerivationNonUnifJumpRates}.

The discretization given by Eq.~\eqref{eq:JumpRate} has
the following properties:
\begin{itemize}[itemsep=0pt, topsep=1pt]
\item 
Converges
at second-order for smooth potentials, and can handle
discontinuous potentials~\cite{ElstonPeskinJTB2003}. 
\item 
Satisfies a discrete version of detailed balance
(zero net flux at equilibrium), 
which helps reduce artificial drift due to numerical
discretization errors~\cite{ElstonPeskinJTB2003}.  
\item 
Is consistent with the standard 
second-order-accurate discretization of the Laplacian operator, 
in that $a_{ij}$ converges to $ D / h^2$ as $V(x_j) - V(x_i)$ 
approaches zero.
\item \emph{Can be extended to higher
dimensions.  The jump rates in each coordinate are then given
by Eq.~\eqref{eq:JumpRate}.  }
\item Can incorporate a spatially dependent
diffusion coefficient $D(x)$.  For example, in the case that $D(x)$ is
continuous, the constant $D$ in Eq.~\eqref{eq:JumpRate} can be replaced by
$[D(x_i) + D(x_j)]/2$~\cite{ElstonPeskinJTB2003}.
\end{itemize}
\vspace{3pt}

\subsection{\comment{Generating Sample Paths}} \label{S:SamplePaths}
To make use of the \comment{WPE} discretization, a mesh is defined
within each protective domain so that every molecule is located at a
mesh point. \comment{Rather than numerically solve the master
  equation~\eqref{eq:MasterEqn}, and then sample this solution to
  determine next event times, we generate realizations of the jump
  process described by Eq.~\eqref{eq:MasterEqn}.  Each molecule then
  undergoes a continuous-time random walk on the mesh, with the
  transition rate from a mesh point $x_i$ to a neighboring point
  $x_j$ given by $a_{ij}$ of Eq.~\eqref{eq:JumpRate}.}  Exact sample
paths of the molecules' random walks are generated using the
event-driven Stochastic Simulation Algorithm \comment{(SSA),
  specifically Gillespie's ``direct method" version of the
  SSA~\cite{GillespieJPCHEM1977}.}  In this method, the times of the
hops are sampled from exponential distributions. There is no fixed
time step. By varying the mesh width, the resolution of this process
can be adjusted depending on the desired trade-off between
computational efficiency and accuracy.  
Our specific approach for choosing the mesh width and the locations of
mesh points is described in more detail in
Subsection~\ref{S:ChoosingMeshWidth}.

Any protective domain endpoints on the interior of the overall domain
are absorbing, as are endpoints that coincide with an absorbing
boundary of the overall domain.  If one endpoint of a protective
domain is located at a reflecting boundary of the overall domain that
endpoint is made reflecting. \comment{Since protective domains are
  proper subintervals of the overall domain, each protective domain
  will have at least one endpoint on the interior of the overall
  domain.Thus,} each protective domain either has absorbing Dirichlet
boundaries ($\rho = 0$) at both endpoints, or an absorbing boundary at
one endpoint and a reflecting \comment{zero-flux} boundary ($ \rho
\frac{d V}{d x} + \frac{\partial \rho}{ \partial x} = 0$) at the
other.  

For a newly constructed protective domain containing a single
molecule, we determine the molecule's next event time by
\comment{using the SSA to sample} an exact random-walk path for the
molecule to hop on the mesh points until it reaches an absorbing
boundary of the protective domain.  The time that \comment{the}
molecule reaches an absorbing endpoint is the first-passage time, and
the endpoint that the molecule reaches is the first-passage location.
For pair protective domains with two molecules, we perform random
walks for each molecule until either: (i) one molecule reaches an
absorbing boundary of the protective domain; or (ii) the distance
between the two molecules is equal to the reaction radius
$r_\textrm{R}$.  The mesh width for pair protective domains is always
chosen to exactly divide $r_\textrm{R}$, so that the reaction occurs
when the two molecules are exactly one reaction radius apart.  A
no-passage location at any specified time before the next event time
can be obtained by finding the last time in the sample path less than
or equal to the specified time and taking the location of the molecule
at that time.

\begin{figure}
\centering
\scalebox{.78}{\includegraphics{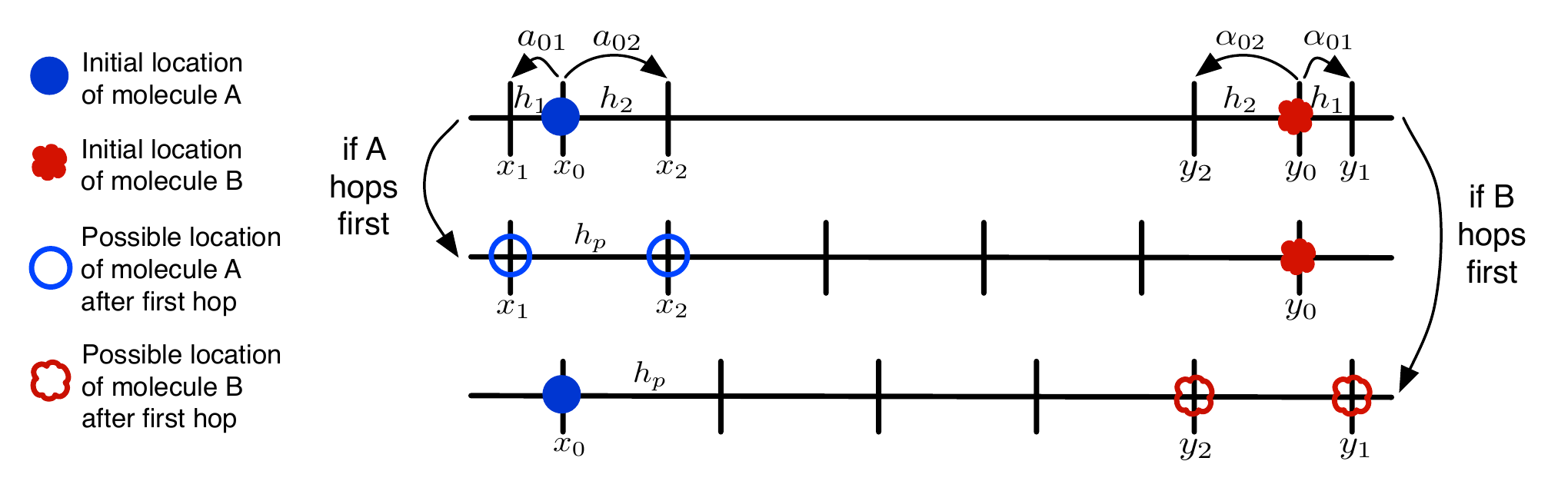}}
\caption{\label{fig:NonUnifSubLattice} \small The top row shows the
  non-uniform sub-meshes that are defined when two molecules are first
  placed in a pair protective domain. The second (resp. third) row
  shows the uniform mesh that is chosen if molecule $A$ (resp. $B$)
  hops first.  The \comment{user-specified} uniform mesh width, $h_p$,
  is chosen to exactly divide $r_\textrm{R}$.  \comment{Let $d$ be the
    initial distance between the two molecules at the time that the
    pair protective domain is defined. Then, $h_2$ is defined to be
    the remainder of the quotient $d/h_p$, and $h_1 = h_p - h_2$.
    By choosing $h_1$ and $h_2$ in this way,} the distances $|y_0 -
  x_1|$, $|y_0 - x_2|$, $|y_1 - x_0|$, and $|y_2 - x_0|$ are all
  exactly divisible by $h_p$.  The rates $a_{01}$ and $a_{02}$ are
  given by Eq.~\eqref{eq:NonUnifJumpRate}.  \comment{The SSA is used
    to simulate a single hop of one of the molecules} to a new point
  on its sub-mesh, \comment{after which} the new distance between the
  two molecules will be one of the four distances \comment{divisible
    by $h_p$ listed above}.  Then, a new mesh of uniform width $h_p$
  can be defined so that both molecules lie exactly on mesh points,
  \comment{and the generation of a sample path using the SSA proceeds
    until one molecule exits the pair protective domain or the
    distance between the two molecules reaches $r_\textrm{R}$}.}
\end{figure}

\subsection{\comment{Discretization for Non-uniform Mesh Cells}} \label{S:NonUnif}
Non-uniform mesh cells are used when needed to conform to a 
boundary or to move molecules onto a uniform mesh where the 
mesh width exactly divides the reaction radius, as will be described 
in Subsection~\ref{S:ChoosingMeshWidth}.  Let $x_0$ be the initial 
location of a molecule on a non-uniform mesh, with $x_1$ and 
$x_2$ denoting the locations of the neighboring mesh points 
in either direction. Note that we may have either $x_1 < x_0 < x_2$ 
or $x_2 < x_0 < x_1$ (see Fig.~\ref{fig:NonUnifSubLattice}, top row). 
Let $h_1 = |x_0 - x_1|$ and $h_2 = |x_0 - x_2|$.  The jump rates 
from $x_0$ to $x_j$ for  $ j = 1,2$ are given by
{\renewcommand{\arraystretch}{1.5}
\begin{equation} \label{eq:NonUnifJumpRate}
a_{0j}  = \left\{
\begin{array}{cl} 
\frac{2 D}{h_j (h_1+h_2) } \frac{V(x_j) - V(x_0) }{\exp[V(x_j) - V(x_0) ] -1} & \text{for }  V(x_i)\ne V(x_j) \\
\frac{2 D}{h_j (h_1+h_2) }   & \text{otherwise. }  \end{array} \right. 
\end{equation}}
The non-uniform discretization in Eq.~\eqref{eq:NonUnifJumpRate} is 
derived in Appendix~\ref{S:DerivationNonUnifJumpRates} by 
modifying the WPE discretization of the
Fokker-Planck equation~\cite{ElstonPeskinJTB2003}.  In the case 
of constant $V(x)$, the non-uniform discretization~\eqref{eq:NonUnifJumpRate} reduces to the 
non-uniform spatial discretization of the Laplacian at a Dirichlet boundary
  given by equation~(20) of~\cite{FedkiwJCP2002}.
  For solving the Poisson equation with Dirichlet boundary conditions, 
  using a uniform interior mesh and non-uniform mesh cells
  at the boundaries, this discretization is second-order accurate~\cite{FedkiwJCP2002}. 
To our knowledge, Eq.~\eqref{eq:NonUnifJumpRate} gives a 
new discretization of the Fokker-Planck equation for non-uniform 
meshes.


\subsection{Choosing the Mesh within Protective Domains} \label{S:ChoosingMeshWidth}

\comment{In this subsection we describe how a mesh is defined within
  each protective domain, given user-specified mesh widths.  The
  convergence tests in Section~\ref {S:methConv} will demonstrate that
  DL-FPKMC converges to the underlying SDLR model as the mesh widths
  are decreased.}

Single-molecule protective domains with absorbing boundaries are
chosen to be symmetric about the location of the molecule.
\comment{Let $r_\textrm{PD}$ be the distance from the molecule to
  either endpoint of the protective domain.}  A maximum mesh width for
all single protective domains, $h_s^{\textrm{max}}$, is specified
\comment{by the user}.
Then, for each individual protective domain, the mesh width $h_s$ is
\comment{calculated according to the formula
\begin{equation*}
h_s =  \frac{r_\textrm{PD} } {\textrm{ceil} (r_\textrm{PD} / h_s^{\textrm{max}} ). }
\end{equation*}
In this way, $h_s$ is always} chosen to be the largest value less than
or equal to $h_s^{\textrm{max}}$ that exactly divides
\comment{$r_\textrm{PD}$.  Generally $h_s$ will satisfy
  $h_s^{\textrm{max}}/2 < h_s \le h_s^{\textrm{max}}$, unless
  $r_\textrm{PD} \le h_s^{\textrm{max}}/2$, in which case $h_s=
  r_\textrm{PD} \le h_s^{\textrm{max}}/2$.  In practice, $h_s$ will
  almost always be strictly less than $h_s^{\textrm{max}}$, since it
  is unlikely that $h_s^{\textrm{max}}$ will exactly divide
  $r_\textrm{PD}$.  After calculating $h_s$, a uniform mesh with
  spacing $h_s$} is constructed so that the molecule and both
endpoints of the protective domain lie exactly on mesh points.  Having
the endpoints lie on mesh points allows enforcement of the absorbing
Dirichlet boundary conditions at the endpoints \comment{without
  modification of} the jumps
rates in Eq.~\eqref{eq:JumpRate}. 

\begin{figure}
\centering
\includegraphics[width=16cm]{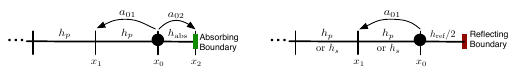}
\caption{\label{fig:NonUnifBoundary} \small Non-uniform mesh cell at
  an absorbing Dirichlet boundary of a pair protective domain (left
  panel), or a reflecting boundary of a pair or single protective
  domain (right panel). Arrows are shown only where the jump rates
  differ from those given by Eq.~\eqref{eq:JumpRate} for uniform mesh
  cells.  \comment{The point labeled $x_0$ is the mesh point closest
    to the boundary (not necessarily the initial location of a
    molecule when the protective domain is first defined). When using
    a non-uniform mesh cell at an absorbing boundary, $h_\textrm{abs}$
    is defined to be the distance from the absorbing boundary to
    nearest mesh point of the uniform mesh, $x_0$.  Note that
    $h_{\textrm{abs}} \le h_p$, because otherwise the uniform mesh
    would have additional mesh points.}  \comment{Then,} the jump
  rates $a_{01}$ and $a_{02}$ \comment{in the absorbing boundary case}
  are given by the non-uniform rates in Eq.~\eqref{eq:NonUnifJumpRate}
  with $h_1 = h_p$ and $h_2 = h_{\textrm{abs}}$.  At a reflecting
  boundary, \comment{$h_{\textrm{ref}}$ is defined to be twice} the
  distance from the boundary to the nearest mesh point, $x_0$.
  \comment{Then,} $h_{\textrm{ref}}/2 \le h_p$ or $h_s$.  In this
  case, \comment{a molecule cannot jump} from $x_0$ toward the 
  boundary due to the reflecting \comment{zero-flux}  boundary 
  condition. The jump rate $a_{01}$, going away
  from the reflecting boundary, is given by
  Eq.~\eqref{eq:NonUnifJumpRate} with $h_1 = h_p \textrm{ or } h_s$,
  and $h_2 = h_{\textrm{ref}}$.  }
\end{figure} 

For single-molecule protective domains with one absorbing endpoint and
one reflecting endpoint, the mesh width $h_s$ is chosen to be the
largest value less than or equal to $h_s^{\textrm{max}}$ that exactly
divides the distance from the molecule to the absorbing endpoint.  A
mesh is defined so that the molecule and the absorbing endpoint lie
exactly on mesh points.  The mesh is uniform with the exception of one
non-uniform cell used immediately adjacent to the reflecting boundary,
as shown in Figure~\ref{fig:NonUnifBoundary} (right panel).

In pair protective domains, the mesh width $h_p$ is a
\comment{user-}specified value chosen to exactly divide
$r_\textrm{R}$.  Each time that two molecules are placed in a new pair
protective domain, the initial distance between the molecules will not
necessarily be divisible by $h_p$.  Rather than perturbing the
molecules, non-uniform mesh cells are used to move one of the
molecules, as shown in Figure~\ref{fig:NonUnifSubLattice}, so that a
uniform mesh of width $h_p$ can be defined with both molecules lying
exactly on mesh points.  Since this uniform mesh is chosen based on
the locations of the two molecules, the endpoints of the protective
domain may not conform with the mesh.  In this case, one non-uniform
mesh cell is used at each endpoint, which may be absorbing or
reflecting (see Figure~\ref{fig:NonUnifBoundary}).

Since $h_p$ is always chosen to exactly divide $r_{\textrm{R}}$, it
necessarily follows that $h_p \le r_{\textrm{R}}$.  For single
molecule protective domains we allow $h_s^{\textrm{max}}$ to be larger
than $r_{\textrm{R}}$.  In the convergence studies \comment{of the
  next section}, we set $h_s^{\textrm{max}} = k h_p$ where $k \ge1$,
and hold the ratio of $h_s^{\textrm{max}}$ to $h_p$ constant as both
are reduced to study convergence. As discussed \comment{above}, the
actual mesh widths $h_s$ used in single protective domains are almost
always strictly less than $h_s^{\textrm{max}}$, and non-uniform mesh
widths are used in both single and pair protective domains. For these
reasons, we keep track of the mean of the mesh widths that are
actually used in each simulation. In calculating this mean, each mesh
width is weighted by the number of times that it is actually used in a
sample path. Then, for all simulations performed with fixed
$h_s^{\textrm{max}}$ and $h_p$, we calculate an overall mean mesh
width by taking the arithmetic mean of the means for each simulation.


\section{Convergence of DL-FPKMC in One Dimension} \label{S:methConv}

\comment{In this section we perform convergence studies of DL-FPKMC}
for the annihilation reaction $A + B \to \varnothing$, where the
molecular species $A$ and $B$ undergo drift-diffusion subject to
various potentials on the interval $[0,1]$.  Our results demonstrate
both the convergence and accuracy of our method as the mesh widths,
$h_s^{\textrm{max}}$ and $h_p$, in the discretization are decreased.
We denote by $M^A(t)$ and $M^B(t)$ the number of molecules of $A$ and
$B$, respectively, at time $t$.  In the first set of convergence
studies \comment{(Subsection~\ref{S:2ParticleConvResults}),} only two
molecules are simulated, $M^A(0) = M^B(0) = 1$.  Each simulation runs
until the two molecules have reacted.  A large number ($10^7$) of
simulations are performed in order to minimize the statistical error,
so that the error due to the spatial discretization and the rate of
convergence can be studied.  In the next set of convergence studies
\comment{(Subsection~\ref{S:MultParticleConvResults}),} multiple
molecules each of $A$ and $B$ are simulated, $M^A(0) = M^B(0) = 10$ or
$M^A(0) = M^B(0) = 50$, and each simulation runs until all the
molecules have reacted. We will denote the $i^{\textrm{th}}$ molecule
of species $A$ by $A_i$, and the location of $A_i$ at time $t$ by
$Q_i^A(t)$. $B_j$ and $Q_j^B(t)$ are defined analogously. In the case
$M^A(0) = M^B(0) = 1$, we will drop the subscripts $i$ and $j$.

\subsubsection{\comment{Potential Functions and Parameters}}
The convergence studies are performed using three different potential
functions: (i) zero potential, $V_{\textrm{zero}}(x)=0$ (which results
in pure diffusion); (ii) a cosine potential with two energy wells,
$V_{\cos}(x)= \cos(4 \pi x)$; and (iii) a step potential with one
step,
\begin{equation*} V_{\textrm{step}}(x) = \left\{
    \begin{array}{rl} 2 & \text{if } x < \frac{1}{2}\\
      0 & \text{if } x \ge \frac{1}{2} \  . \end{array} \right.
\end{equation*}
The step potential is used to demonstrate that our DL-FPKMC algorithm
with the WPE discretization of the Fokker-Planck equation can handle
discontinuous potentials. Note that adding a constant to any potential
would not change the results, since the Fokker-Planck equation depends
on the derivative of the potential but not the potential itself. In
particular, any constant potential would produce the same results as
$V(x)=0$.

In all the convergence studies, the length $L$ of the overall domain
is $1$ unit, the boundaries of the overall domain are reflecting, and
the diffusion coefficient $D$ is $1$ unit$^2/$sec for both $A$ and
$B$. The values used for the reaction radius $r_\textrm{R}$ will be
specified in each subsection. We will use the notation
$\mathcal{U}(a,b)$ for the uniform random distribution on the interval
$(a,b)$.  The initial locations $Q_i^A(0)$ and $Q_j^B(0)$ are drawn
from $\, \mathcal{U}(a,b)$, where $(a,b) \subseteqq (0,L)$ will be
specified in each subsection. If $| Q_{i*}^A(0) - Q_{j*}^B(0)| \le
r_{\textrm{R}}$ for some $i^*$ and $j^*$, then $A_{i*}$ and $B_{j*}$
will react immediately, at $t=0$.  For $ t > 0 $, a reaction occurs if
$| Q_{i*}^A(t) - Q_{j*}^B(t)| = r_{\textrm{R}}$.

\comment{ Note, for most of the convergence studies we run in this
  section, these initial conditions correspond to non-equilibrium
  spatial distributions.  In the absence of reactions, the equilibrium
  probability density for a molecule to be at location $x \in (0, L)$
  is given by the Gibbs-Boltzmann distribution
\begin{equation} \label{eq:GibbsBoltzmann}
  \frac{e^{-V(x)}}{\int_0^{L} e^{-V(x)} \, dx}.
\end{equation}
For non-constant $V(x)$, the systems with uniformly distributed
initial locations on $(0,L)$ are not initially in spatial equilibrium,
since the molecule positions are not distributed according to
Eq.~\eqref{eq:GibbsBoltzmann}.  Furthermore, the systems with initial
locations restricted to subintervals of the overall domain are not in
spatial equilibrium for any of the potentials used here.}

\subsubsection{\comment{Comparison of DL-FPKMC simulation 
results in the two-molecule case to analytic and numerical solutions}} 
In general, to describe a stochastic reaction-drift-diffusion system
of many molecules by the probability density of having a given number
of molecules at specified positions, a large coupled system of partial
integro-differential equations is required~\cite{IsaacsonRDMENote}.
In the special case of only two molecules, $M^A(0) = M^B(0) = 1$, with
both molecules having the same diffusion coefficient $D$,
\comment{the SDLR model for} the reaction-drift-diffusion system $A + B
\to \varnothing$ in 1D can be described by a single 2D PDE: a
Fokker-Planck equation (or a diffusion equation when $V$ is constant).
Let $\rho(x,y,t)$ denote the \comment{joint} probability density for
finding molecule $A$ at location $x$ and molecule $B$ at location $y$
at time $t> 0$.  \comment{We consider the case where the molecules
  move in an interval domain of length $L$ with zero-flux boundary
  conditions at both endpoints.  The domain for the corresponding 2D
  PDE is then $\Omega = \{(x,y): 0 < x , y < L \ \text{and} \ |x-y| >
  r_{\textrm{R}} \}$, and} $\rho(x,y,t)$ evolves according to the
equations
\begin{align} \label{eq:FPEbvp}
\frac{\partial \rho(x,y,t)}{\partial t} = D \nabla \cdot  \left( \rho(x,y,t) 
\left( \frac{d V(x)}{d x}  + \frac{d V(y) }{d y} \right) + \nabla \rho(x,y,t) \right),   
&  \qquad \textrm{on } \comment{ \Omega},  \notag \\
\rho(x,y,t) = 0,  &  \qquad  \textrm{on }   
\comment{ \partial \Omega \cap \{ |x-y| = r_{\textrm{R}} \} },     \\
 \comment{ \rho (x,y,t)  \frac{\partial}{\partial \veta} (V(x)+V(y)) + }
 \frac{\partial \rho (x,y,t)}{\partial \veta } = 0,  & \qquad \textrm{on }
	  \comment{ \partial \Omega  \backslash  \{ |x-y| = r_{\textrm{R}} \} } , 
 \notag
\end{align}  
where $\veta = \veta(x,y)$ denotes the outward pointing normal at the
point $(x,y)$, \comment{and the gradient and divergence
operators are in the $(x, y)$ coordinates.}
The 2D domain for Eq.~\eqref{eq:FPEbvp} is illustrated
in Appendix~\ref{S:AnalytlcAndNumericalPDEsolns}, \comment{left panel
  of }Figure~\ref{fig:PDEdomains}.  When the initial locations of the
two molecules in the DL-FPKMC simulations are drawn from $\,
\mathcal{U}(0,L)$, the corresponding initial condition for the 2D
Fokker-Planck or diffusion equation is a constant, $\rho(x,y,0) =
1/L^2$.  Note, in the following we define $\rho$ on $0 \leq x, y \leq
L$ by also defining $\rho(x,y,t) = 0$ for $\abs{x-y} < r_{\textrm{R}}$
and $t > 0$. The 2D Fokker-Planck and diffusion
equations~\eqref{eq:FPEbvp} can both be solved numerically by finite
difference methods, and the diffusion equation can be solved
analytically using an eigenfunction expansion.  These numerical and
analytic solutions are discussed in
Appendix~\ref{S:AnalytlcAndNumericalPDEsolns}, and provide a baseline
against which we compare the results of two-molecule DL-FPKMC
simulations in Subsection~\ref{S:2ParticleConvResults}.

Let $T$ denote the random variable for the time at which the two
molecules react.  Using the solution $\rho(x,y,t)$ of
Eq.~\eqref{eq:FPEbvp}, we can calculate the survival probability,
\begin{align}  \label{eq:SurvProb}
S(t) & = \Pr \brac{ T > t } =  \int_0^L \int_0^L \rho(x,y,t) \, dx \, dy ,
\end{align}
and the mean reaction time,
\begin{equation}   \label{eq:MeanRxnTime}
\avg{T} = -\int_0^\infty t \, S'(t) \, dt = \int_0^\infty S(t) \, dt .
\end{equation} 
Note that $1-S(t)$ is the reaction time \comment{cumulative}
distribution function \comment{(CDF), and $-S'(t)$ is the
  corresponding density function.}

\subsubsection{\comment{Statistical Error and Discretization Error}}
In what follows, we will use the term `statistical error' to refer to
the difference between the empirical value of a statistic (e.g., mean
reaction time) estimated from the DL-FPKMC simulations and the upper,
or lower, bound of the $99\%$ confidence interval for the statistic.
By `discretization error,' we will mean the difference between the
empirical value from the DL-FPKMC simulations and the actual value.
In the two-molecule case, actual values are known exactly from the
analytic solution $\rho(x,y,t)$ \comment{of Eq.~\eqref{eq:FPEbvp}}
when $V$ is constant, and are estimated from the numerical PDE solver
described in Appendix~\ref{S:AnalytlcAndNumericalPDEsolns} when $V$ is
not constant. \comment{Note that the measured discretization error is
  the sum of two unknown quantities: sampling error, and the true
  discretization error due to the spatially discretized nature of the
  method.}

Since we perform a large number of simulations ($10^7$) at each mesh
width in the two-molecule case, the statistical error is quite small,
generally between $0.04\%$ and $0.19\%$ for 
\comment{reaction time} statistics.  Although we perform fewer
simulations ($4 \times 10^4$) when using multiple molecules each of
species $\textrm{A}$ and $\textrm{B}$, the statistical error is still
reasonably small, generally between $0.4\%$ and $1\%$.  Our results
show that as the mesh width is decreased, the discretization error
rapidly decreases to below the statistical error.  This demonstrates
that the DL-FPKMC algorithm converges and accurately resolves the
underlying reaction-drift-diffusion processes.


\subsection{Two-Molecule Convergence Studies}   \label{S:2ParticleConvResults}

In this subsection we consider the reaction $A + B \to \varnothing$
for a system with $M^A(0) = M^B(0) = 1$. Here, the initial locations
$Q^A(0)$ and $Q^B(0)$ are drawn randomly from $\, \mathcal{U}(0,1)$;
the reaction radius $r_{\textrm{R}}$ is $0.02$ units; and the pair
threshold $r_{\textrm{pair}}$ is equal to $2 r_{\textrm{R}}$ (i.e. the
molecules are placed in a pair protective domain when $| Q^A(t) -
Q^B(t) | \le r_{\textrm{pair}}$).  We study the convergence of
\comment{several statistics}
as the mesh widths $h_s^{\textrm{max}}$ and $h_p$ are reduced.

For $V=0$, the errors in the DL-FPKMC simulation results are
\comment{generally} calculated relative to the exact analytic
solution, as determined from the eigenfunction expansion in
Appendix~\ref{S:AnalytlcAndNumericalPDEsolns}.  For $V_{\cos}$ and
$V_{\textrm{step}}$, the errors are relative to the numerical solution
from the PDE solver described in
Appendix~\ref{S:AnalytlcAndNumericalPDEsolns}.  \comment{In the PDE
  solver, we use either the Crank--Nicolson method or the
  Twizel--Gumel--Arigu method for the time-stepping as explained in
  the appendix}.

As the mesh widths in the DL-FPKMC simulations are decreased,
\comment{the empirical values of the statistics studied approach the
  actual values} for all three potentials.  In each of
\comment{Figures~\ref{fig:MeanRxnTimes}, \ref{fig:RelErrMeanRxnTimes},
  and~\ref{fig:RelErrSurvProb}--\ref{fig:RelErrRxnLocCDF},} the
\comment{first} panel \comment{will show} results for $V=0$, the
\comment{second} panel for $V_{\cos}$, and the \comment{third} panel
for $V_{\textrm{step}}$.

\subsubsection{\comment{Mean Reaction Times} }
 
Let $\E_{\textrm{emp}}[T]$ be the empirical mean reaction time
calculated from the DL-FPKMC simulations. Define
$\E_{\textrm{upp}}[T]$ \comment{and $\E_{\textrm{low}}[T]$ to be,
  respectively, the upper and lower bounds} of the $99\%$ confidence
interval for the empirical mean reaction time. We denote by
$\E_{\textrm{act}}[T]$ the exact analytic mean reaction time in the
case that $V=0$, or the mean reaction time determined from the
numerical PDE solution in the $V \ne 0$ cases \comment{(see
  Eqs.~\eqref{eq:MeanRxnTime} and~\eqref{eq:AnalyticMeanRxnTime})}.
We calculate the relative error by
\begin{equation} \label{eq:RelErrMeanRxnTm}
\frac{ | \E_{\textrm{act}}[T] - \E_{\textrm{emp}}[T] | }{ \E_{\textrm{act}}[T] }  
\pm \frac{ \comment{ | \E_{\textrm{bd}}[T] - \E_{\textrm{emp}}[T] | }}{ \E_{\textrm{act}}[T] }
\comment{  \quad \text{ where } \E_{\textrm{bd}}[T] = \E_{\textrm{upp}}[T] \text{ or } \E_{\textrm{low}}[T].  }
\end{equation}
\comment{ The $99\%$ confidence interval for $\E_{\textrm{emp}}[T]$ is
  symmetric about the mean, so $\E_{\textrm{upp}}[T] -
  \E_{\textrm{emp}}[T] = \E_{\textrm{emp}}[T] - \E_{\textrm{low}}[T]
  $.  If $\E_{\textrm{act}}[T] $ is contained in the interval $(
  \E_{\textrm{low}}[T] , \E_{\textrm{upp}}[T] )$, then the
  discretization error $| \E_{\textrm{act}}[T] - \E_{\textrm{emp}}[T]
  | $ is less than the statistical error $| \E_{\textrm{bd}}[T] -
  \E_{\textrm{emp}}[T] |$.
}

\begin{figure}
\centering
\includegraphics{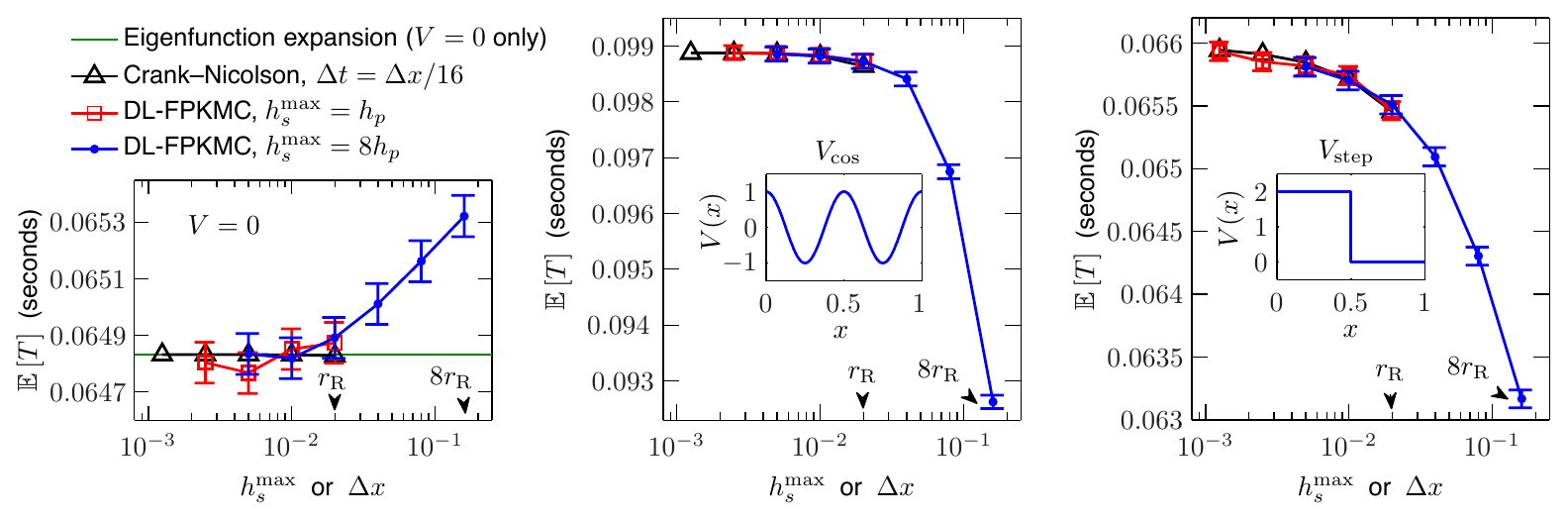}
\caption{\label{fig:MeanRxnTimes} \small Convergence of mean reaction
  time, $\avg{T}$, for the two-molecule $ A + B \to \varnothing $
  reaction, as the mesh width is decreased. Note the expanded scales
  of the vertical axes.  Each DL-FPKMC data point shows
  $\E_{\textrm{emp}}[T]$ from $10^7$ simulations, \comment{with}
  $99\%$ confidence intervals. 
  \comment{The resulting values of $\avg{T}$ from solving Eq.~\eqref{eq:FPEbvp} 
  using our Crank--Nicolson numerical PDE solver} are shown here in the 
  $V=0$ case to demonstrate the \comment{solver's accuracy (see
    Appendix~\ref{S:AnalytlcAndNumericalPDEsolns} for more detail)}.
  The mean reaction time calculated from the numerical PDE solution
  was resolved to an absolute error tolerance of \comment{less than} 
  $10^{-5}$ for $V_{\cos}$ and \comment{less than} $10^{-4}$ for $V_{\textrm{step}}$. 
  \comment{The PDE solver error tolerances are smaller than the 
    DL-FPKMC statistical errors, allowing the PDE solver estimate for 
    $\avg{T}$, determined with the finest value of $\Delta x$, to be used 
    in the absence of an analytic value when $V$ is non-constant. }}%
\end{figure}
\begin{figure}
\centering
\includegraphics{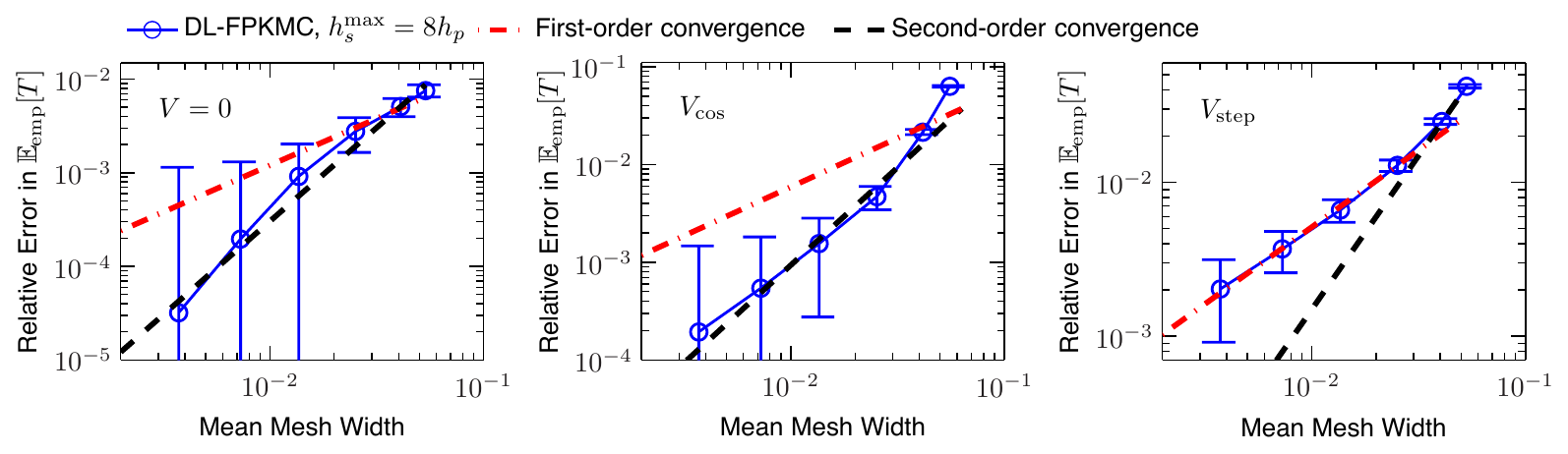}
\caption{\label{fig:RelErrMeanRxnTimes} \small Relative
  errors~\eqref{eq:RelErrMeanRxnTm} in the empirical mean reaction
  time, $\E_{\textrm{emp}}[T]$, for the two-molecule $A + B \to
  \varnothing$ reaction.  \comment{ In the $V=0$ case, the errors in
    $\E_{\textrm{emp}}[T]$ from the DL-FPKMC simulations are
    calculated relative to the exact analytic $\E_{\textrm{act}}[T]$
    given by Eq.~\eqref{eq:AnalyticMeanRxnTime}.  For $V \neq 0$, 
    the errors in $\E_{\textrm{emp}}[T]$ are calculated relative to the
    $\E_{\textrm{act}}[T]$ determined from the numerical solution
    of Eq.~\eqref{eq:FPEbvp} using the Crank--Nicolson method with 
    the finest value of $\Delta x$ (see
    Appendix~\ref{S:AnalytlcAndNumericalPDEsolns}).}  Note that the
  vertical axes have different scales in each panel.  \comment{Here 
    $h_s^{\textrm{max}}=8 h_p$, and the ratio is held constant as both 
    are decreased by successive halving.} Each data point is based on
    $10^7$ simulations with fixed \comment{values of} 
    $h_s^{\textrm{max}}$  and $h_p$. The error bars are determined by
    Eq.~\eqref{eq:RelErrMeanRxnTm} using the $99\%$ confidence
    intervals for the empirical mean reaction times.  \comment{ Note,
      the error bars are symmetric (and similar in size for the
      different mesh widths), but appear asymmetric (and larger for
      finer meshes) due to the log scale.  Relative errors are not
      plotted for the simulations with $h_s^{\textrm{max}}$ = $h_p \le
      r_{\textrm{R}}$, because the discretization errors are small
      compared to the statistical errors (cf.
      Fig.~\ref{fig:MeanRxnTimes}, red lines with square markers).}
  }%
\end{figure}
Figure~\ref{fig:MeanRxnTimes} \comment{shows} $\E_{\textrm{emp}}[T]$
plotted against $h_s^{\textrm{max}}$ as the mesh widths are varied,
\comment{and in Figure}~\ref{fig:RelErrMeanRxnTimes} the relative
errors in $\E_{\textrm{emp}}[T]$ \comment{calculated by
  Eq.~\eqref{eq:RelErrMeanRxnTm} are} plotted against the mean mesh
width.  The insets in Figure~\ref{fig:MeanRxnTimes} show the
respective potentials. \comment{As the mesh widths are reduced,} the
discretization errors decrease to less than the corresponding
statistical errors. Note that the statistical errors are very small
since $10^7$ simulations were performed at each mesh width.

We estimate the rate of convergence to be approximately second-order
for $V_{\cos}$ and approximately first-order for $V_{\textrm{step}}$.
This is consistent with the convergence rates of the WPE
discretization of the Fokker-Planck equation for smooth and
discontinuous potentials. In the $V=0$ case, it is difficult to draw a
conclusion about the rate of convergence since the discretization
errors are small relative to the statistical errors; however, for the
same reason, we can conclude that \comment{the results are} very
accurate in this case.

\subsubsection{\comment{Survival Probabilities} }

\begin{figure}
\centering
\includegraphics{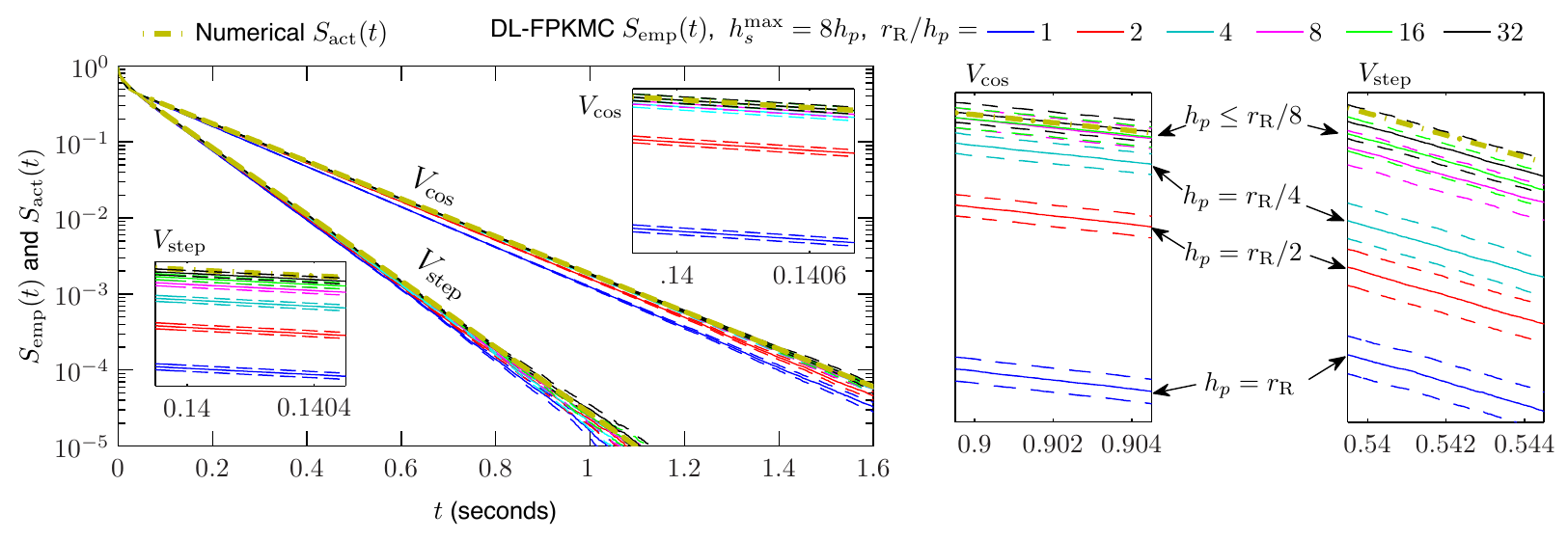}
\caption{\label{fig:SurvProb} \comment{ \small Convergence of
    empirical survival probabilities, $S_{\textrm{emp}}(t)$, for the
    two-molecule $ A + B \to \varnothing $ reaction with $V_{\cos}$
    and $V_{\textrm{step}}$.  Each empirical survival probability
    function is based on $10^7$ simulations. The dashed lines show 
    $99\%$ confidence bounds. (Note that $r_{\textrm{R}} = 0.02$ units
    is constant;  $r_{\textrm{R}}  / h_p$ changes as $h_p$ and 
    $h_s^{\textrm{max}}$ are successively halved.)
    $S_{\textrm{emp}}(t)$ from the $V=0$ case is not plotted because
    it is difficult to visually distinguish between the results for
    different mesh widths.  In the case $V=0$ and $h_p =
    r_{\textrm{R}}/32$, the analytic $S_{\textrm{act}}(t)$ is
    contained within the $99 \%$ confidence bounds of
    $S_{\textrm{emp}}(t)$ at \emph{every} time point $t_i$ at which
    $S_{\textrm{act}}(t)$ was evaluated ( $> 6 \times 10^4$ time
    points); for $h_p = r_{\textrm{R}}/8$ and $h_p =
    r_{\textrm{R}}/16$, $S_{\textrm{act}}(t)$ is contained within the
    confidence bounds of $S_{\textrm{emp}}(t)$ at more than $96 \%$ of
    the time points. This demonstrates that the discretization error
    is small compared to the statistical error.  }  }%
\end{figure}
For each mesh width, the empirical survival probability
\comment{$S_{\textrm{emp}}(t)$} and the associated $99\%$ confidence
bounds are calculated using the MATLAB function `\texttt{ecdf}'.
\comment{The confidence bounds determined by `\texttt{ecdf}'
  are symmetric to within numerical precision at all but a few points,
  which are at the tails of the distributions.}
$S_{\textrm{act}}(t)$ \comment{will denote} the analytic \comment{survival
probability for $V = 0$, and the} numerical survival
  probability \comment{determined from the PDE solver for $V \neq 0$ (see
  Eqs.~\eqref{eq:SurvProb}, \eqref{eq:AnalyticSurvProb} and
  Appendix~\ref{S:AnalytlcAndNumericalPDEsolns} for more
  information)}.  \comment{Figure~\ref{fig:SurvProb} plots the
  empirical survival probabilities $S_{\textrm{emp}}(t)$ with $99 \%$
  confidence bounds, and shows convergence of $S_{\textrm{emp}}(t)$ to
  $S_{\textrm{act}}(t)$ as the mesh widths are reduced.}

\comment{To quantify the magnitude of the error, we calculate the
  distances between $S_{\textrm{emp}}(t)$ and $S_{\textrm{act}}(t)$
  using the following methods.}  The discrete $L^1$, $L^2$, and
$L^\infty$ norms of a function $u(t)$ are given by \comment{
\begin{equation}  \label{eq:Norms}
 || u(t) ||_{L^1} = \sum_i | u(t_i) | \Delta{t_i},
 \qquad
  || u(t) ||_{L^2} = \left\{ \sum_i u(t_i)^2 \Delta{t_i}  \right\} ^{\frac{1}{2} },
   \qquad
    || u(t) ||_{L^\infty} = \max_i | u(t_i) | .
\end{equation}
} The relative error for $S_{\textrm{emp}}(t)$ we report in each norm
is then given by
\begin{equation} \label{eq:RelErrSurvProb}
\frac{ || S_{\textrm{act}}(t) - S_{\textrm{emp}}(t) || }{|| S_{\textrm{act}}(t) ||},
\end{equation}
where the norms are evaluated on the interval $t\in [0, \comment{\
  S^{-1} ( 10^{-6} ) } ] $ \comment{with} $\Delta{t_i} = t_{i+1} -
t_{i}$. The time points, $t_i$, used in evaluating the norms
correspond to those at which the numerical PDE solutions were
calculated
(see Appendix~\ref{S:AnalytlcAndNumericalPDEsolns} for more
information).  In the $V=0$ case, the analytic expression for the
survival probability was evaluated at those same $t_i$'s.  In all
cases, the empirical DL-FPKMC survival probabilities were linearly
interpolated to obtain values at every $t_i$.

\comment{ The $L^\infty$ error is equivalent to the Kolmogorov
  distance between distributions~\cite{CaoPetzoldJCP2006}, which is a
  statistical distance. Another statistical distance is the
  Kullback--Leibler (KL) divergence~\cite{Kullback1951information},
  which is a measure of the information lost when a distribution
  $G(t)$ is used to approximate a (true) distribution $F(t)$.  For
  continuous probability distributions, let $g(t)$ and $f(t)$ be the
  corresponding densities. Then, the KL divergence of $G(t)$ from
  $F(t)$ is given by
\begin{equation} \label{eq:KLcontinuous}
KL(F || G) = \int f(t) \ln \left( \frac{ f(t) }{ g(t) } \right) dt.
\end{equation}
For discrete distributions
\begin{equation} \label{eq:KLdiscrete}
KL(F || G) = \sum_i F_i \ln \left( \frac{ F_i }{G_i } \right).
\end{equation}
Since empirical densities obtained from Monte Carlo simulations are
noisy, we approximate $KL( S_{\textrm{act}} || S_{\textrm{emp}} )$ of
Eq.~\eqref{eq:KLcontinuous} by binning the data into a finite number
of time intervals and then applying Eq.~\eqref{eq:KLdiscrete}. We use
either $10$ or $20$ bins: $9$ or $19$ evenly-sized bins on the
interval $t\in \brac{0, S^{-1}(0.01) }$, and one bin for $t >
S^{-1}(0.01)$. This approximation of the KL divergence becomes noisier
as more bins are used. 
 Note, while the KL divergence can be interpreted as a measure of the
 difference between two distributions, it is does not define metric.
}

\begin{figure}
\centering
\includegraphics{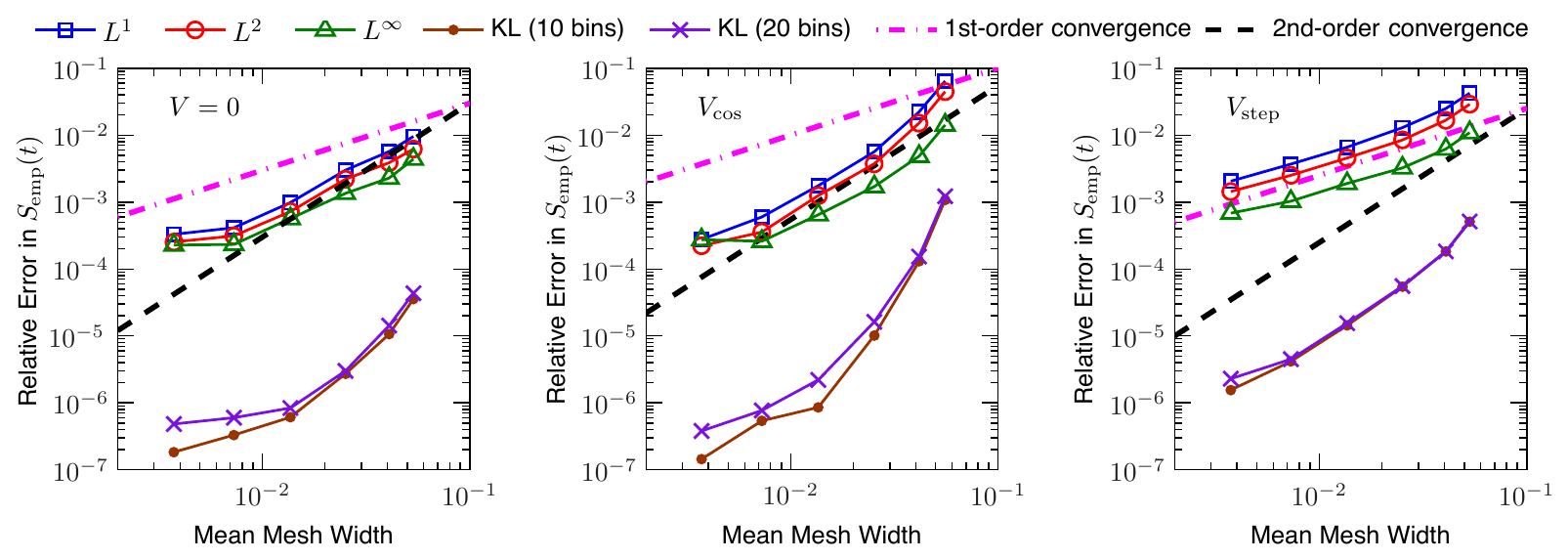}
\caption{\label{fig:RelErrSurvProb} \small Errors in the empirical
  survival probability $S_{\textrm{emp}}(t)$, \comment{relative to
    $S_{\textrm{act}}(t)$,} for the two-molecule $ A + B \to
  \varnothing $ reaction: \comment{relative errors measured in
    norm~\eqref{eq:RelErrSurvProb} and KL
    divergence~\eqref{eq:KLdiscrete}. Here} $h_s^{\textrm{max}} 
    = 8 h_p$ \comment{as both are decreased}. Each empirical 
  survival probability function is based on $10^7$ 
  simulations. \comment{Note 
    that the slight plateauing of the errors for finer
    meshes is due to sampling error, as illustrated in
    Fig.~\ref{fig:SurvProb} by the size of the confidence bounds
    relative to the smaller distance between $S_{\textrm{emp}}(t)$ and
    $S_{\textrm{act}}(t)$.}
}%
\end{figure}

\comment{Figure~\ref{fig:RelErrSurvProb} shows the relative errors
  measured in norm~\eqref{eq:RelErrSurvProb} and the KL
  divergence~\eqref{eq:KLdiscrete}, plotted against the mean mesh
  width.  The magnitude of the relative errors and the estimated
  convergence rates for $S_{\textrm{emp}}(t)$ are similar to those for
  $\E_{\textrm{emp}}[T]$. For each measure we see that the error in
  the survival probabilities decreases as the mesh widths are
  reduced. For finer mesh sizes, the discretization error becomes
  small compared to the statistical error (see Fig.~\ref{fig:SurvProb}).
  As such, the observed rate of convergence decreases and the graphs
  become less regular in Fig.~\ref{fig:RelErrSurvProb}.}

\subsubsection{\comment{Discrete Joint Spatial Probabilities}}

\comment{ Recall that the joint density, $\rho(x,y,t)$, is the
  solution of Eq.~\eqref{eq:FPEbvp}.
  We now study the convergence of the joint probability that $( Q^A(t)
  , Q^B(t) )$ is contained in one of $N_x \times N_x$ subregions of
  the domain $\Omega$ of Eq.~\eqref{eq:FPEbvp}. For $1 \le i, j \le
  N_x$, define the discrete joint probability function for the
  locations of the two molecules at time $t$ by
\begin{align} \label{eq:DiscreteJointSpatialProb}
p_{\textrm{act}}(i,j,t) &= \Pr \brac{ \frac{ i-1}{N_x } L < Q^A(t) < \frac{i}{N_x } L, \ \frac{j-1}{N_x } L < Q^B(t) < \frac{j}{N_x } L   } \\
&= \int_{(j-1)L/N_x  }^{ j L/N_x }   \int_{(i-1)L/N_x  }^{ i L/N_x }   \rho(x,y,t) \, dx \, dy. \notag
\end{align}
We study the convergence of the empirical $p_{\textrm{emp}}(i,j,t)$
from the DL-FPKMC simulations at two fixed times, $t =
S_\textrm{act}^{-1}(0.75)$ and $t = S_\textrm{act}^{-1}(0.5)$, where
$p_{\textrm{emp}}(i,j,t)$ is determined by binning the locations of
the molecules at time $t$ into the $N_x \times N_x$ subregions.

\begin{figure}
\centering
\includegraphics{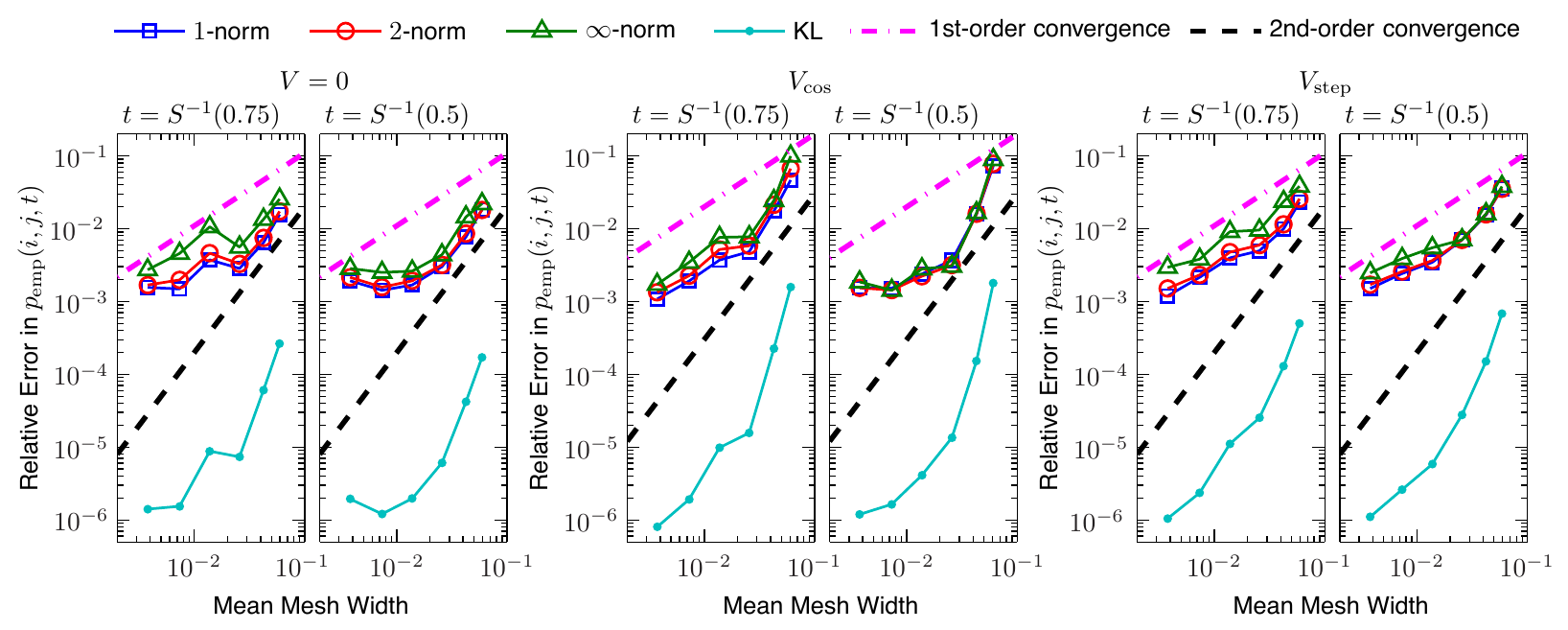}
\caption{\label{fig:RelErrLocs} \small \comment{ Errors in the
    empirical joint spatial probabilities $p_\textrm{emp}(i,j,t)$ at
    $t = S_\textrm{act}^{-1}(0.75)$ and $t = S_\textrm{act}^{-1}(0.5)$
    for the two-molecule $ A + B \to \varnothing $ reaction.  
    For all graphs we divide the domain into $25$ equally-sized 
    bins ($1 \leq i, j \leq 5$). The
    errors are calculated relative to $p_\textrm{act}(i,j,t)$, as
    determined from the PDE solver 
    with Twizell--Gumel--Arigu time-stepping described in
    Appendix~\ref{S:AnalytlcAndNumericalPDEsolns}.   
    Here $h_s^{\textrm{max}} = 8 h_p$ as both are decreased.
    Each data point is based on $10^7$ simulations.  At $t =
    S_\textrm{act}^{-1}(0.75)$ (resp. $t = S_\textrm{act}^{-1}(0.5)$)
    the two molecules have not yet reacted in approximately $75 \%$
    (resp. $50 \%$) of the simulations.  We also recorded
    $p_\textrm{emp}(i,j,t)$ at $t = S_\textrm{act}^{-1}(0.25)$,
    however the convergence results were more strongly affected by noise, since
    the molecules had already reacted by that time in most of the
    simulations.  }}%
\end{figure}
We calculate the relative error in $p_{\textrm{emp}}(i,j,t)$ at a
fixed time $t$ by
\begin{equation} \label{eq:RelErrJointProb}
  \frac{ || p_{\textrm{act}}(i,j,t) - p_{\textrm{emp}}(i,j,t) || }{|| p_{\textrm{act}}(i,j,t) ||},
\end{equation}
where the norms in the numerator and denominator are given by
\begin{align}  \label{eq:2Dnorms}
  || q(i,j,t)  ||_1 &= \sum_{i=1}^{N_x}  \sum_{j=1}^{N_x} | q(i,j,t) | ,
  \qquad
  || q(i,j,t)  ||_2 = \left\{ \sum_{i=1}^{N_x}  \sum_{j=1}^{N_x} q(i,j,t)^2 \right\} ^{\frac{1}{2} }, \\
  & \text{ or } \quad
  || q(i,j,t)  ||_{\infty} = \max_{1 \le i \le N_x}  \max_{1 \le j \le N_x}| q(i,j,t) | , \notag
\end{align}
with $q = p_{\textrm{act}}- p_{\textrm{emp}}$ for the numerator
and $q = p_{\textrm{act}}$ for the denominator.

We also calculate the KL divergence~\eqref{eq:KLdiscrete} of
$p_{\textrm{emp}}(i,j,t)$ from $p_{\textrm{act}}(i,j,t)$.  Since the
KL divergence is only defined for probabilities that sum to one, we
include an extra ``bin" representing that probability that the two
molecules have already reacted, which is given by $ 1- S(t) = 1 -
\sum_{i=1}^{N_x} \sum_{j=1}^{N_x} p(i,j,t)$.  The relative errors and
KL divergence are shown in Fig.~\ref{fig:RelErrLocs}. We see similar
convergence behavior as in the previous subsection, with the errors
clearly decreasing for coarser mesh widths, but plateauing as sampling
error becomes dominant at smaller mesh widths.}

\subsubsection{\comment{Reaction Location Distributions} }

\comment{ For a reaction between $A$ and $B$, the location of the
  reaction will be taken to be $Q^{\textrm{rxn}} = ( Q^A(t) + Q^B(t)
  )/2$ where $t$ is the time of the reaction. Then $
  F^{\textrm{rxn}}(x) = \Pr \brac{ Q^{\textrm{rxn}} \le x }$ is the
  cumulative distribution function for the reaction locations.  In
  studying convergence of the empirical
  $F_{\textrm{emp}}^{\textrm{rxn}}(x)$ as the mesh widths are reduced,
  we consider only those reactions that occur after $t=0$.  (When $|
  Q^A(0) - Q^B(0)| \le r_{\textrm{R}}$, the reaction occurs
  immediately at $t=0$ regardless of the mesh width.) For the
  parameters used here approximately $3.96 \%$ of the reactions occur
  at $t=0$.

  For $V=0$, the exact distribution $
  F_{\textrm{act}}^{\textrm{rxn}}(x)$ can be calculated analytically
  from the solution of Eq.~\eqref{eq:FPEbvp}.  See
  Eq.~\eqref{eq:AnalyticRxnLocCDF} of
  Appendix~\ref{S:AnalytlcAndNumericalPDEsolns} for the analytic
  result.  In this case, we study convergence by comparing
  $F_{\textrm{emp}}^{\textrm{rxn}}(x)$ from the DL-FPKMC simulations
  to the analytic $ F_{\textrm{act}}^{\textrm{rxn}}(x)$. For the $V
  \ne 0$ cases, we examine the successive pairwise differences between
  $F_{\textrm{emp}}^{\textrm{rxn}}(x)$ from the DL-FPKMC simulations
  as the mesh widths are decreased.  In all cases, the distributions
  $F^{\textrm{rxn}}(x)$ are evaluated at $10^4$ evenly spaced points.

\begin{figure}
\centering
\includegraphics{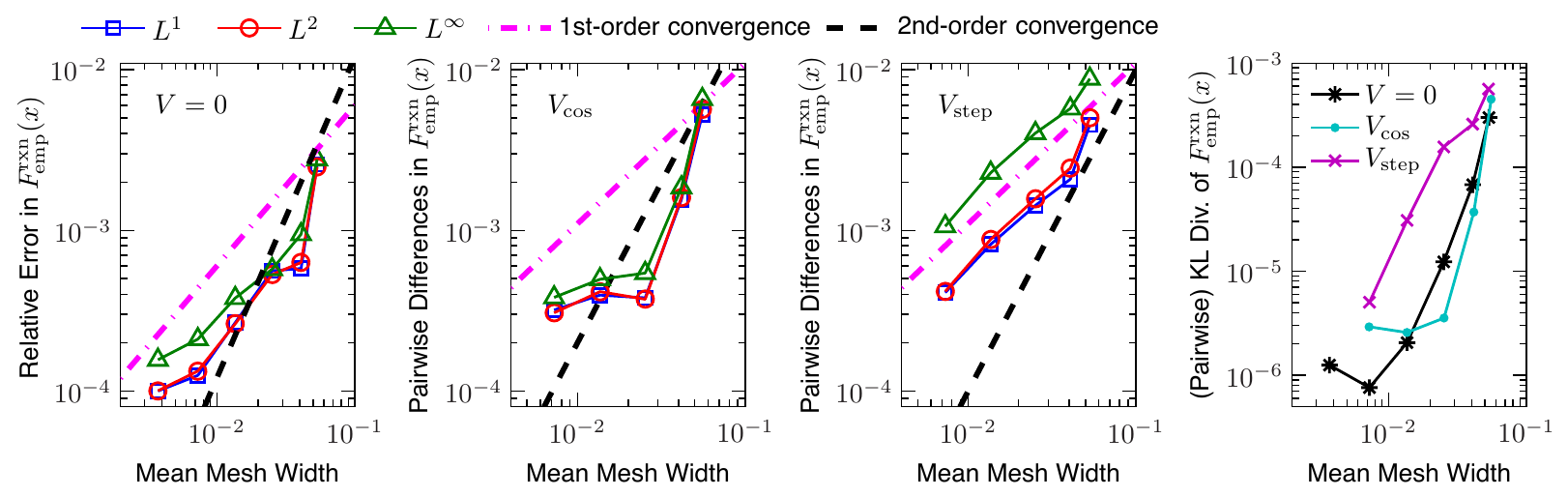}
\caption{\label{fig:RelErrRxnLocCDF} \small \comment{Errors in
    $F_{\textrm{emp}}^{\textrm{rxn}}(x)$, the empirical reaction
    location CDF, for the two-molecule $ A + B \to \varnothing $
    reaction: errors relative to the analytic
    $F_{\textrm{act}}^{\textrm{rxn}}(x)$ for $V=0$, and pairwise
    relative differences between $F_{\textrm{emp}}^{\textrm{rxn}}(x)$
    from successive mesh widths for $V_{\cos}$ and
    $V_{\textrm{step}}$. The pairwise relative differences are plotted
    against the mean mesh width for the coarser of the two meshes.
    Here $h_s^{\textrm{max}} = 8 h_p$ as both are decreased.
    Each empirical CDF is based on the simulations in which the two
    molecules react after time $t=0$, out of $10^7$ total simulations.
    The plateaus in the errors or differences at the finer mesh widths 
    are due to sampling error.  For example, in the case of $V_{\cos}$,
    the $99 \%$ confidence bounds for
    $F_{\textrm{emp}}^{\textrm{rxn}}(x)$ from the two finest mesh
    widths overlap at more than $96 \%$ of the $10^4$ points at which
    they were evaluated.  }}%
\end{figure}

  For $V=0$, the relative errors in each norm~\eqref{eq:Norms} are
  calculated by $|| F_{\textrm{act}}^{\textrm{rxn}}(x) -
  F_{\textrm{emp}}^{\textrm{rxn}}(x) || \, / \, ||
  F_{\textrm{act}}^{\textrm{rxn}}(x) ||$.  For $V_{\cos}$ and
  $V_{\textrm{step}}$, the pairwise relative differences are
  calculated by the preceding formula with
  $F_{\textrm{emp}}^{\textrm{rxn}}(x)$ from the finer of the two
  meshes replacing $F_{\textrm{act}}^{\textrm{rxn}}(x)$. 
  The corresponding density is approximated by binning the
  reaction locations into 20 evenly-sized bins, so that the KL
  divergence~\eqref{eq:KLdiscrete} can be calculated
  (comparing to the analytic values for $V=0$, and
  comparing the empirical values from successive
  mesh widths for $V_{\cos}$ and $V_{\textrm{step}}$).  The results are
  shown in Fig.~\ref{fig:RelErrRxnLocCDF}.  Note that the relative errors
  and pairwise relative differences are all less than $1\%$.
  Although the results for the spatial statistics are more affected
  by noise (Figs.~\ref{fig:RelErrLocs},~\ref{fig:RelErrRxnLocCDF}), 
  the convergence rates appear similar to those for the time statistics 
  (cf. Figs.~\ref{fig:RelErrMeanRxnTimes}, ~\ref{fig:RelErrSurvProb}).
  }


\subsection{Results of  Multiple-Molecule Convergence Studies} \label{S:MultParticleConvResults}

In the convergence studies for multiple molecules undergoing the
reaction $A + B \to \varnothing$, we start the simulations with
either $20$ molecules ($M^A(0) = M^B(0) = 10$) or $100$ molecules
($M^A(0) = M^B(0) = 50$).  In the $20$-molecule simulations:
$r_{\textrm{R}}$ = $0.02$ units; $r_{\textrm{pair}} =
2r_{\textrm{R}}$; $Q_{i}^A(0) \sim \mathcal{U}(0.1,0.4)$ and
$Q_{j}^B(0) \sim \mathcal{U}(0.6,0.9)$ for $ 1 \le i, \, j \le 10$.  In
the $100$-molecule simulations: $r_{\textrm{R}}$ = $0.001$ units;
$r_{\textrm{pair}} = 4r_{\textrm{R}}$; $Q_{i}^A(0)$ and $Q_{j}^B(0)
\sim \mathcal{U}(0,1)$ for $ 1 \le i, \, j \le 50$.  Two molecules
$A_{i*}$ and $B_{j*}$ are placed in a pair protective domain if 
they are closer to each other than to any other molecules of the 
opposite type, i.e. if
\begin{equation*} 
| Q_{i*}^A(t) - Q_{j*}^B(t) |  = \min_i | Q_{i}^A(t) - Q_{j*}^B(t) |  = \min_j | Q_{i*}^A(t) - Q_{j}^B(t) |,
\end{equation*}
and if the distance between them satisfies
\begin{align} \label{eq:PairConditions}
| Q_{i*}^A(t) - Q_{j*}^B(t) |  &\le r_{\textrm{pair}}  \notag \\
\text{ and}&  \notag \\
 | Q_{i*}^A(t) - Q_{j*}^B(t) |  &\le r_{\textrm{R}} + \min(  \min_{i \ne i*}( | Q_{i*}^A(t) - Q_{i}^A(t) | , \ \min_{j \ne j*} | Q_{j*}^B(t) - Q_{j}^B(t) | ).
\end{align}
The last condition in Eq.~\eqref{eq:PairConditions} was added to prevent 
the length of the
protective domain of a molecule, say $A_i$, from 
approaching zero
when $A_i$ is close to another molecule of the same type, say
$A_{i*}$, where $A_{i*}$ would otherwise have been placed in a pair
with $B_{j*}$.

\begin{figure}
\centering
\includegraphics[width=15cm]{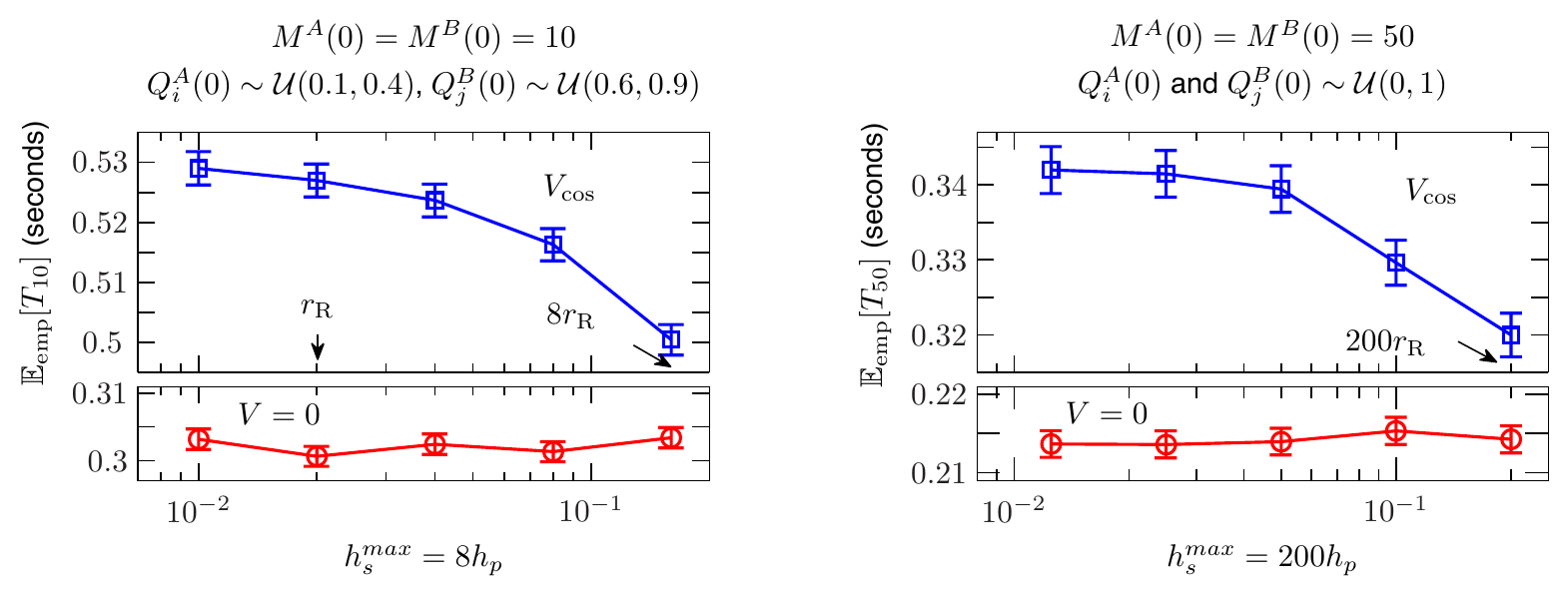}
\caption{\label{fig:MeanTimeLastRxn} \small 
 Mean time for all molecules to react via the reaction 
 $A + B \to \varnothing$. Error bars show 99\% confidence 
 intervals, based on $4 \times 10^4$ simulations per data point.
 \emph{Left panel}: $r_{\textrm{R}} = 0.02$ units.
 \emph{Right panel}:  $r_{\textrm{R}} = 0.001$ units. }
\end{figure}

\begin{figure}
\centering
\includegraphics{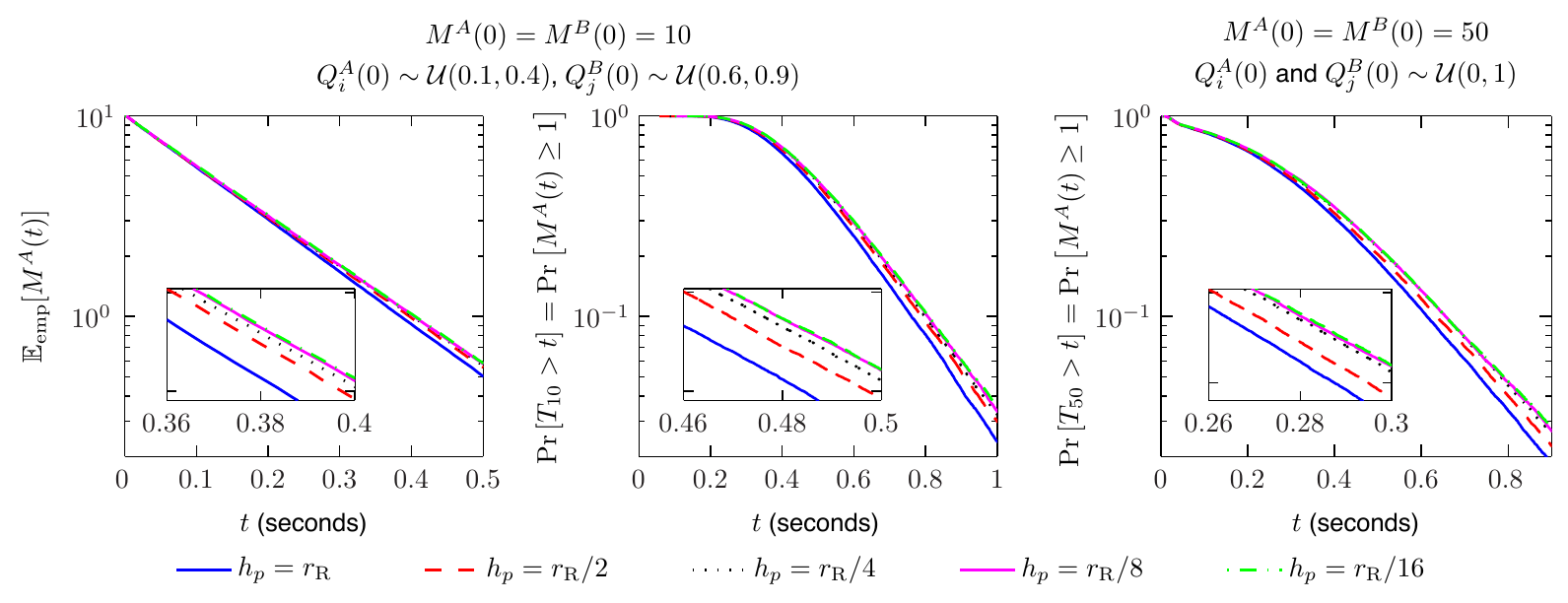}
\caption{\label{fig:NmoleculeDistributionsVcos} \small 
 Mean number of molecules of $A$ remaining at time $t$,
 $\E_{\textrm{emp}}[M^A(t)]$ (left panel), and probability that 
 at least one molecule of $A$ remains at time $t$, 
 $\Pr \brac{ M^A(t) \ge 1}$ (center and right panels).
 $V= \cos(4 \pi x)$. Each graph is based on $4 \times 10^4$ 
 simulations. $\E_{\textrm{emp}}[M^A(t)]$ is not plotted in the 
 case $M^A(0) = M^B(0) = 50$, because the results for the 
 different mesh widths are essentially indistinguishable until  
 $\E_{\textrm{emp}}[M^A(t)] \lesssim 4$.
 \emph{Left and center}: $r_{\textrm{R}} = 0.02$ units, 
 $h_s^{max} = 8 h_p$.
 \emph{Right}:  $r_{\textrm{R}} = 0.001$ units,
 $h_s^{max} = 200 h_p$.  }%
\end{figure}

Each simulation runs until all of the molecules have reacted.  Let
$T_n$ denote the random variable for the time at which the
$n^{\textrm{th}}$ reaction occurs. Since we are working with the
irreversible reaction $A + B \to \varnothing$, we have that $M^A(t) =
M^A(0) - n$ for $t \in [T_n, T_{n+1})$, and similarly for $M^B(t)$.
Figure~\ref{fig:MeanTimeLastRxn} shows the empirical mean time for all
molecules to react; this is $\E_{\textrm{emp}}[T_{10}]$ when $M^A(0) =
M^B(0) = 10$ and $\E_{\textrm{emp}}[T_{50}]$ when $M^A(0) = M^B(0) =
50$. With these values of $M^A(0)$ and $M^B(0)$, the SDLR model
described in Section~\ref{S:overview} will correspond to a large
coupled system of partial integro-differential equations for the
probability densities of having a specified number of molecules at
specified locations.  It is no longer feasible to solve these
equations numerically to obtain high-accuracy solutions for assessing
the empirical convergence of our DL-FPKMC method.
As such, we now estimate the accuracy of reaction time statistics by
comparing DL-FPKMC results at coarser mesh widths to results obtained
with the finest mesh width.

When $V(x)= \cos(4 \pi x)$ (Fig.~\ref{fig:MeanTimeLastRxn}, 
top panels), we see convergence as
the mesh width is decreased.  The percent difference between
$\E_{\textrm{emp}}[T_{10}]$ for the coarsest mesh width compared 
to the finest mesh width is approximately $5.4\%$; this percent 
difference for $\E_{\textrm{emp}}[T_{50}]$ is approximately $7\%$.  
These differences are comparable in size to the explicitly calculated 
discretization error of approximately $6.4\%$ at the coarsest mesh 
width \comment{for $\E_{\textrm{emp}}[T]$} in the two-molecule case (previous subsection).

When $V=0$ (Fig.~\ref{fig:MeanTimeLastRxn}, 
bottom panels), the confidence intervals for all mesh widths
overlap, indicating that the discretization error is less than the
statistical error even for the coarsest mesh width. The statistical
errors when $V=0$ are between $0.49\%$ and $0.81\%$.  If the 
unknown discretization error here is comparable in size to the 
known discretization error \comment{in $\E_{\textrm{emp}}[T]$} 
in the two-molecule case, which was 
approximately $0.87\%$ at the coarsest mesh width, that would 
provide an explanation for why can not observe convergence 
when $V=0$.

$\E_{\textrm{emp}}[M^A(t)]$ is the mean number of $A$ molecules 
remaining  at time $t$, and $\Pr \brac{ M^A(t) \ge 1}$ is the probability 
that at least one molecule of $A$ remains at time $t$. 
Figure~\ref{fig:NmoleculeDistributionsVcos} shows convergence of 
$\E_{\textrm{emp}}[M^A(t)]$ and  $\Pr \brac{ M^A(t) \ge 1} $ as the 
mesh width is decreased, in the case $V(x)= \cos(4 \pi x)$. Results are not 
plotted for the $V=0$ case, because the confidence bounds for different 
mesh widths overlap; this indicates that the results are resolved to within
the statistical error, even for coarser mesh widths.


\section{Running Time Analysis}  \label{S:RunningTime}

In this section we demonstrate that the running time of the
  DL-FPKMC algorithm when simulating the $A + B \to \varnothing$
  reaction scales linearly with the number of molecules in the system.
  We also compare the computational performance of DL-FPKMC to a
  second method in which all molecules hop on a fixed, global, uniform
  lattice. Both methods were implemented in MATLAB, and an attempt was
  made to take advantage of reasonable and standard optimizations.
  That said, we make no claim that our implementation of either method
  provides optimal computational performance. All DL-FPKMC and fixed lattice
  simulations were performed in MATLAB on a Sun Fire X4600 M2 x64
  server. The server was configured with four AMD Opteron Model 8220
  processors (2.8 GHz dual-core) and 16 GB of RAM.

  In both methods, bimolecular reactions occur when the distance
  between two reactants is exactly equal to $r_\textrm{R}$.  To
  enforce this condition in the fixed lattice method, the lattice
  spacing $h$ will be chosen to equal $r_\textrm{R}$. We set
  $V(x)=0$, so the  spatial hopping rates in the fixed lattice method are
  simply $D/h^2$.  The results of the convergence studies in
  Section~\ref{S:methConv} indicate that in DL-FPKMC taking $h_p
  \approx r_\textrm{R}$ and $h_s^{\textrm{max}} \approx L/50$ is
  sufficient to resolve the reaction and diffusion processes to within
  statistical error for $n \le 10^7$ simulations.  For biologically
  relevant parameter values, e.g. $L = 10 \mu$m and $r_\textrm{R} =
  1$nm, this will allow $h_s^{\textrm{max}}$ in DL-FPKMC to be a couple
  orders of magnitude larger than $h$ in the fixed lattice method
  without compromising accuracy.\footnote{Since
    $V(x)=0$, in the fixed lattice method each possible spatial hop is
    equally likely.  This allows optimization of the underlying
    Stochastic Simulation
    Algorithm~\cite{GillespieJPCHEM1977,KalosKMC75}.  Incorporating
    drift into the fixed lattice method would give spatially-varying
    hopping rates, so additional computational cost might be required
    to sample which molecule hops and in which direction.}
 
  To our knowledge, the only reaction system that has previously been
  simulated \emph{in one dimension} using any FPKMC-type method is the
  $A + A \to \varnothing$~\cite{KalosDSMC06,OppelstrupPRE2009}
  reaction, in which any two molecules annihilate as soon as they
  collide. In contrast, for the reaction system $A + B \to
  \varnothing$, molecules of the same type do not react.  As mentioned
  in Subsection~\ref{StepsOfAlgorithm}, we allow the protective domains
  of non-reacting molecules to overlap to prevent the size of
  protective domains from going to zero when two non-reacting
  molecules approach each other. We expect that this issue could also
  be addressed by allowing protective domains to contain more than two
  molecules.  In order to simulate the same underlying
  process with both DL-FPKMC and the fixed lattice method,
  non-reacting molecules in the fixed lattice simulations are allowed
  to occupy the same lattice site and cross each other.
 
\subsection{\comment{Parameters}}
 
In comparing the DL-FPKMC and fixed lattice methods, the following
parameters values are used: $L = 10 \mu \text{m}$, $r_\textrm{R} =
1\text{nm}$, $D = 10 \mu \text{m}^2/ \text{sec}$, and $V(x) = 0$.  The
overall simulation domain, the interval $(0, L)$, has reflecting
boundaries.  The initial number of molecules of $A$ and $B$ are equal,
$M^A(0) = M^B(0)$, and each simulation runs until all molecules have
reacted. The initial locations of molecules are uniformly distributed
over the interval $(0, L)$. In both DL-FPKMC and the fixed lattice
method, if the initial distance between a molecule of $A$ and a
molecule of $B$ is less than or equal to $r_\textrm {R}$, then they
react immediately. All later reactions occur when the distance between
an $A$ molecule and a $B$ molecule equals $r_\textrm{R}$.  In the
fixed lattice simulations, $h = r_\textrm{R} = 1\text{nm}$.  In
DL-FPKMC, $h_s^\textrm{max} = L/50 = 200 \text{nm}$, $h_p =
r_\textrm{R} = 1\text{nm}$, and $r_\textrm{pair} = 50 \text{nm}$.  To
check that using a coarser value for $h_s^\textrm{max}$ than $h_p$ is
still sufficient to obtain accurate results, we ran $10^6$ simulations
with $M^A(0) = M^B(0) =1$. In this case the exact mean reaction time,
$\E_{\textrm{act}}[T]$, is known analytically, see
Eq.~\eqref{eq:AnalyticMeanRxnTime}. Using the preceding parameters,
the resulting $\E_{\textrm{emp}}[T]$ from the DL-FPKMC simulations
agrees with $\E_{\textrm{act}}[T]$ to within statistical error, which
is approximately $0.34\%$.

\subsection{\comment{Results of Running Time Analysis}}

\begin{figure}
\centering
\includegraphics{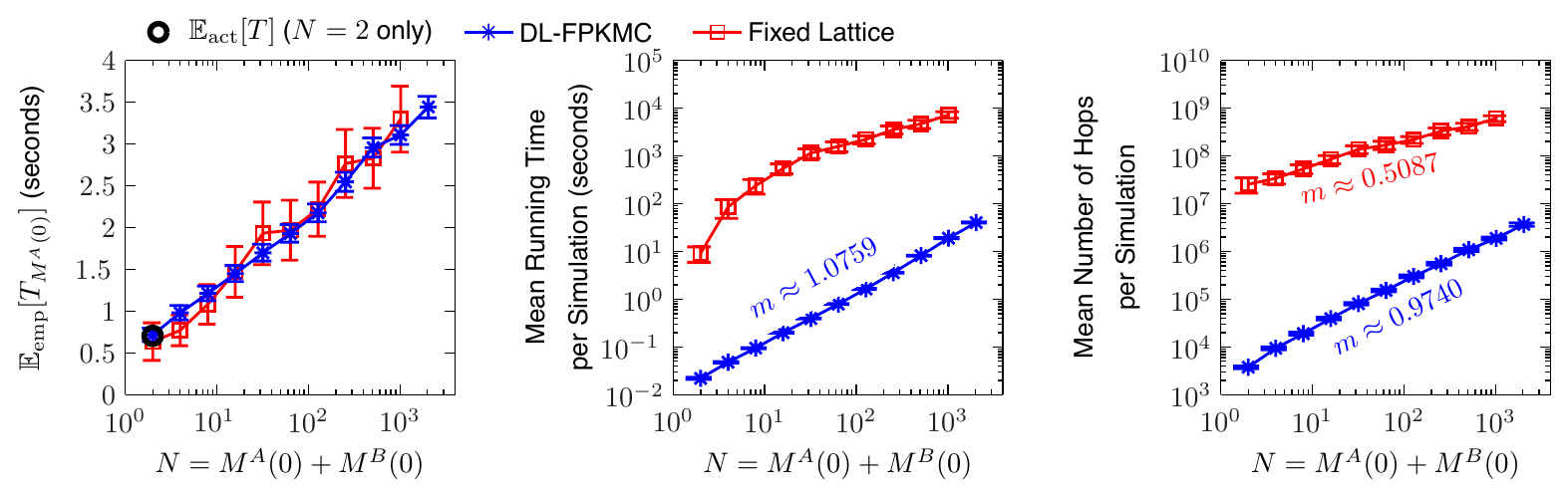}
\caption{\label{fig:RunningTime} \small 
Comparison of DL-FPKMC and fixed lattice methods for the reaction 
 $A + B \to \varnothing $ as the number of molecules present initially,
 $N = M^A(0) + M^B(0)$, is increased. 
 $M^A(0) = M^B(0)$, $Q_{i}^A(0)$ and $Q_{j}^B(0)
\sim \mathcal{U}(0,L)$ for $ 1 \le i, \, j \le M^A(0)$,
 $L = 10 \mu \text{m}$, $r_\textrm{R} = 1\text{nm}$, 
 $D = 10 \mu \text{m}^2/ \text{sec}$, and $V(x) = 0$.
 Each DL-FPKMC data point is based on $10^3$ simulations,
 and each fixed lattice data point is based on $10^2$ simulations.
 Error bars indicate $99 \%$ confidence intervals. 
}%
\end{figure}

Figure~\ref{fig:RunningTime} compares the \comment{simulation} 
results and computational
performance of DL-FPKMC to those of the fixed lattice method as the
total number of molecules in the system, $N = M^A(0) + M^B(0)$, is
varied. As shown in the left panel, the resulting values of
$\E_{\textrm{emp}}[T_{M^A(0)}]$, the mean time for all molecules to
react within the simulations, from the two methods agree to within
statistical error.  The running time (center panel) for DL-FPKMC is
two to three orders of magnitude faster than the running time for the
fixed lattice.  DL-FPKMC also requires fewer hops than the fixed
lattice by two to three orders of magnitude (right panel).  The slope
of the line for DL-FPKMC on the log-log plot of running time (center
panel) indicates that the method is approximately $O(N)$, i.e., the
running time scales linearly with the number of molecules, $N = M^A(0)
+ M^B(0)$.  Although the fixed lattice method has better scaling with
$N$, it would not appear to become more efficient than DL-FPKMC until
$N$ is substantially larger than $10^5$. At that large a number of
molecules, it is common to transition to more macroscopic stochastic
reaction-diffusion models.

As the number of molecules in the system is increased, the proportion of
running time spent on different steps of the DL-FPKMC algorithm
changes.  For the parameter values used in this section, the changes
in proportions of running time as $N$ is increased from $2$ to $2048$ are
as follows (approximated using the MATLAB `\texttt{profile}'
function):
\begin{itemize}  [itemsep=0pt, topsep=1pt]
\item Generating sample paths: decreases from  $67 \%$ to $38 \%$,
\item Identifying neighboring molecules and defining new protective domains 
for updated molecules based on locations of neighbors:
increases from  $ 17\%$ to $ 52 \%$,
\item Determining which molecules (possibly none) to no-passage update
after each first-passage update:
decreases from  $ 4\%$ to $ 3\%$,
\item Sorting the event queue: decreases from  $7 \%$ to $3 \%$.
\end{itemize}
These percentages could vary substantially for different parameter
values and mesh widths, but we expect the overall trends would hold.

As was described in Subsection~\ref{S:SamplePaths}, the DL-FPKMC
algorithm generates a sample path for each molecule within its
protective domain until a first-passage or reaction event occurs.
Whenever a molecule is no-passage updated before its next event time,
part of its sample path goes unused.  The mean number of hops per
simulation shown in Figure~\ref{fig:RunningTime} (right panel) is
based on the total number of hops in all sample paths that are
generated, not just the hops that are actually used.  Particularly
long paths tend to be generated for protective domains touching an
overall domain boundary, since molecules can only exit through one
endpoint of such protective domains.  If we cap the length of single-molecule
protective domains that are near the boundaries, then the total number of
hops decreases, but the number of times that the protective domains
are updated increases and the overall running time also increases.
Some guidelines for and difficulties in optimizing the partitioning of
space among protective domains in FPKMC are discussed
in~\cite{DonevJCP2010}, but we have not yet attempted to address this
difficult optimization problem.

\subsection{\comment{Expectations for Future Studies}}

It is expected that the implementation of DL-FPKMC in higher
dimensions can be done in such a way as to maintain the $O(N)$
scaling for the following reasons:
\begin{itemize} [itemsep=0pt, topsep=1pt]
\item FPKMC in two and three dimensions can be implemented to have
  $O(N)$ scaling~\cite{DonevJCP2010}.  One might expect identifying
  neighboring molecules to be more costly in higher dimensions,
  \comment{however, in practice the near-neighbor list (NNL) method
    allows for constant per event costs, leading to $O(N)$ scaling
    for the overall method~\cite{DonevJCP2010};}
\item The only algorithmic difference between DL-FPKMC and FPKMC is
  that DL-FPKMC generates sample paths within protective domains using
  continuous-time random walks, whereas FPKMC samples from analytic
  solutions of the diffusion equation;
\item Since the fraction of running time spent generating sample paths
  in DL-FPKMC decreases as $N$ increases, the scaling of overall
  running time as $N$ increases is expected to depend mainly on other
  steps of the algorithm. These steps can all be implemented in the
  same way in DL-FPKMC as in FPKMC.
\end{itemize}

\comment{See Table 5 of~\cite{DobrzynskiBioinf07} for the computational 
cost of several reaction-diffusion methods, including methods for simulating the SDLR
model and the RDME.  The original, approximate GFRD~\cite{WoldeGFRD05} 
and Smoldyn~\cite{AndrewsBrayPhysBio2004} exhibit 
$O(N)$ scaling with the total number of molecules for diffusive movements
and $O(\sum_{N_R} \prod_{S \in R} N_S )$ scaling for reactive distances, where
$N_R$ is the number of reaction channels and $N_S$ is the number of molecules
of a given species~\cite{DobrzynskiBioinf07}. We expect FPKMC and DL-FPKMC 
will have similar scaling with number of reaction channels.}

\comment{ Although $O(N)$ scaling would eventually cause the
  computational cost to become too high for very large $N$, many
  systems of biological interest could still be simulated at low cost.
  FPKMC has been used to simulate systems with $10^8$ particles in one
  dimension~\cite{KalosDSMC06} and $216 \times 10^6$ particles in
  three dimensions~\cite{OppelstrupPRE2009}.  For many relevant
  biological systems, the number of molecules would less than 1000, or
  even in the single digits in many cases.  For example, the
  chemotaxis system in bacteria shows ``sensitivity to concentrations
  as low as 3 ligands per cell volume"~\cite{Sourjik2012}.  In
  mammalian cells, the number of molecules of a particular mRNA was on
  the order of tens or hundreds of molecules per
  cell~\cite{Raj2006PLoS}.  Ref.~\cite{Beck2011} measured the copy
  numbers of approximately 7300 proteins in a common human tissue
  culture cell line (U2OS), and found that proteins involved in
  signaling, cell communication, regulation of cellular processes,
  catalysis of post-translational modifications, and lipid metabolism
  tend to be present at copy numbers of less than 500 per cell.  }


\section{Applications} \label{S:Applications}

\comment{In Subsection~\ref{S:ComparisonPotentials} we
investigate the effects of drift due to several potentials on reaction
time and location statistics. In Subsection~\ref{S:DNAproteinApplication} 
we present a simplified model
of a protein-polymer system, in which two reacting molecules 
undergo drift-diffusion along a polymer, and may also unbind
from the polymer and diffuse in three dimensions.  We study the
interaction between polymer geometry, binding potentials along the 
polymer,  and unbinding rates.}

\subsection{Comparison of Potentials}  \label{S:ComparisonPotentials}

To demonstrate the contrasting effects that can be produced by 
different drifts, we consider the reaction $A + B \to \varnothing$ where 
the molecules diffuse within various potential energy landscapes.  
We consider the following three cases: (i) zero potential, $V(x)=0$; (ii) a 
one-well potential, $V(x)= \cos(2 \pi x)$; and (iii) a two-well potential, 
$V(x)= \cos(4 \pi x)$. We use a domain of length $L = 1$ unit with reflecting 
boundaries, and diffusion coefficient $D = 1$ unit$^2/$sec for both 
molecular species $A$ and $B$.

\begin{figure}
\centering
\includegraphics{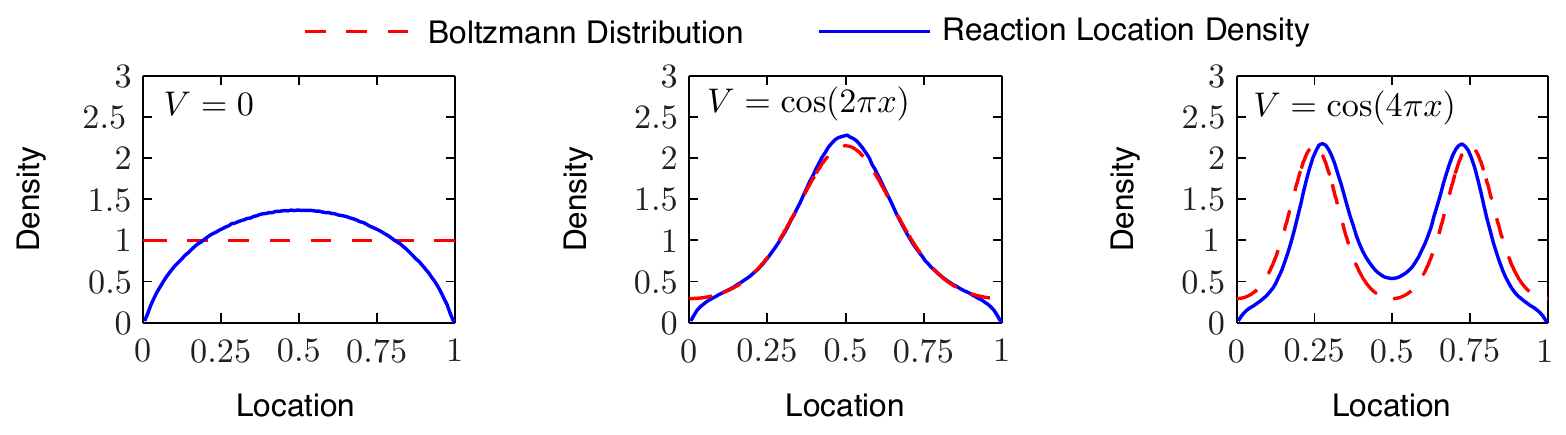}
\caption{\label{fig:RxnLocDensity} \small 
Reaction locations from $ A + B \to
  \varnothing $ DL-FPKMC simulations with  $M^A(0) = M^B(0)  = 1$, and 
  the Gibbs-Boltzmann distribution for each potential.  The reaction location is 
  $( Q^A(t) + Q^B(t) )/2$, where $t$ is the time of the reaction.
  Each graph of reaction locations is based on $10^7$ simulations.  
  $r_{\textrm{R}}= 0.02$. $h^{\textrm{max}}_s = h_p = r_{\textrm{R}}/8$.  
  $Q^A(0)$ and  $Q^B(0) \sim \mathcal{U}(0,1)$. The plotted densities 
  were determined by binning the reaction locations into $100$ bins. }
\end{figure}

In the absence of reactions, the equilibrium probability density for a
molecule to be at location $x \comment{ \in (0, L) }$ is given by the \comment{Gibbs-Boltzmann
distribution in Eq.~\eqref{eq:GibbsBoltzmann}}.  We compare the
Gibbs-Boltzmann distributions for each of the three potentials to the
reaction locations from the $A + B \to \varnothing$ DL-FPKMC
simulations in the particular case of two molecules, $M^A(0) = M^B(0)
= 1$.  The results are shown in Figure~\ref{fig:RxnLocDensity}.  The
potentials serve to spatially ``confine'' molecules, in the sense that
molecules are most likely to be found in locations where the potential
energy is lowest.  Consequently, the reactions are most likely to
occur in such low energy locations.

\begin{figure}
\centering
\includegraphics[width=16cm]{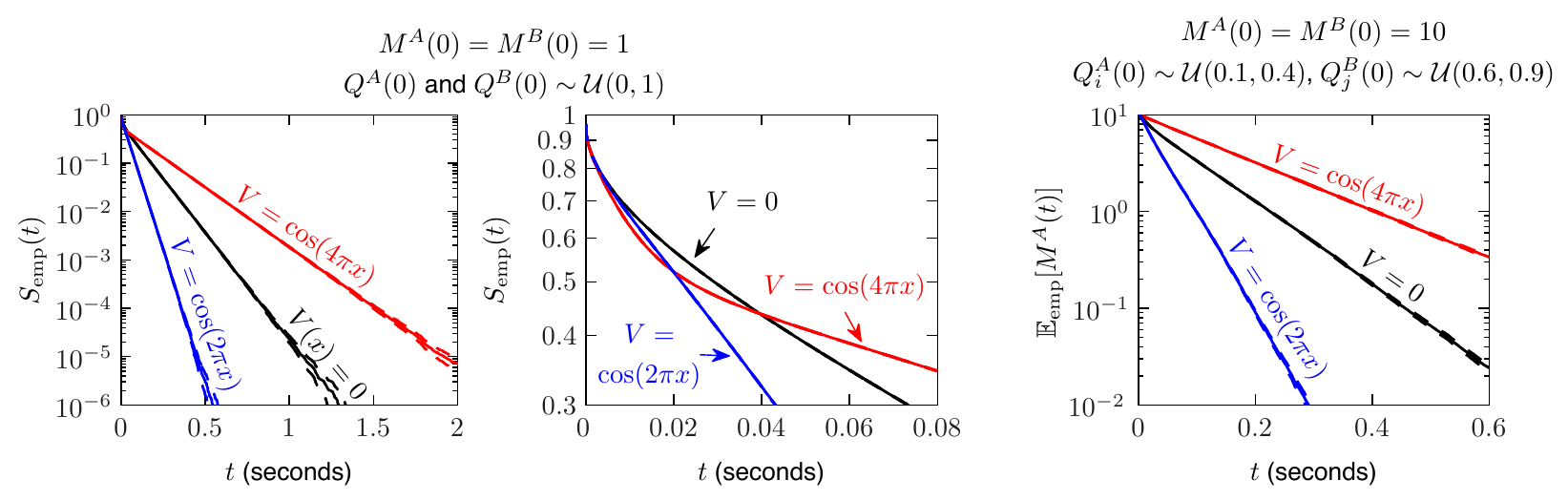}
\caption{\label{fig:RxnTimeCDFs} \small \emph{Left and center panels}:
  Empirical survival probability $S_{\textrm{emp}}(t)$ for the two-molecule 
  $A + B \to \varnothing$ reaction. Both panels show the same
  survival probabilities, however the axes have different scales.
  The dashed lines show $99\%$ confidence bounds based on $10^7$
  simulations. $r_{\textrm{R}} = 0.02$.  $h^{\textrm{max}}_s = h_p =
  r_{\textrm{R}}/8$.  \emph{Right panel}: Empirical mean number of
  molecules of $A$ remaining at time $t$, for the reaction $ A + B \to
  \varnothing$ with $M^A(0) = M^B(0) = 10$.  The dashed lines show
  $99\%$ confidence bounds based on $4 \times 10^4$ simulations.
  $r_{\textrm{R}} = 0.02$.  $h_p = r_{\textrm{R}}/16 $.
  $h^{\textrm{max}}_s = 8 h_p$.}
\end{figure}

As would be expected, the mean reaction time for the two-molecule $A +
B \to \varnothing$ reaction is faster with the one-well potential
($\sim 0.03481$ sec) than with no potential ($\sim 0.06483$ sec),
while slower with the two-well potential ($\sim 0.09887$ sec).  Figure
\ref{fig:RxnTimeCDFs} (left and center panels) compares the survival
probabilities $S_{\textrm{\textrm{emp}}}(t)$ for the three different
potentials.  The semi-log graphs of $S_{\textrm{\textrm{emp}}}(t)$ for
all three potentials appear linear except at short times, indicating
that the reaction time distributions could be described by 
exponential distributions for larger times
(Fig.~\ref{fig:RxnTimeCDFs}, left panel).  
For $V=0$, this is as to be expected since the survival probability 
distribution is known analytically and is given by the eigenfunction 
expansion~\eqref{eq:AnalyticSurvProb} in 
Appendix~\ref{S:AnalytlcAndNumericalPDEsolns}.
At short times
(Fig.~\ref{fig:RxnTimeCDFs}, center panel), the graphs  of
$S_{\textrm{\textrm{emp}}}(t)$ for the different potentials
are not linear and behave differently from each other.  Initially, reactions occur
more quickly with the two-well potential than with either the one-well
potential or no potential.  However, after $S_{\textrm{emp}}(t)$ has
decreased to about $50\%$, reactions occur more slowly with the
two-well potential than with the other potentials.  As would be
expected, if the initial locations of the two molecules are in the
same well of a potential, then they tend to react more quickly;
whereas, if the two molecules start in different energy wells, then
the time until they react tends to be longer.

Figure~\ref{fig:RxnTimeCDFs} (right panel) shows $\E_{\textrm{emp}}[M^A(t)]$,
the mean number of molecules  of $A$ remaining at time $t$, when 
$M^A(0) = M^B(0)  = 10$, $Q_{i}^A(0) \sim \mathcal{U}(0.1,0.4)$  and 
$Q_{j}^B(0) \sim \mathcal{U}(0.6,0.9)$. In this case, $\E_{\textrm{emp}}[M^A(t)]$ 
could be described by exponential distributions for all three potentials. 
In the case of $M^A(0) = M^B(0)  = 50$ with $Q^A_i(0)$ and  
$Q^B_j(0) \sim \mathcal{U}(0,1)$, $\E_{\textrm{emp}}[M^A(t)]$ is not plotted 
but behaves very similarly to $S_{\textrm{emp}}(t)$ in the two-molecule case
with $Q^A(0)$ and  $Q^B(0) \sim \mathcal{U}(0,1)$
(Fig.~\ref{fig:RxnTimeCDFs},  left and center panels).


\subsection{Role of Biopolymer Geometry  \comment{and Binding Potentials}
in Protein Diffusive Search} \label{S:DNAproteinApplication}

\begin{figure}
\centering
\includegraphics[width=2.8in]{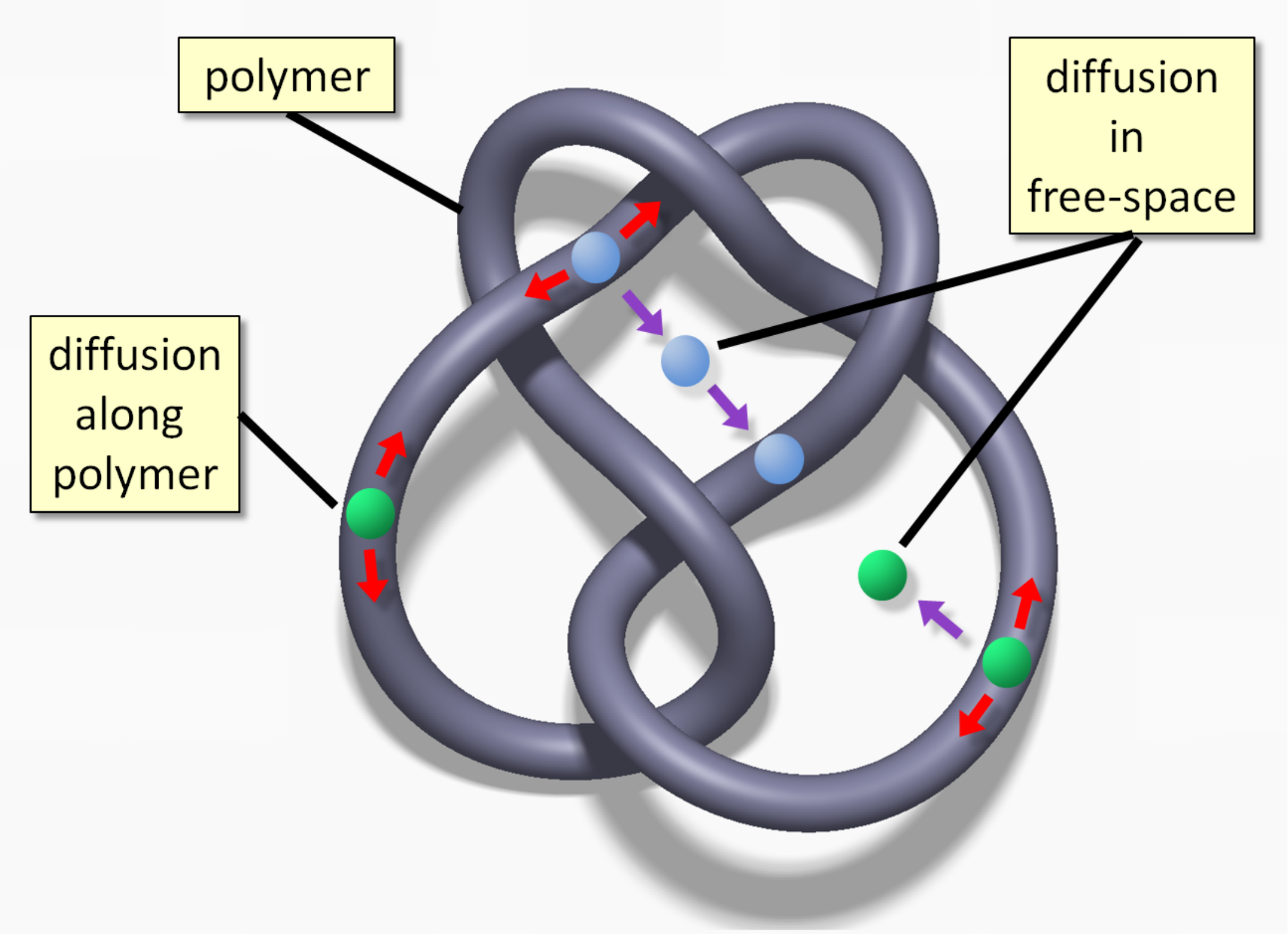}
\caption[Protein-DNA Schematic]{ \small Biopolymer diffusion-excursion search
  and protein-protein interactions.  The geometry of an individual biopolymer 
  or network is expected to be influential in the search kinetics
  of proteins which diffuse in one dimension along the biopolymers, but may
  also detach to undergo excursions in three dimensions.}
\label{fig_ProteinDNASchematic}
\end{figure}
Many proteins diffuse effectively in one dimension by adhering loosely
to biopolymer filaments, such as actin, microtubules, or
DNA~\cite{AlbertsMOLECELLBIO}.  In the case of DNA regulatory
proteins, single molecule experiments and theoretical work indicate
that observed rapid kinetic rates are achieved by a search process
involving a combination of one-dimensional diffusion of the protein
sliding along the DNA in conjunction with three-dimensional diffusive
excursions~\cite{Wang2006,Riggs1970, Shimamoto1999, Gorman2008,
  MirnyDNASlide04, MirnySlideGlob09, Hu2006, Berg1981}.  The geometry
of the individual biopolymer or network of biopolymers is expected to
play an influential role in such search kinetics.  We introduce a model
to investigate this phenomena in the diffusion of proteins that
interact with biopolymers having basic knotted closed-loop
configurations, see Figure~\ref{fig_ProteinDNASchematic}.  In our
model we also allow for position-dependent forces acting on the
proteins, such as might occur from heterogeneities in DNA sequence or
protein binding sites.
We perform simulations of two proteins undergoing a search process
until they encounter each other on a biopolymer. The process involves
one-dimensional drift-diffusion along the polymer in combination with
absorption/desorption events associated with excursions from the
polymer.  For example, this process could model the formation of a
regulatory complex at a non-specific DNA binding site.  Movement along
the polymer is simulated using the drift-diffusion DL-FPKMC method 
introduced in the preceding sections.  The effects of the biopolymer
geometry are taken into account through the absorption/desorption
statistics computed by \comment{numerically} solving a diffusion equation in three
dimensions for each specified biopolymer configuration.  We remark
that the model and methods we introduce are readily extended to more
complex geometries, polymer networks, and filament bundles.

\subsubsection{Model of the Biopolymer Drift-Diffusion Process and Three-Dimensional Excursions}
\label{S:ModelPolymerExcursions}
The biopolymer is represented geometrically as a one-dimensional
filament embedded within a three-dimensional space
$\Omega$.  The heterogeneity
along the biopolymer strand, for example changes in DNA sequence,
density of bound proteins, or other chemical factors, are taken into
account through the potential in the drift-diffusion
process.  To account for diffusive excursions in three dimensions, we
introduce desorption and absorption events for the protein from the
biopolymer.  Let $\lambda_{\textrm{off}}$ be the first-order rate at
which a protein desorbs from the biopolymer.  We assume the protein
must undergo a conformational change before regaining affinity to bind
to the biopolymer (otherwise it would instantly re-absorb at the
desorption position).  This conformational change is considered as a
first-order reaction with the rate $\mu_a$. To model this process, let
$p(s_2, t_2 | s_1, t_1)$ denote the probability density for a protein
that desorbs at location $s_1$ at time $t_1$ to re-absorb at location
$s_2$ at time $t_2$.

This probability density can be calculated by considering the two
stages of the protein desorption and re-absorption to the biopolymer.
In the first stage, the protein desorbs at location $s_1$ along the
biopolymer and performs three-dimensional diffusion until the
occurrence of the conformational change at a time $\tau_a$ after
desorption.  At this time the probability density the protein is at
the location $\mb{X}$ is given by $p_s (\mb{X} \, | s_1, \tau_a) =
\frac{1}{(2\pi Dt)^{3/2}} \exp\left\lbrack -\frac{(\mb{X} -
    \mb{Z}(s_1))^2}{4D\tau_a}\right\rbrack$, where $\mb{Z}(s_1)$
denotes the position on the biopolymer from which the protein desorbed
and $D$ the three-dimensional diffusion coefficient of the protein.
In the second stage, the protein has affinity for the biopolymer and
diffusively re-absorbs to the biopolymer at the location $s_2$ at time
$t_2$.  We denote this probability density by $p_a(s_2, t_2 | \mb{X},
t_1 + \tau_a)$.
This gives the density $p(s_2, t_2 | s_1, t_1) =$ 
$\int_\Omega \int_0^{t_2-t_1} p_a(s_2, 
t_2 | \mb{x}, t_1 + t) \, p_s (\mb{x} | s_1, t) \, \mu_a e^{-\mu_a t}
\, dt \, d\mb{x}$, \comment{which takes into account
the role played by the geometry of the biopolymer conformation.}
 
For simplicity we shall assume the biopolymer geometric conformation
is fixed, resulting in a desorption and absorption process that is
stationary.  This gives $p_a(s_2, t_2 | s_1, t_1 + \tau_a) = p_a(s_2,
t_2 - t_1 - \tau_a | s_1, 0)$.  In practice, we tabulate $p_a(s_2, t
\, | \, \mb{X}, 0)$ at select locations, $\mb{X}$, as a one-time
off-line calculation using a numerical diffusion equation solver.
This provides a very efficient method to simulate the excursions.  To
obtain $\tau_a$, we first generate a random exponential time at which
the affinity conformational change of the protein occurs.  We then
sample the protein location, $\mb{X} = \mb{x}$, using the density
$p_s(\mb{X} \, | s_1, \tau_a)$.  From our pre-constructed table, we
can then sample $p_a(s_2, t_2 - t_1 - \tau_a \, | \, \mb{x} )$ to find
a time of re-absorption, $t_2$, and the protein association location,
$s_2$.
The protein will eventually be re-absorbed, because
 the three-dimensional spatial domain has reflecting boundaries and re-absorption
 to the polymer is modeled using a sink term which has support on
 a set of positive measure (see Subsection~\ref{S:PolymerSink}).

The modeled search process of the two proteins proceeds by
drift-diffusion along the biopolymer \comment{and three-dimensional 
diffusive excursions} until the proteins encounter each other 
\comment{on the polymer}.  The drift-diffusion along the biopolymer is 
handled by our one-dimensional DL-FPKMC algorithm
with periodic boundary conditions.  Upon a desorption event from
location $s_1$,  the protein is repositioned on the biopolymer at
location $s_2$ at the re-attachment time $t_2$ according to the
density $p(s_2, t_2 | s_1, t_1)$.  
This repositioning process models the three-dimensional diffusive 
excursion of the protein until re-absorbing to biopolymer.  
We remark that reflecting boundary conditions \comment{or potential barriers} 
in our DL-FPKMC method could be used to model obstructions on the biopolymer.
 Either potential \comment{sinks} or absorbing Dirichlet boundary conditions 
 could be used to model irreversible binding sites on the polymer.

\subsubsection{Numerical Methods for Excursion-Time Probability Distribution and Reattachment Locations}
\label{S:PolymerSink}
We numerically solve the three-dimensional diffusion equation with
diffusivity $D$ for the position of a detached molecule at a given
time. Denote by $\rho(\vx,t)$ the corresponding probability density
that solves the diffusion equation.  Reattachment sites on a
biopolymer are modeled using sink terms in the equation. We let
$s$ label the parameterization variable along the biopolymer,
  and $s_0$ the detachment position.  From the probability mass
absorbed at the locations of the sinks, we can obtain the probability
densities $p(s,t|s_0,0)$, where $p(s, t | s_0, 0)$ is as
  defined in the previous subsection.
   In the diffusion equation we use a sink term of
the form
\begin{equation*} 
g(\vy,t) = - \paren{ \int \lambda(s)q(\vy-\vx(s)) \, ds } \rho (\vy,t),
\end{equation*}
where $\lambda$ is the intensity of the sink and $q$ is an averaging
kernel function. The rate of the reattachment flux is given by
\begin{equation*}
\frac{dQ}{dt}(s,t) = \int\frac{\lambda(s)q(\vy-\vx(s))}{\int \lambda(s')q(\vy-\vx(s'))ds'}g(\vy,t)\, d\vy 
= -\int\lambda(s)q(\vy-\vx(s))\rho(\vy,t) \, d\vy. 
\end{equation*} 
In the discretized form this becomes
\begin{equation*}
g(\vy_m,t) = -\sum_k \lambda_k q(\vy_m-\vx(s_k))\Delta s_k \rho(\vy_m,t)
\end{equation*}
and the rate of the flux becomes
\begin{equation*}
\frac{dQ_k(t)}{dt} = - \sum_m \lambda_k q(\vy_m-\vx_k)\Delta{s}_k \rho (\vy_m,t) \Delta{\vy}_m.
\end{equation*}

The particular choice that we make for the kernel $q$ is the
four-point Peskin-$\delta$ kernel~\cite{Peskin2002}. For $\vx =
(x,y,z)$ this is given by $q(\vx) =
\phi_p(x/\Delta{x})\phi_p(y/\Delta{x})\phi_p(z/\Delta{x}),$ where
\begin{eqnarray}
 \phi_p(u) = \frac{1}{16}\left\{
\begin{array}{lrl}
3-2u + \sqrt{1+4u-4u^2}  &   0 \leq &u \leq 1,\\
5-2u - \sqrt{-7+12u-4u^2} &   1 \leq &u \leq 2,\\
0                      &   2 \leq &u,
\end{array} \right. \nonumber 
\label{eq_PeskinDelta}
\end{eqnarray}
with $\phi_p(-u) = \phi_p(u)$. This particular function is chosen to
reduce numerical error induced by the off-lattice shifts of the polymer 
adsorption locations relative to the underlying discretization lattice 
used in the three-dimensional diffusion solver.  We take the absorption
parameter $\lambda$ to have two stages. Namely, a ``non-sticky" state
$\lambda = 0$ and a ``sticky" state $\lambda(s) = \lambda_k =
\lambda_\textrm{absorb}$. The latter is taken to be uniform along
the biopolymer. As discussed in Subsection~\ref{S:ModelPolymerExcursions}, 
the protein becomes sticky after undergoing a conformational change, the time 
for which is modeled as an exponentially distributed random variable.

\begin{figure}
  \centering
  \includegraphics[width=3.5in]{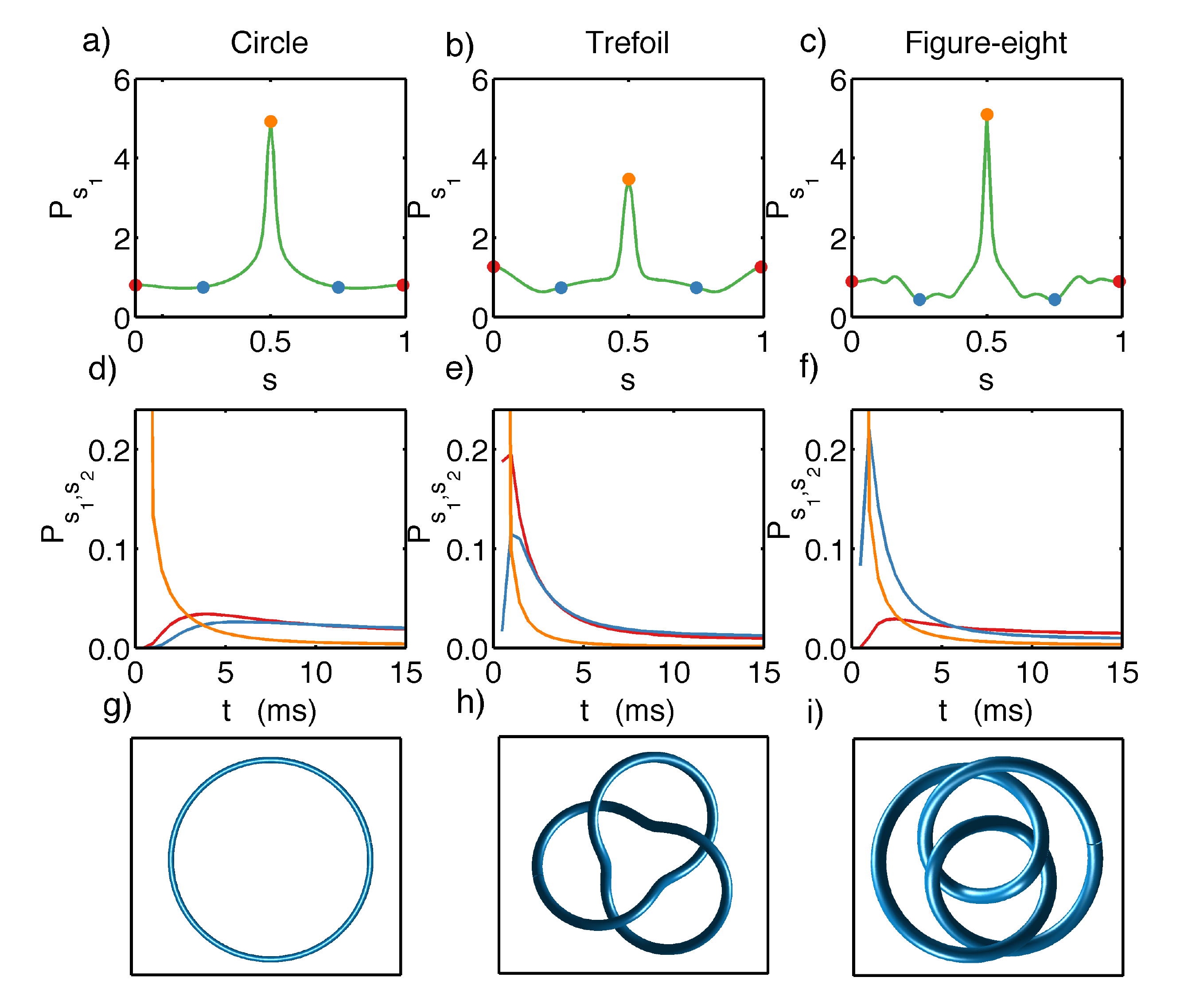}
  \caption[Protein-DNA Reattachment Probabilities]{ \small Reattachment
    probability distributions for the different topologies. (a--c)
    Probability density of absorbing at a given location
    $\mb{Z}(s)$ when the protein desorbed at location $\mb{Z}(0.5)$.
    $P_{s_1}(s) = \int_0^\infty p(s,t|s_1 = 0.5,0)\, dt$. (d--f) 
    Conditional probability densities for a protein that desorbed
    at location $\mb{Z}(0.5)$ to absorb at a point $s_2$.  For select
    locations these are indicated in the first row by a colored dot
    with the same color convention holding for the plots in the second
    row. $P_{s_1,s_2}(t) = p(s_2, t\, | \, s_1=0.5, 0)$.  (g--i) 
    Biopolymer geometries from left to right: a circular
    loop, trefoil knot, and figure-eight knot.}
  \label{fig_DNAAttachProb}
\end{figure}
We consider three distinct geometric configurations for the
biopolymer, each with the same arc-length but a different topology.
Biopolymers represented by a circular unknotted loop, a trefoil knot,
and a figure-eight knot are studied.  We use parameterizations given
by $X(s) = c*(x(s), y(s), z(s))$ for $0\leq s \leq 1$, where $c$ is
chosen so that the arc-length of each of the knots is
$L_{\textrm{biopolymer}} = 1000$nm.  The $x$, $y$, and $z$
parameterizations used for the three polymer configurations are:
unknotted loop, $x(s) = \cos(2\pi s)$, $y(s) = \sin(2\pi s)$, $z(s) =
0$; trefoil knot, $x(s) = (2+\cos(6\pi s))\cos(4 \pi s)$, $y(s) =
(2+\cos(6 \pi s))\sin(4\pi s)$, $z(s) = \sin(6\pi s)$; and
figure-eight knot, $x(s) = (2+\cos(4\pi s))\cos(6 \pi s)$, $y(s) =
(2+\cos(4 \pi s))\sin(6\pi s)$, $z(s) = \sin(8\pi s)$.  For each
configuration we use the diffusion equation solver to tabulate the
joint reattachment time and location distribution (see
Fig.~\ref{fig_DNAAttachProb}). \comment{In the solver, the overall
  spatial domain containing the polymer is taken to be a box with
  sides of length $500$nm.}

\subsubsection{Simulation Results: Diffusion-Excursion Search with Different Biopolymer Geometries}

\begin{table}
\centering
\caption{\label{tab:biopolyParams} \small Diffusion-excursion search parameters.}
{\renewcommand{\arraystretch}{1.1} 
\renewcommand{\tabcolsep}{0.2cm}
\begin{tabular}{ p{3.3cm}   p{4cm} |   p{0.5cm}   p{3.3cm}    p{2.8cm} } 
\hline
3D Parameters & Value  & & 1D Parameters & Value \\ 
\hline
3D diffusion coefficient  &  $ 2.183823 \, \mu \text{m}^2 \text{sec}^{-1}$ & &
1D diffusion coefficient  & $0.01 \, \mu \text{m}^2 \text{sec}^{-1}$ \\
$ \mu_a$ & $5.1728\times 10^7 \mbox{ns}^{-1}$ & &
$\lambda_{\textrm{off}}$ & $0.02$  to $200 \, \text{sec}^{-1}$ \\
$\Delta{x}$  & $12.5 \mbox{nm}$   & &
$r_{\textrm{R}}$ & $20 \,$nm\\
$\Delta{t}$ & $4292.9 \mbox{ns}$          &  & $h_s^{\textrm{max}}$ & $10 \,$nm \\
$N_{\textrm{biopolymer}}$         &  $100$        &  &  $h_p$ & $5 \,$nm \\
$N_x,N_y,N_z$         & $40$    & &  $r_{\textrm{pair}}$ & $ 4 \, r_{\textrm{R}} = 80 \, \text{nm} $  \\
\hline
\end{tabular}}
\vspace{.5cm}
\end{table}

\begin{figure}
\centering
\includegraphics{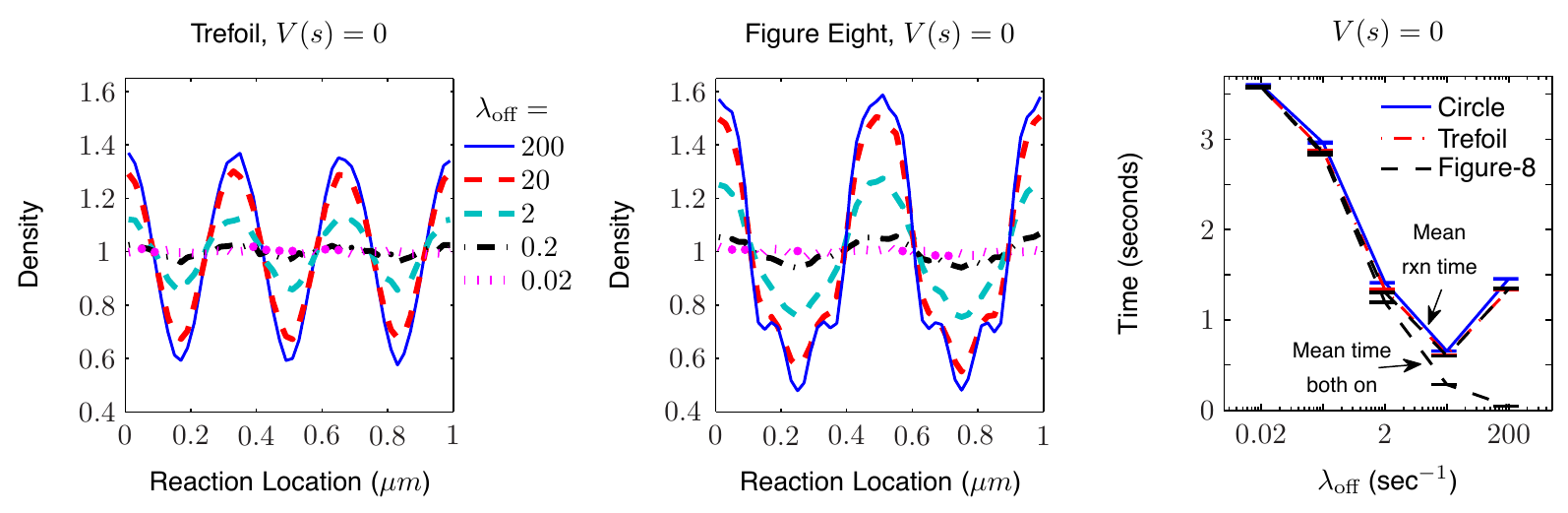}
\caption{\label{fig:VaryingOffRateV0} \small 
 \comment{\emph{Left and center panels:}}
  Reaction location densities for $ A + B \to \varnothing $ with  
  $M^A(0) = M^B(0)  = 1$ and $V(s)=0$.
  The reaction location is the midpoint between the locations of 
  $A$ and $B$ at the time of the reaction.
  Each graph is based on $10^6$ simulations.  
  The plotted densities were determined by binning the
  reaction locations into $50$ bins.
 Note that no potentials are imposed in this figure. The 
  non-uniformity in the reaction location density is a result of the
  polymer geometry.
  \comment{\emph{Right panel:}}
  Mean reaction time and mean time that both molecules are 
  on the biopolymer. Each data point is based on $10^6$ simulations.  
  Error bars indicate $99\%$ confidence intervals. \comment{The mean 
  time that both molecules are on the biopolymer is shown for the 
  figure-eight conformation only, but the results for the circle and trefoil
  conformations are similar.} 
  }
\end{figure}

\begin{figure}
\centering
\includegraphics{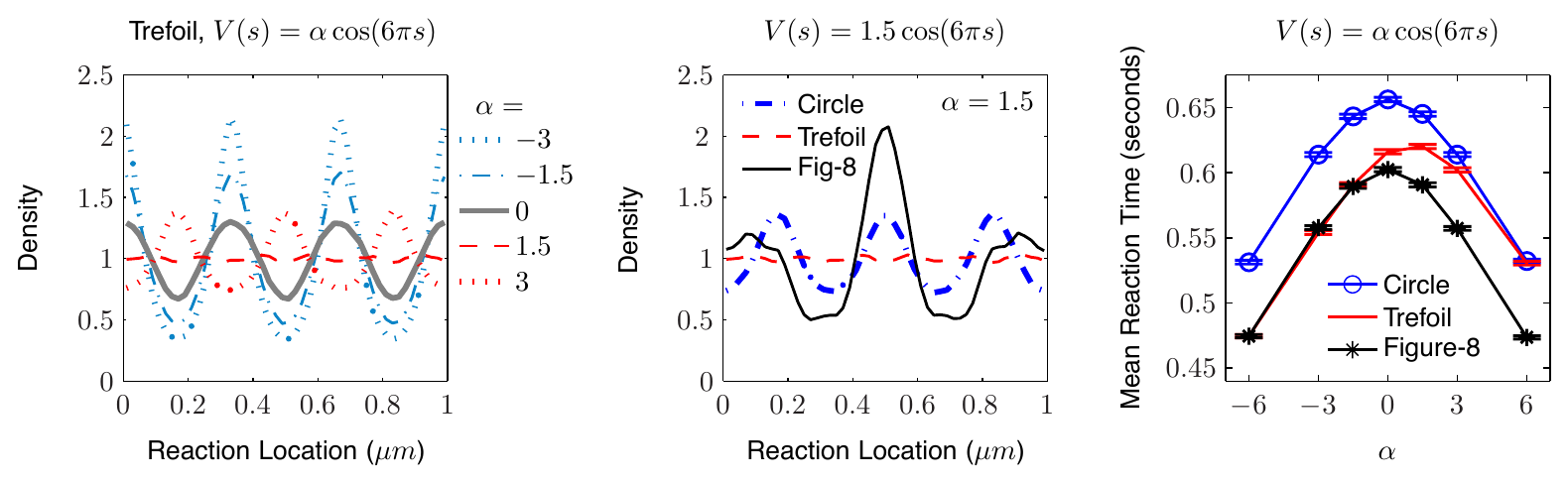}
\caption{\label{fig:vAcos6px} \small
 \comment{\emph{Left and center panels:}} 
  Reaction location densities for $ A + B \to \varnothing $, with  $M^A(0) = M^B(0)  = 1$,
  $V(s) = \alpha \cos{(6 \pi s)}$, \comment{and
  $\lambda_{\textrm{off}} = 20$ sec$^{-1}$}.
  Each graph is based on $10^6$ simulations.  
  The plotted densities were determined by binning the
  reaction locations into $50$ bins. 
  \comment{\emph{Right panel:}}
  Mean reaction times.  Each data point is based
  on $10^6$ simulations.  Error bars indicate $99\%$ confidence
  intervals. 
  }
\end{figure}

We assume one molecule each of protein species $A$ and $B$ are initially
present on the biopolymer. The initial positions are drawn from a
uniform distribution over the length of the biopolymer.  The two
molecules undergo drift-diffusion along the biopolymer and diffusion
in the three-dimensional space by the method developed in the
preceding subsections.  Detachment times are sampled from an
exponential distribution with rate $\lambda_{\textrm{off}}$. A
reaction between the two molecules can only occur when both are on the
biopolymer.

For the parameters we take \comment{the} three-dimensional diffusion coefficient
\comment{to be} $D = k_B{T}/\gamma$, where $\gamma = 6 \pi \eta a$ is the Stokes drag, 
\comment{$a = 10\,$nm is the radius of the protein}, $T = 298.15$ 
is room temperature, $k_B$ is Boltzmann constant, and the
dynamic fluid viscosity, $\eta$, is chosen to be 10 times the
viscosity of water.    The parameter for the exponential distribution
that represents when the protein becomes sticky is denoted by $\mu_a =
6D/ \comment{ r^2}$, \comment{where $r = L_\textrm{biopolymer} / 2\pi$}.
For convenience, the absorption strength, $\lambda_\textrm{absorb} =
0.005 \textrm{ns}^{-1}$, is chosen to be larger than the timescale
required for numerical stability in the diffusion equation solver.
The other parameters are specified in Table~\ref{tab:biopolyParams}.

\comment{
For each biopolymer conformation, the reaction
location distributions and the mean reaction times 
are investigated as the detachment
rate, $\lambda_{\textrm{off}}$, is varied
(Fig.~\ref{fig:VaryingOffRateV0}). These simulations are 
for the case that $V(s) = 0$ along the polymer. 
As $\lambda_{\textrm{off}}$ increases, the effects of the trefoil and
figure-eight biopolymer geometries on the reaction locations
become more pronounced
(Fig.~\ref{fig:VaryingOffRateV0}, left and center panels). For the
unknotted circle the reaction locations are approximately uniform over
the biopolymer regardless of the choice of $\lambda_{\textrm{off}}$.
The mean reaction times attain a minimum at 
approximately $\lambda_{\textrm{off}} = 20$ 
and are slightly faster with the trefoil and figure-eight biopolymer
conformations than with the circle conformation
(Fig.~\ref{fig:VaryingOffRateV0}, right panel).  The mean total time that both
molecules are on the biopolymer decreases monotonically as
$\lambda_{\textrm{off}}$ is increased 
(Fig.~\ref{fig:VaryingOffRateV0}, right panel).
Since the proteins spend more time diffusing in the
  three-dimensional space (off of the polymer) when $\lambda_{\textrm{off}}$ 
  is large,  the polymer geometries affect the re-absorption locations more strongly 
  (see Fig.~\ref{fig_DNAAttachProb}), which in turn affects the reaction
   locations (Fig.~\ref{fig:VaryingOffRateV0}).}

A natural question is how diffusion in a potential energy landscape
may enhance or counteract the effects of biopolymer geometry on the
reaction locations.  From Subsection~\ref{S:ComparisonPotentials}, we
expect the density of reaction locations to decrease in areas where
the potential is large and to increase in areas where the potential is 
small.  To test this idea, we investigated the trefoil knot
conformation with potentials of the form $V(s) = \alpha \cos{(6 \pi
  s)}$ for several values of $\alpha \in \brac{-6,6}$.
The resulting reaction location densities as $\alpha$ is varied are
shown in Figure~\ref{fig:vAcos6px} (left panel).  As expected,
the potential enhances the effects of the trefoil knot geometry when
$\alpha<0$, and counteracts the effects when $\alpha>0$.  When
$\alpha = 1.5$, the effects of the potential and the trefoil geometry
essentially cancel each other out. For comparison, we also ran
simulations with $\alpha = 1.5$ for each of the other two biopolymer
geometries (Fig.~\ref{fig:vAcos6px}, center panel).

When the potential $V(s) = \alpha \cos{(6 \pi s)}$ is used, 
the mean reaction time is fastest 
when the amplitude $| \alpha |$ is large (see 
Fig.~\ref{fig:vAcos6px}, right panel).
As discussed in
Subsection~\ref{S:ComparisonPotentials}, molecules that are in the
same energy well of a potential tend to react more quickly.  However,
if the proteins can only move on the biopolymer, potentials with large
amplitudes as used here would represent large energy barriers, greatly
slowing the time for the two molecules to find each other.  For
example, with $\alpha = -3$ and $\lambda_{\textrm{off}} = 0$, the mean
reaction time is approximately $80 \pm 9$ seconds, based on $10^3$
simulations (compare to Fig.~\ref{fig:vAcos6px}, right
panel).  The three-dimensional diffusion excursion provides an
alternative path for proteins to circumvent the energy barriers
present on the biopolymer.  This could be an important mechanism in
protein-protein interactions associated with biopolymers.
As shown in Figure~\ref{fig:vAcos6px} 
(right panel), for the circle and figure-eight conformations, the maximum 
mean reaction time occurs when $\alpha = 0$, and the mean reaction 
times with positive and negative $\alpha$ of the same magnitude are similar. 
In contrast, for the trefoil conformation, the maximum mean reaction time 
occurs when $\alpha = 1.5$, which is the same value of $\alpha$ that resulted
in the most uniform distribution of reaction locations (as was shown in
Fig.~\ref{fig:vAcos6px}, left panel). This demonstrates that 
the effect of the potential is not independent of the polymer geometry.

Overall, these results demonstrate the importance of incorporating
drift-diffusion and other spatial heterogeneities, such as biopolymer
geometry, when investigating biological chemical processes.  They 
illustrate the power of the presented DL-FPKMC 
method in capturing such
effects.  Our approach and methods can be readily extended to more
complex chemical kinetics, biopolymer geometries, or filament
networks.


\section{Conclusion}

We have presented a new Dynamic Lattice First-Passage
Kinetic Monte Carlo method (DL-FPKMC) that is capable of incorporating
the roles played by drift and spatial heterogeneities into the
stochastic dynamics of chemically reacting molecular species.  We have
shown that our numerical method is convergent for smooth potentials
with approximately second-order accuracy and for discontinuous
potentials with approximately first-order accuracy.  Unlike the
standard lattice RDME, our method retains bimolecular reactions as the
lattice spacings within protective domains are taken to
zero, and converges to the underlying SDLR model. In higher dimensions
we expect the use of Walk on Rectangles~\cite{LejayExactWOS2006}
techniques will facilitate the incorporation of complex geometries
into our method, while still allowing the use of basic Cartesian-grid
meshes within each protective domain (hence avoiding the need for
unstructured or embedded boundary meshes as commonly used for RDME
models.)

We have further demonstrated how our method can be utilized in
practice.  In particular, the many examples we investigated
demonstrate that drift can have a significant effect on reaction
locations, reaction time distributions, and number of reactions that
occur in a system.  To demonstrate how our approach might be used for
more complicated systems, we investigated the process of
diffusion-excursion of proteins interacting with a biopolymer, taking
into account geometric effects due to polymer shape.  Such processes
are thought to occur in the interactions between regulatory proteins
and DNA.  We considered a simplified model where the biopolymer has
basic knotted configurations.  Our results demonstrate that the
combined effects of the biopolymer-protein interactions, biopolymer
geometry, and drift-diffusion play a significant role in the chemical
kinetics.  The results also illustrate that capturing important
features of biological systems may require models that account for
drift-diffusion effects and spatial heterogeneities.  We expect the
DL-FPKMC methods we have introduced to be useful in performing more
realistic simulations of the chemical kinetics of biological systems.


\section{Acknowledgments}

A.J.M. was supported by NSF DMS-0602204 EMSW21-RTG and Boston
University.  S.A.I. and A.J.M. were supported by NSF DMS-0920886 and
NSF DMS-1255408.  P.J.A. and J.K.S. were supported by NSF CAREER -
0956210.  J.S. was supported by the CCS Summer Undergraduate Research
Fellowship (SURF) program at UCSB and from P.J.A.'s startup funds
(UCSB Math).  \comment{The authors would like to thank the reviewers
  for their helpful comments and suggestions. In particular,
  Section~\ref{S:RunningTime} and
  Appendix~\ref{S:ProtectiveDomainsAppendix} were added in response to
  reviewer comments.}


\appendix

\section{\comment{Constructing Protective Domains (PDs)}}   \label{S:ProtectiveDomainsAppendix} 

\comment{ This appendix describes our method for constructing
  protective domains (PDs). It should be noted that this is just one 
  of many possible approaches, and it is still an open question how 
  best to optimize the partitioning of space among PDs. For the 
  reader interested in this question we refer and defer
  to~\cite{DonevJCP2010}.  Regardless of the method used for
  constructing the PDs, our lattice approach of Section~\ref{S:oneDimMethod} 
  can be used to propagate molecules within their respective PDs.

  The PDs should be defined in a such a way that the distance
   between any two bimolecular reactants in separate PDs will remain strictly
   greater than their reaction radius, $r_\textrm{R}$, for as long as
   they remain in their respective PDs.
  This is necessary to ensure that the movements of molecules in
  separate PDs are independent.
  We allow the PDs of non-reacting molecules to overlap.

  The following is the approach for defining PDs that resulted in the
  scaling demonstrated in Section~\ref{S:RunningTime}.  Steps 3
  through 5 are specific to one dimension, but could be extended to
  higher dimensions by following a similar procedure in each
  coordinate. When constructing PDs for all molecules initially, or
  when updating the PDs of more than one molecule, we begin with 
  the molecule closest to the left endpoint of the overall domain and 
  then proceed to the right. An alternative approach would be to begin 
  with a seed molecule and proceed outward.

  \begin{enumerate}  [itemsep=0pt, topsep=1pt]

  \item For each molecule in need of a protective domain, 
    identify the nearest potential reaction
    partners. In one dimension, this can by done by keeping a list of
    all molecules ordered by location. In higher dimensions, the
    near-neighbor list (NNL) method described in~\cite{DonevJCP2010}
    can be used.

  \item Determine which molecules will be placed in pair PDs.  Two potential 
  reaction partners are placed in a pair if:
	
    \begin{itemize}  [itemsep=0pt, topsep=0pt]
      
    \item the two molecules are closer to each other than to any
      other potential reaction partners, and
      
    \item the distance between the two molecules is less than a pair threshold,
      $r_\textrm{pair}$, which is a parameter chosen by the user.	
      
    \end{itemize}
    
    For the reaction system $ A + B \to \varnothing $, we also enforce
    the condition in Eq.~\eqref{eq:PairConditions} when determining if
    two molecules will be placed in a pair. No other PDs are allowed
    to overlap with a pair PD.

  \item For each molecule, identify it's ``limiting neighbor," which we define in the
    following way:
    
    \begin{itemize}  [itemsep=0pt, topsep=0pt]
          
    \item For each molecule that will be in a pair PD,
      the limiting neighbor is the next nearest molecule of any type
      outside the pair.
      
    \item For a molecule that will be in a single PD,
      the limiting neighbor is the nearest molecule that is either a
       potential reaction partner or in a pair with another molecule.
      
    \item Let $d_\textrm{nbr}$ be the distance from a molecule to it's
      limiting neighbor.  If a PD for the limiting neighbor has
      already been defined, let $d_\textrm{nbrPD}$ be the distance
      from the molecule in question to the nearest endpoint of it's
      limiting neighbor's PD.
	
    \item Note, a molecule is not necessarily the limiting neighbor of
      it's limiting neighbor.
        
    \end{itemize}
		
  \item Determine the size of each pair PD and then the size of each 
    single PD as follows, ignoring overall domain boundaries until Step 5:

    \begin{itemize}  [itemsep=0pt, topsep=0pt]
	
    \item Pair PDs will be symmetric about the midpoint of the two 
      molecules' locations.  Single PDs will be symmetric about the location 
      of the molecule.
      
      \item  Let $r_\textrm{PD}$ for a pair PD be the 
      distance from either molecule to the nearest endpoint of the PD.
      For a single PD, $r_\textrm{PD}$ will denote the distance from the
      molecule to either endpoint of the PD.
         
    \item Define Condition 1 to be that the limiting neighbor is a potential
    reaction partner, and Condition 2 to be that a PD for the limiting neighbor
    has not yet been defined.
    
    For a single PD, calculate $r_\textrm{PD}$ by
          
      \begin{equation} \label{rPD}
	r_\textrm{PD}  = \left\{
          \begin{array}{cl} 
            (d_\textrm{nbr} -  r_\textrm{R}) / 2   
            & \text{ if Conditions 1 and 2 hold}   \\
            d_\textrm{nbrPD} -   r_\textrm{R}  & \text{ if only Condition 1 holds} \\
                        d_\textrm{nbr}/ 2   
            & \text{ if only Condition 2 holds}   \\
            d_\textrm{nbrPD} & \text{ otherwise.}     \end{array}  \right.
      \end{equation} 
      
    For a pair PD, calculate the quantity in Eq.~\eqref{rPD} for each molecule
    and then set  $r_\textrm{PD}$ for the pair to be the minimum of the two quantities.
      
  \item We recommend capping the size of pair PDs, so that two
    molecules will not remain in a pair if they have moved
    sufficiently far away from each other.  For example,
    $r_\textrm{PD}$ for a pair PD could be set to the minimum of the
    value calculated above and half the initial distance between the
    two molecules in the pair.\footnote{\comment{The size of single PDs may be
      capped also.  Generally speaking, the size of the PDs should be
      made as large as possible. However, making a particular
      molecule's PD as large as possible will result in less space for
      the PDs of neighbor molecules. Introducing a cap on the size of
      all PDs may result in a more equitable partitioning of space
      among them.
      We have not yet determined what value for the cap, if any,
      results in the most efficient performance of the overall
      algorithm. An optimal cap would mostly likely take into account
      the potential field (in DL-FPKMC), and would therefore vary
      spatially. See~\cite{DonevJCP2010} for discussion of the case
      where different molecules have different diffusion coefficients
		      in FPKMC.}}
                                
    \end{itemize}   
  \item If a PD as defined in Step 4 extends beyond an endpoint of the
    overall domain, truncate the PD so that one endpoint of the PD will coincide
    with the overall domain endpoint. Such PDs will no longer be
    symmetric. An alternative approach would be to treat the overall domain
    boundary as a neighbor in Steps 1 through 4, and only allow a molecule's PD to touch the
    domain boundary if the distance to the boundary is less than $r_\textrm{pair}$.
    
  \end{enumerate}
}

\section{Derivation of Non-Uniform Jump Rates}   \label{S:DerivationNonUnifJumpRates} 

To derive the non-uniform rates in
Eq.~\eqref{eq:NonUnifJumpRate}, we use the fluxes from
  the WPE discretization~\cite{ElstonPeskinJTB2003}. As illustrated in
  Figure~\ref{fig:NonUnifSubLattice} (top row, left), $x_1 < x_0 <
  x_2$ are the locations of mesh points with non-uniform spacing $h_j
  = |x_0 - x_j|$. The jump rates $a_{0j}$ from $x_0$ to $x_j$,
  $j=1,2$, are derived in this appendix.  The solution $\rho(x,t)$ of
  the Fokker-Planck equation~\eqref{eq:FPE}, gives the probability
  density of being at location $x$ at time $t$.  Let
  $\rho^{\textrm{eq}}(x)$ denote the equilibrium value of $\rho(x,t)$.
 Define $p_i(t)$ to be the probability of being at the mesh point
$x_i$ at time $t$ in the discrete master equation model.  We
consider the point $x_1$ to represent the interval
$(x_1-\frac{h_1}{2}, x_1+\frac{h_1}{2})$ in the sense that
\begin{equation} \label{eq:p1}
p_1(t) \approx  \int_{ x_1-\frac{h_1}{2} }^{ x_1+\frac{h_1}{2} }   
 \rho(x,t) \, dx \approx \rho(x_1, t) h_1.
\end{equation}
Similarly, $x_0$ represents $(x_0-\frac{h_1}{2}, x_0+\frac{h_2}{2})$
and $x_2$ represents $(x_2-\frac{h_2}{2}, x_2+\frac{h_2}{2})$, so
\begin{equation} \label{eq:p0andp2}
p_0(t) \approx   \rho(x_0, t) \frac{h_1+h_2}{2} \qquad \text{and} \qquad
p_2(t) \approx   \rho(x_2, t) h_2.
\end{equation}

Let $J(x,t)$ denote the flux
\begin{equation*}
J(x,t) = - D  \left(\rho(x,t) \frac{\partial V(x)}{\partial x}   + 
\frac{\partial \rho(x,t) }{\partial x} \right).
\end{equation*}
For convenience, we define
{\renewcommand{\arraystretch}{1.5}
\begin{equation*}
A_{ik}(t)  = \left\{
\begin{array}{cl} 
\frac{V(x_k) - V(x_i)} {\exp[V(x_k) - V(x_i) ] -1}  
 & \text{for }  V(x_k)\ne V(x_i) \\
1 & \text{otherwise. }  \end{array} \right. 
\end{equation*} } 
Based on the WPE discretization, we approximate the unidirectional 
outward flux from $x_0$ to $x_j$ by
\begin{equation*}
J_{0j}(t)  = \frac{D }{h_j} A_{0j}(t)  \rho(x_0,t),
\end{equation*}
and the unidirectional inward flux  from $x_j$ to $x_0$ by
\begin{equation*}
J_{j0}(t)  = \frac{D }{h_j} A_{j0}(t)  \rho(x_j,t).
\end{equation*}
Then, the net flux from $x_0$ to $x_j$ is $J_{0j}(t) - J_{j0}(t)$.
In the case that $x_j$, for $j=1$ or $2$, lies on an absorbing
Dirichlet boundary, $\rho(x_j, t) = 0$, and so we would have
$J_{j0} = 0$ throughout the following calculation. The
Fokker-Planck PDE~\eqref{eq:FPE} at $x_0$ is approximated by
\begin{align}
\frac{\partial \rho(x_0,t)}{\partial t} 
&=  - \frac{\partial }{\partial x}  J(x_0,t)
\approx \frac{2}{h_1+h_2} \int_{x_0-\frac{h_1}{2} }^{x_0+\frac{h_2}{2} } 
 - \frac{\partial }{\partial x} J(x,t)  \, dx  \notag \\
&=  \frac{2}{h_1+h_2} \left(  - J(x_0+\frac{h_2}{2},t) +J(x_0-\frac{h_1}{2},t) \right)  \label{eq:fluxDiff} \\
& \approx \frac{2}{h_1+h_2} \left( (J_{20} - J_{02}) - (J_{01} - J_{10}) \right) \notag \\
& =  \frac{2}{h_1+h_2} 
\left( \frac{D}{h_2} (A_{20} \, \rho(x_2,t) - A_{02} \, \rho(x_0,t)) -  
\frac{D}{h_1} (A_{01} \, \rho(x_0,t) - A_{10} \, \rho(x_1,t)) \right). \label{eq:WPEDensity}
 \end{align} 
 When $x_0$, $x_1$, and $x_2$ are interior points, Taylor series
 expansion shows the spatial truncation error to be $O(h_1-h_2)$ at
 $x_0$. For $h_1=h_2=h$ the discretization is second-order accurate
 and we recover the WPE discretization. Note, as described in
 Section~\ref{S:oneDimMethod}, we only use a non-uniform mesh for
 mesh cells bordering a boundary, or to move molecules onto a uniform
 mesh in a newly formed pair protective domain.  When one of $x_1$ or $x_2$
 corresponds to a (non-uniform) Dirichlet boundary point, and all
 interior mesh cells are uniform, the spatial discretization
 in Eq.~\eqref{eq:WPEDensity} in the case that $V(x) =0$ is known to 
 be second-order accurate for the Poisson equation, see Ref.~\cite{FedkiwJCP2002}.  
 For a reflecting boundary at $x_0-\frac{h_1}{2}$ or $x_0+\frac{h_2}{2}$, 
 the corresponding flux term in Eq.~\eqref{eq:fluxDiff} is zero (since we
 assume reflecting boundaries are at the edges of mesh cells).  The 
 corresponding terms involving $A_{0 j}$ and $A_{j 0}$ then drop out
 of Eq.~\eqref{eq:WPEDensity}.

 Multiplying \comment{Eq.~\eqref{eq:WPEDensity}} through by $\frac{h_1+h_2}{2}$ and making the
 substitutions in Eqs.~\eqref{eq:p1} and~\eqref{eq:p0andp2} yields the master
 equation
 \begin{equation*}
 \frac{d \, p_0(t)}{d t}  = 
 \left(\frac{D}{h_2^2} A_{20} \, p_2(t)  -  \frac{2D}{h_2(h_1+h_2)} A_{02} \, p_0(t)  \right) -  
\left( \frac{2D}{h_1(h_1+h_2)} A_{01} \, p_0(t) - \frac{D}{h_1^2} A_{10} \, p_1(t)  \right).
 \end{equation*}
 Thus, we obtain the non-uniform jump rates~\eqref{eq:NonUnifJumpRate}
 \begin{equation*}
a_{0j} = \frac{2 D}{h_j (h_1+h_2) } A_{0j}
\qquad j = 1,2.
\end{equation*}
 
Whenever $x_j$ is \emph{not} a boundary point, the jump
   rate $a_{j0}$ in the opposite direction agrees with the uniform
   rate in Eq.~\eqref{eq:JumpRate} with $h=h_j$.  Hence, the system still
   satisfies a discrete detailed balanced condition at equilibrium
   (similar to the uniform WPE discretization):
\begin{align*} 
a_{0j}\,  \rho^{\textrm{eq}}(x_0) \frac{h_1+h_2}{2} &=  a_{j0} \rho^{\textrm{eq}}(x_j) h_j ,
\end{align*}  
since $\rho^{\textrm{eq}}(x)$ is given by the Gibbs-Boltzmann distribution 
$ \rho^{\textrm{eq}}(x)  \propto  \exp[-V(x)] $ and
\begin{equation*}
  \frac{\rho^{\textrm{eq}}(x_j)}{\rho^{\textrm{eq}}(x_0)} = \frac{\exp[-V(x_j)] }{\exp[-V(x_0)] }  = \frac{1}{   \exp [V(x_j) -V(x_0)] } .
\end{equation*}


\section{Analytic and Numerical Solutions for the Two-Molecule 
Annihilation Reaction $A + B \to \varnothing$}   \label{S:AnalytlcAndNumericalPDEsolns} 

\begin{figure}
\centering
\includegraphics{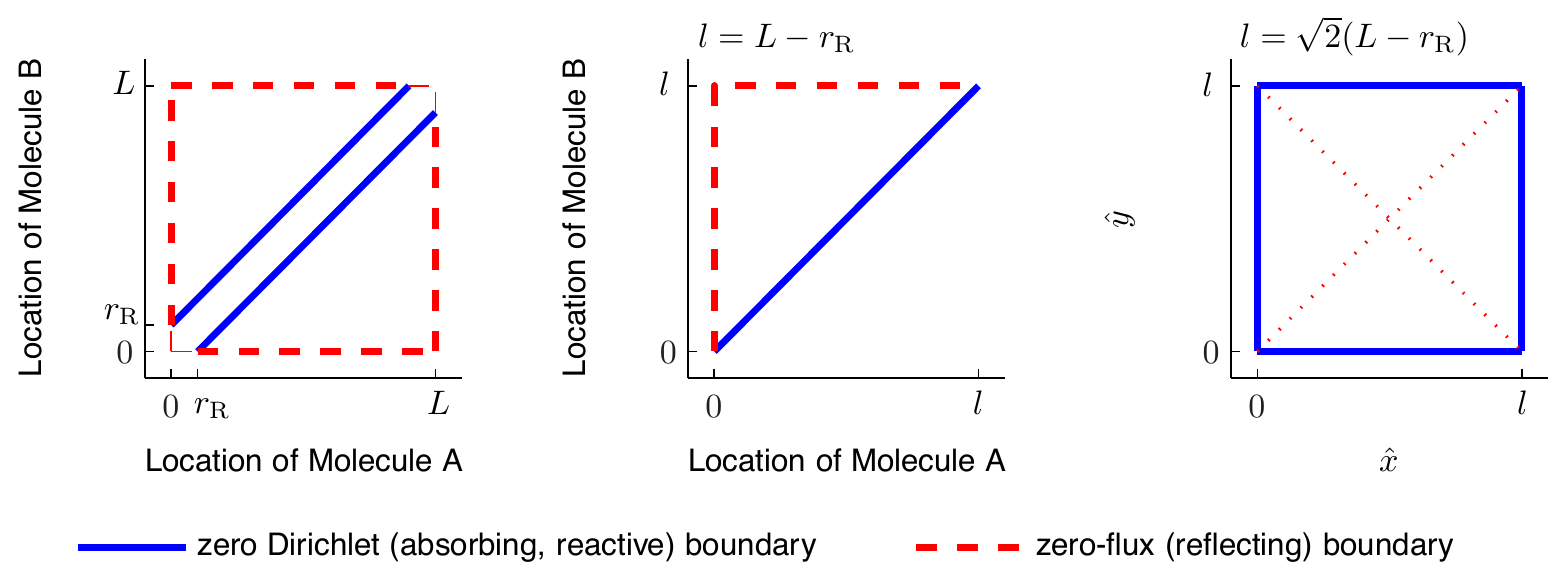}
\caption{\label{fig:PDEdomains} \small 
\emph{Left}: 2D domain equivalent to 1D simulation domain
of length $L$ in which the two molecules are located. 
The zero Dirichlet boundaries on the diagonals
 correspond to the reaction occurring between the two
molecules when they are one reaction radius apart. 
The zero-flux boundaries on the outer edges 
correspond to the reflecting endpoints of the 1D simulation domain.
\emph{Center}: PDE solver domain.
\emph{Right}: Eigenfunction expansion domain. }
\end{figure}

In this Appendix we discuss the numerical solution of the
Fokker-Planck equation and the analytic solution of the diffusion
equation to which we compared the DL-FPKMC simulation results in
Section~\ref{S:2ParticleConvResults}. The 1D reaction-drift-diffusion
system of two molecules undergoing the reaction $ A + B \to
\varnothing$ can be described by Eq.~\eqref{eq:FPEbvp} on the 2D
domain in the left panel of Figure~\ref{fig:PDEdomains}.  The
solutions of the Fokker-Planck equation on the two disjoint triangular
components of this domain are independent of each other.  By symmetry,
solving the Fokker-Planck equation on the domain in the left panel of
Figure~\ref{fig:PDEdomains} can be reduced to solving the same
equation on the single triangular domain in the center panel.  We have
written a PDE solver to solve the Fokker-Planck equation on the
triangular domain. The PDE solver uses the rates in
Eq.~\eqref{eq:JumpRate} from the WPE
discretization~\cite{ElstonPeskinJTB2003} of the Fokker-Planck
equation. This discretization is second-order accurate for smooth
potentials and first-order accurate for discontinuous potentials.  The
mesh in the PDE solver is uniform.

The solution to the diffusion equation on the triangular domain in
Figure~\ref{fig:PDEdomains} (center) can be recovered by solving on a
square domain formed by adjoining four copies of the triangular domain 
at their reflecting (zero Neumann) edges, as shown in the right panel of 
Figure~\ref{fig:PDEdomains}. On the square domain, an eigenfunction 
expansion for the solution can be determined analytically.  
This provides a check for the DL-FPKMC simulations in the $V=0$
case, as well as a check for the PDE solver in that case.  The
solution $\rho(x,y,t) $ to the diffusion equation with constant
initial condition $\rho_0$ and diffusion coefficient $D$ on a square
domain with sides of length $l$ is given by~\cite{HandbookPDE2002}:
\begin{align*}
\rho(x,y,t) = \frac{16 \rho_0}{\pi^2} & \left[ \sum_{n=0}^\infty \frac{1}{2n+1} 
\sin \left( \frac{(2n+1) \pi x}{l} \right) 
\exp \left( \frac{-\pi^2 (2n+1)^2 D t}{l^2}   \right) \right] 
\nonumber \\
\times & \left[ \sum_{m=0}^\infty \frac{1}{2m+1} 
\sin \left( \frac{(2m+1) \pi y}{l} \right) 
\exp \left( \frac{-\pi^2 (2m+1)^2 D t}{l^2}   \right) \right] .
\end{align*}
The survival probability $S(t)$, the probability
that the two particles have not yet reacted by time $t$, is 
\begin{align}  \label{eq:AnalyticSurvProb}
S(t) &=  \int_0^l \int_0^l \rho(x,y,t) \, dx \, dy
= \frac{64 \, \rho_0 \, l^2}{\pi^4}  \left[ \sum_{n=0}^\infty 
\left( \frac{1}{2n+1} \right)^2
\exp \left( \frac{-\pi^2 (2n+1)^2 D t}{l^2}   \right) \right] ^2 .
\end{align}
The mean reaction time is
\begin{align}   \label{eq:AnalyticMeanRxnTime}
\avg{T} = -\int_0^\infty t \, S'(t) \, dt  &= \int_0^\infty  S(t) \, dt 
= \frac{64 \,\rho_0 \, l^4}{D \, \pi^6}  \sum_{n=0}^\infty 
\sum_{m=0}^\infty
\left( \frac{1}{2n+1} \right)^2
\left( \frac{1}{2m+1} \right)^2
\frac{1}{(2n+1)^2 + (2m+1)^2} \notag  \\
&\comment{
= \frac{16 \,\rho_0 \, l^4}{D \, \pi^5}  \sum_{n=0}^\infty 
\left( \frac{1}{2n+1} \right)^4
 \left( \frac{\pi}{2} - \frac{\tanh(\frac{\pi}{2} (2n+1))}{(2n+1)} \right) }
\end{align}
Evaluating $\avg{T}$ at the parameter values corresponding 
to the two-molecule DL-FPKMC simulations 
\comment{in Section~\ref{S:2ParticleConvResults}}
gives a mean reaction time of  0.064831881311 seconds.

\comment{By symmetry the reaction location distribution can be
  calculated using the flux across any one of the four sides of the
  domain in Fig.~\ref{fig:PDEdomains} (right panel), since each side
  corresponds to a reactive boundary. Let $\veta = \veta(x,y)$ be the
  outward pointing normal at the point $(x,y) $. Using the boundary
  where $y=0$, the reaction location distribution is
\begin{align} \label{eq:AnalyticRxnLocCDF}
F^{\textrm{rxn}}(x) &=   4 \int_0^x \int_0^\infty  
 \left[ - \frac{\partial}{\partial \veta}  \rho(\xi,y,t)  \right]_{y=0} dt \, d\xi
 =  4  \int_0^x  \int_0^\infty  
 \left[ \frac{\partial}{\partial y}  \rho(\xi,y,t)  \right]_{y=0} dt  \, d\xi \\
 &=  \frac{64 \,\rho_0 \, l^2}{D \, \pi^4}  \sum_{n=0}^\infty 
\sum_{m=0}^\infty
\left( \frac{1}{2n+1} \right)^2
\frac{1}{(2n+1)^2 + (2m+1)^2} 
\left(1 - \cos \left( \frac{(2n+1) \pi x}{l} \right) \right)  \notag \\
 &= \frac{16 \,\rho_0 \, l^2}{D \, \pi^3}  \sum_{n=0}^\infty 
\left( \frac{1}{2n+1} \right)^3
\tanh \left(\frac{\pi (2n+1)}{2}  \right)
 \left(1 - \cos \left( \frac{(2n+1) \pi x}{l} \right) \right).     \notag
 \end{align}
 }

For each of the potential functions, we ran the PDE solver using the
Crank--Nicolson method in time, with spatial step sizes ranging from $\Delta
x=r_{\textrm{R}} = 0.02$ to $\Delta x=r_{\textrm{R}}/16$ and time
steps $\Delta t= \Delta x /16$.  
For $V=0$, we can check the results of the PDE solver against the 
analytic solution; in this case, the numerical mean reaction times
determined using the Crank--Nicolson method converge at approximately
second-order ($m \approx 2.0025$) to the
analytic mean reaction time.  For $V \ne 0$, there is no
analytic solution to which the numerical solutions can be compared;
however, the decrease in the pairwise differences between the mean
reaction times determined from the numerical PDE solutions at
successive step sizes indicates convergence
($m \approx 1.9968$ for $V_{\textrm{cos}}$ and $m \approx 0.9844$ for $V_{\textrm{step}}$).  
In the DL-FPKMC convergence results in Section~\ref{S:2ParticleConvResults}, the
empirical mean reaction times are compared to the analytic mean
reaction time in the case $V=0$, and to the numerical mean reaction
times determined using the Crank--Nicolson PDE solver with the finest
spatial step size, $\Delta x=r_{\textrm{R}}/16$, in the cases of
$V_{\textrm{cos}}$ and $V_{\textrm{step}}$.

The survival probabilities calculated using the Crank--Nicolson method
show small numerical oscillations during the first few time steps.
This is due to the incompatibility of the initial condition with the
zero Dirichlet boundary condition.  In order to numerically resolve
the survival probabilities more accurately at short times, we re-ran
the PDE solver using the Twizell--Gumel--Arigu (TGA)
method~\cite{TwizellL0Meth} from $t=0$ to $t=0.07$ seconds. The TGA
method is a second order, $L_0$ stable time discretization.  We also
used a finer time step, $\Delta t = (\Delta x)^2$ where 
$\Delta x=r_{\textrm{R}}/16$, to further improve the accuracy at short 
times when the survival probabilities change most rapidly.  We then 
determined the numerical survival probabilities using the results from
the TGA method for $t=0$ to $t=0.07$ seconds and using the results
from Crank--Nicolson method for $t>0.07$ seconds. At $t=0.07$, the 
absolute difference in survival probabilities between the two methods 
is less than $10^{-8}$ for $V_{\textrm{cos}}$ and less than $10^{-7}$ 
for $V_{\textrm{step}}$.  
 
 To check that using the TGA method 
 with a finer time step improved the 
 accuracy of the numerical survival probabilities
 over the Crank--Nicolson method, we compared
 the numerical survival probabilities for $V=0$ to the analytic
 survival probability  on the interval from $t=0$ to $t=0.07$.
 Table~\ref{tab:CN_TGA} shows the absolute errors between the numerical
 survival probability calculated using each method
 and the analytic survival probability.
 
 \vspace{.3cm}
\begin{table}
\centering
\caption{\label{tab:CN_TGA} \small Improved errors in numerical survival probabilities.}
 {\renewcommand{\arraystretch}{1.2} 
\renewcommand{\tabcolsep}{0.2cm}
\begin{tabular}{ p{2.2cm}  p{5cm}  p{4cm} } 
\hline  
& Crank--Nicolson & TGA  \\ 
& $\Delta x = r_{\textrm{R}}/1$6 , $\Delta t = \Delta x/16$ &
$\Delta x=r_{\textrm{R}}/16$ , $\Delta t = (\Delta x)^2$  \\ 
\hline  
$L^1$ Error &  4.1117e-6 &  4.7909e-8\\ 
$L^2$ Error &  8.2513e-5  & 1.6171e-6\\ 
$L^\infty$ Error& 5.8634e-3 & 2.2463e-4 \\
\hline
\end{tabular}}
\end{table}

\small{
\bibliographystyle{plain}

}
\end{document}